\documentclass{article}

\usepackage{amsfonts,amssymb}
\usepackage[section]{placeins}

\newtheorem{Theorem}{Theorem}[section] 
 
\newtheorem{Lemma}[Theorem]{Lemma}
\newtheorem{Corollary}[Theorem]{Corollary}
\newtheorem{Definition}[Theorem]{Definition}

\def\axioms{\begin{tabular*}{0.8\textwidth}{@{\extracolsep{\fill}}lr}}
\def\endaxioms{\end{tabular*}\smallskip}

\def\B{{\bf B}}
\def\T{{\bf T}}   % non-strict betweenness

\def\implies{\ \rightarrow\ }
\def\defined{\downarrow}
\def\Left{\mbox{\it Left\,}}
\def\Right{\mbox{\it Right\,}}
\def\Point{\mbox{\it Point\,}}
\def\MakePoint{\mbox{\it MakePoint\,}}

\def\Line{\mbox{\it Line\,}}
\def\Arc{\mbox{\it Arc\,}}
\def\Circle{\mbox{\it Circle\,}}

\def\Ray{\mbox{\it Ray\,}}
\def\SameOrder{\mbox{\it SameOrder\, }}

\def\Parallel{\mbox{\it Parallel\,}}
\def\CircleCenter{\mbox{\it center\,}}

\def\on{\mbox{\it on\,}}

\def\max{\mbox{\it max\,}}

\def\pointOnOne{\mbox{\it pointOn1\,}}
\def\pointOnTwo{\mbox{\it pointOn2\,}}

\def\SquareRoot{\mbox{\it SquareRoot\,}}

\def\IntersectLines{\mbox{\it IntersectLines\,}}
\def\IntersectLineCircleOne{\mbox{\it IntersectLineCircle1\,}}
\def\IntersectLineCircleTwo{\mbox{\it IntersectLineCircle2\,}}
\def\IntersectCirclesOne{\mbox{\it IntersectCircles1\,}}
\def\IntersectCirclesTwo{\mbox{\it IntersectCircles2\,}}
\def\IntersectCirclesSame{\mbox{\it IntersectCirclesSame\,}}
\def\IntersectCirclesOpp{\mbox{\it IntersectCirclesOpp\,}}
 
\def\Reflect{\mbox{\it Reflect\,}} 

\def\Rotate{\mbox{\it Rotate\,}}
\def\SameSide{\mbox{\it SameSide\,}}
\def\OppositeSide{\mbox{\it OppositeSide\,}}
\def\Extend{\mbox{\it Extend\,}}

\def\HilbertMultiply{\mbox{\it HilbertMultiply\,}}

\def\Reciprocal{\mbox{\it Reciprocal\,}}

\def\Add{\mbox{\it Add\,}}

\def\Midpoint{\mbox{\it Midpoint\,}}

\def\CircleThree{\mbox{\it Circle3\,}}
\def\project{\mbox{\it Project}}
\def\para{\mbox{\it Para}}
\def\Perp{\mbox{\it Perp}}

\def\R{{\mathbb R}}

\def\N{{\mathbb N}}
\def\K{{\mathbb K}}

\def\A{{\mathcal A}}

\def\Q{{\mathbb Q}}

\def\iff{{\leftrightarrow}}

\def\ECG{{\bf ECG}}

\def\Rrec{{\mathbb R_{\bf rec}}}
\def\F{{\mathbb F}}
\def\matrixtwo #1#2#3#4{\left(
\begin{array}{ll} 
#1 & #2\\ 
#3 & #4
\end{array} \right)}

\usepackage{graphicx}
\usepackage{pstricks}
\psset{unit=3cm}
\def\RotationFigure{%
\psset{unit=2.8cm}
\pspicture(0.15,-0.05)(2,1.8)
\qline(-0.2,0)(1.732,0)
\psdot(0,0)
\put(0,-0.1) {$o$}
\psdot(1,0)
\put(1,-0.1) {$a$}
\psdot(1.62,0)
\put(1.6,-0.1){$p$}
\qline(-0.1732,-0.1)(0.98,0.5623)  %oe
\qline(1.02,0.58)(1.5,0.866) %eL
\put(1.55,0.85) {$L$}
\qline(-0.1,-0.1732)(1,1.732)
\psdot(0.924,1.6)  %q
\put(0.8,1.6){$q$}
\qline(1,0)(1,0.5473)  %ae
\pscircle(1,0.5773){0.03}
\qline(1,0.6073)(1,1.15) 
\put(0.95,1.2){$K$}
\put(1.03,0.4773) {$e$}
\qline(0.978,0.597)(0.5,0.867) %ez
\psdot(0.5,0.867)  %z 
\put(0.4,0.867) {$z$}
\psline[linestyle=dashed](1,0)(0.5,0.867)
\psdot(0.75,0.4335)  %f
\put(0.68,0.3){$f$}
\endpspicture
\psset{unit=3cm}}

\def\AdditionFigure{%
\pspicture(3, 1.3)(0,-0.6)
\qline(0,0.1)(3,0.1)
\qline(0.1,0)(0.1,1.05)
\psdot(0.1,0.1)
\psdot(0.9,0.1)
\psdot(1.9,0.1)
\psdot(1.9,0.9)
\psdot(2.7,0.1)
\psdot(0.1,0.9)
\psline{->}(0.1,0.9)(1.9,0.9)
\qline(1.9,0.1)(1.9,1.05)
\psarc{->}(0.1,0.1){0.8}{0}{90}
\psarc{->}(1.9,0.1){0.8}{90}{0}
\put(0,-0.02) 0
\put(0.9,-0.02) {$A$}
\put(1.9,-0.02) {$B$}
\put(2.3,0.14) {$W=A+B$}
\put(0.14,0.77) {$U$}
\put(1.93, 0.77) {$V$}
\put(0.05,1.07) {$C$}
\put(1.85, 1.07) {$H$}
\put(3.05,0.05) {$R$}
\endpspicture}

\def\AdditionFigureTwo{%
\pspicture(3, 1.2)
\qline(0,0.9)(3,0.9)
\qline(0.9,0)(0.9,1.05)
\psline{->}(0.9,0.1)(2.8,0.1)
\qline(2.8,0.1)(2.8,1)
\psarc{->}(0.9,0.9){0.8}{180}{270}
\psarcn{->}(2.8,0.9){0.8}{270}{180}
\put(0,0.95) {$A$}
\put(0.92,0.95) 0
\put(1.9,0.95) {$W=A+B$}
\put(2.84,0.95){$B$}
\put(0.92,-0.01) {$U$}
\put(2.85,0.02) {$V$}
\endpspicture}

\def\CommutativeAdditionFigure{%
\pspicture(3, 2.2)
\qline(0,0.1)(3,0.1)
\qline(0.1,0)(0.1,2.03)
\psdot(0.1,0.1)
\psdot(0.9,0.1)
\psdot(1.9,0.1)
\psdot(1.9,0.9)
\psdot(2.7,0.1)
\psdot(0.1,0.9)
\psline{->}(0.1,0.9)(1.9,0.9)
\qline(1.9,0.1)(1.9,1.05)
\psarc{->}(0.1,0.1){0.8}{0}{90}
\psarcn{->}(1.9,0.1){0.8}{90}{0}
\put(-0.1,0) 0
\put(0.9,-0.04) {$A$}
\put(1.9,-0.04) {$B$}
\put(2.2,-0.04) {$W=A+B=B+A$}
\put(0.14,0.77) {$U$}
\put(1.93, 0.77) {$V$}
\put(0.05,2.07) {$C$}
\put(1.85, 1.07) {$H$}
\put(3.05,0.05) {$R$}
\psdot(0.1,1.9)
\put(0.14,1.77) {$X$}
\psdot(0.9,1.9)
\put(0.93,1.77) {$Y$}
\qline(0.9,0.1)(0.9,2.03)
\psarc{->}(0.1,0.1){1.8}{0}{90}
\psarcn{->}(0.9,0.1){1.8}{90}{0}
\psline{->}(0.1,1.9)(0.9,1.9)
\endpspicture}

\def\AssociativeAdditionFigure{%
\hskip-1.6cm
 \pspicture(3, 1.6)
\qline(0,0.05)(3,0.05)
\qline(0.05,0.05)(0.05,1.45)
\psdot(0.05,0.05)
\psdot(0.45,0.05)
\psdot(0.95,0.05)
\psdot(0.95,0.45)
\psdot(1.35,0.05)
\psdot(2.9,0.05)
\put(1.58,-0.08){$C$}
\qline(1.6,0.05)(1.6,1.45)
\psdot(0.05,0.45)
\psline{->}(0.05,0.45)(0.95,0.45)
\qline(0.95,0.05)(0.95,0.525)
\psarc{->}(0.05,0.05){0.4}{0}{90}
\psarcn{->}(0.95,0.05){0.4}{90}{0}
\put(0,-0.08) 0
\put(0.45,-0.08) {$A$}
\put(0.95,-0.08) {$B$}
\put(1.2,-0.08) {$A+B$}
\put(2.1,-0.08){$G=B+C$}
\put(2.8,-0.08){$A+B+C$}
\psdot(2.5,0.05)
\qline(2.5,0.05)(2.5,0.5)
\psdot(0.05,0.95)
\psdot(1.6,0.95)
\psdot(1.6,1.35)
\put(0.07,0.34) {$U$}
\put(0.965, 0.34) {$V$}
\psdot(0.05,1.35)
\put(0.08,1.24) {$X$}
\put(1.63,1.24) {$Y$}
\qline(1.35,0.05)(1.35,1.45)
\psarc{->}(0.05,0.05){1.3}{0}{90}
\psarc{->}(0.05,0.05){0.9}{0}{90}
\psarcn{->}(1.6,0.05){0.9}{90}{0}
\put(0.07,0.83){$Z$}
\psarcn{->}(1.6,0.05){1.3}{90}{0}
\psline{->}(0.05,1.35)(1.6,1.35)
\psline{->}(0.05,0.95)(1.6,0.95)
\psline{->}(0.95,0.45)(2.5,0.45)
\psdot(2.5,0.45)
\psarcn{->}(2.5,0.05){0.4}{90}{0}
\put(1.63,0.84){$E$}
\put(2.54,0.5){$P$}
\endpspicture}

\def\DescartesMultiplicationFigure{%
\pspicture(2, 1.2)
\qline(0,0)(2,0)
\put(0,-0.10) {$B$}
\qline(0,0)(1.8,1.08)
\qline(0.9,0)(0.8,0.48)
\put(0.9,-0.10) {$A$}
\put(0.7,0.48) {$C$}
\qline(1.35,0)(1.2,0.72)
\put(1.35,-0.10) {$D$}
\put(1.1,0.72) {$E$}
\endpspicture}

\def\HilbertMultiplicationFigure{%
\pspicture(2.8,2.2)
\qline(0,0.2)(2.5,0.2)
\qline(0.5,0)(0.5,2.0)
\pscircle(1.45,1.05){1.05}
\psdot(0.5,1.48)
\put(0.3,1.45) {$ab$}
\psdot(0.5,0.615)
\put(0.4,0.56) {$I$}
\psdot(0.8425,0.2)
\put(0.83,0.25){$a$}
\psdot(2.056,0.2)
\put(2.02,0.25){$b$}
\endpspicture}

\def\HilbertSquaringFigure{%
\pspicture(2.4,2)
\qline(0,0.2)(2.5,0.2)
\qline(0.5,0)(0.5,2.0)
\pscircle(1.2,1.0){0.8}
\psdot(0.5,1.38)
\put(0.32,1.42) {$a^2$}
\psdot(0.5,0.615)
\put(0.4,0.56) {$I$}
\psdot(1.2,0.2)
\put(1.15,0.27){$a$}
\endpspicture}

\def\ReciprocalScriptFigure{%
\pspicture(3,2.2)
\qline(0,0.2)(2.5,0.2)
\put(2.6,0.16){$L$}
\qline(0.5,0)(0.5,1.9)
\put(0.45,2.0){$K$}
\pscircle(1.5,1.05){1.00}
\psdot(1.5,1.05)
\put(1.4,1.1){$e$}
\psdot(0.5,1.05)
\put(0.4,1.0) {$I$}
\psdot(0.98,0.2)
\put(0.9,0.1){$a$}
\psdot(2.02,0.2)
\put(2.06,0.1){$b$}
\qline(0.5,1.05)(2.7,1.05)
\put(2.8, 1.07){$H$}
\qline(0.98,0.2)(0.26,1.475)
\put(0.2,1.52){$M$}
\psdot(0.74,0.625)
\put(0.8,0.55){$m$}
\qline(0.74,0.625)(1.88,1.2625)
\put(1.9,1.3){$J$}
\put(2.0,2.0){$C$}
\psline[linestyle=dashed](0.5,1.05)(1.26,1.475)
\put(1.27,1.5){$N$}
\endpspicture}
 
\def\TriangleCircumscriptionFigure{%
\pspicture(3.25,2.1)
\qline(0,0.2)(2.5,0.2)
\put(2.6,0.2){$L$}
\qline(1.45,0)(1.45,1.8)
\put(1.43,1.9){$K$}
\psdot(1.45,1.05)
\put(1.51,0.96){$e$}
\pscircle(1.45,1.05){1.05}
\psdot(0.5,0.615)
\put(0.4,0.56){$c$}
\psdot(0.8425,0.2)
\put(0.83,0.25){$a$}
\psdot(0.6712,0.4075)
\put(0.65,0.48){$p$}
\psline(0.8425,0.2)(0.5,0.615)
\psdot(2.056,0.2)
\put(2.02,0.25){$b$}
\psline(1.53,1.12)(0.19,0.0)
\put(1.6,1.15){$M$}
\psdot(1.45,0.2)
\put(1.52,0.25){$x$}
\endpspicture}

\def\EuclidExampleFigureOne{%
\psset{unit=2.8cm}
\pspicture(-0.4,-0.5)(1.8,1.3)
\psline(0,0.4)(1.48,0.989)   %pe
\qline(-0.2,0)(1.7,0)
\put(1.75,-0.03){$L$}
\psline(0,-0.1)(0,1.1)
\put(0, 1.15){$J$}
\pscircle(1.5,1){0.03}
\put(1.53,0.91){$e$}
\qline(1.5,-0.5)(1.5,0.97)
\qline(1.5,1.03)(1.5,1.1)
\put(1.48,1.15){$K$}
\psline(1.52,1.005)(1.6,1.035)
\put(1.64,1.05){$M$}
\psdot(0.16,0)
\put(0.16,-0.1){$a$}
\psdot(0,0.4)  %p
\put(-0.1,0.38){$p$}
\psline[linestyle=dashed](0,0.4)(1.5,0)  %px
\psline(0.16,0)(0,0.4)  % ap
\psdot(1.5,0)  %x
\put(1.54,-0.1){$x$}
\psdot(0,0)  %f
\put(-0.095,-0.12){$f$}
\psset{unit=3cm}
\endpspicture
}

\def\EuclidExampleFigureTwo{%
\psset{unit=2.8cm}
\pspicture(-0.4,-0.5)(1.8,1.3)
\pspolygon[fillstyle=solid,fillcolor=yellow](0,0.8)(0.75,0.2)(0,0.4)  %bmp
\pspolygon[fillstyle=solid,fillcolor=yellow](1.5,0)(0.75,0.2)(1.5,-0.4)  %xmc
\psline(0,0.4)(1.48,0.989)   %pe
\pspolygon[linecolor=red](0.43,0.57)(0,0.4)(1.5,0)(0,0.8)(0.75,0.2)(1.5,0) %qpxbmx
\psline(0.75,0.2)(1.5,-0.4)  %mc
\qline(-0.2,0)(1.7,0)
\put(1.75,-0.03){$L$}
\psline(0,-0.1)(0,1.1)
\put(0, 1.15){$J$}
\pscircle(1.5,1){0.03}
\put(1.53,0.91){$e$}
\qline(1.5,-0.5)(1.5,0.97)
\qline(1.5,1.03)(1.5,1.1)
\put(1.48,1.15){$K$}
\psline(1.52,1.005)(1.6,1.035)
\put(1.64,1.05){$M$}
\psdot(0.16,0)
\put(0.16,-0.1){$a$}
\psdot(0,0.4)  %p
\put(-0.1,0.38){$p$}
\psdot(0,0.8) %b
\put(-0.1,0.78){$b$}
\psdot(0.33,0.53)  %r 
\put(0.3,0.43){$r$}
\psline(0.16,0)(0,0.4)  % ap
\psdot(1.5,0)  %x
\put(1.54,-0.1){$x$}
\psdot(0.43,0.57)  %q
\put(0.41,0.64){$q$}
\psdot(0,0)  %f
\put(-0.095,-0.12){$f$}
\psdot(0.75,0.2) %m
\put(0.68,0.1){$m$}
\psdot(1.5,-0.4)%c
\put(1.54,-0.4){$c$}
\psset{unit=3cm}
\endpspicture
}

\def\JustifyMultiplicationFigure{%
\pspicture(3,2.2)
\qline(0,0.2)(2.5,0.2)
\qline(0.5,0)(0.5,2.0)
\pscircle(1.45,1.05){1.05}
\psdot(0.5,1.48)
\put(0.35,1.45) {$D$}
\psdot(0.5,0.615)
\put(0.4,0.56) {$I$}
\psdot(0.5,1.043)
\put(0.55,1.1){$M$}
\psdot(0.8425,0.2)
\put(0.83,0.1){$a$}
\psdot(2.056,0.2)
\put(2.07,0.1){$b$}
\psdot(1.425,0.2)
\put(1.41,0.1){$m$}
\psdot(1.425,1.043)
\put(1.47,1.03){$e$}
\psline[linestyle=dashed](0.5,1.043)(1.425,1.043)
\psline[linestyle=dashed](1.425,0.2)(1.425,1.043)
\psline(1.425,1.043)(2.056,0.2)
\psline(1.425,1.043)(0.8425,0.2)
\psline(1.425,1.043)(0.5,0.615)
\psline(1.425,1.043)(0.5,1.48)
\endpspicture}

\def\OtherIntersectionPointFigure{%
\pspicture(3,2.3)
\qline(0.5,0.2)(2.5,0.2)
\put(2.6,0.15){$L$}
\qline(1.45,-0.1)(1.45,2.1)
\psdot(1.45,1.05)
\put(1.5,1.05){$a$}
\pscircle(1.45,1.05){1.05}
\psdot(0.8425,0.2)
\put(0.8,0.3){$p$}
\psdot(2.056,0.2)
\put(2.03,0.3){$q$}
\psdot(1.45,0.2)
\put(1.5,0.3){$x$}
\put(2.4,1.7){$C$}
\psline[linestyle=dashed](2.056,0.2)(1.45,1.05)
\psline[linestyle=dashed](0.8425,0.2)(1.45,1.05)
\endpspicture}

\def\UniformReflectionFigure{%
\pspicture(3,2.3)
\qline(0,1)(2.5,1)
\put(2.6,0.95){$K$}
\put(1.48,0.5){$L$}
\qline(1.45,0.5)(1.45,2.1)
\psdot(1.45,1)
\put(1.48,0.87){$e$}
\psdot(0.6,1)
\put(0.6,0.87){$z$}
\psdot(2.3,1)
\put(2.3,0.87){$x$}
\psdot(1.45,1.85)
\put(1.5,1.93){$y$}
\psarc{->}(1.45,1){0.85}{0}{90}
\psarc{->}(1.45,1){0.85}{9}{180}
\psdot(1,1)
\put(1,0.87){$a$}
\psdot(1.45,1.45)
\put(1.5,1.5){$c$}
\psdot(1.9,1)
\put(1.9,0.87){$b$}
\psline[linestyle=dashed](1.9,1)(1.45,1.45)
\psline[linestyle=dashed](1,1)(1.45,1.45)
\psdot(1.225,1.225)
\psdot(1.675,1.225)
\psline[linestyle=dashed](1.45,1)(1.225,1.225)
\psline[linestyle=dashed](1.45,1)(1.675,1.225)
 \endpspicture
\vskip -1.5cm
 }

\def\SquareRootFigure{%
\pspicture(2.2, 1.2)
\qline(0.1,0.1)(2.1,0.1)
\psarc(1.1,0.1){1}{0}{180}
\qline(1.8071,0.1)(1.8071,0.8071)
\put(0,0) {$H$}
\put(1,0) {$K$}
\put(1.8,0) {$G$}
\put(2.1,0) {$F$}
\put(1.807,0.85) {$I$}
\endpspicture}

\def\HilbertDescartesFigure{%
\pspicture(-0.6,-0.25)(3.4,2)
\pscircle(0.3,0.8){0.9}
\qline(1.2,-0.2)(1.2,1.7)
\psline[linestyle=dashed](1.2,1.55)(-0.6,-0.25)
\psdot(1.2,0.8)
\put(1.3,1.15){$a$}
\psline[linestyle=dashed]{<->}(1.25,0.8)(1.25,1.55)
\psdot(0.785,1.55)
\psline[linestyle=dashed]{<->}(0.785,1.6)(1.2,1.6)
\put(0.95,1.66) 1
\qline(-0.6,1.55)(1.2,1.55)
\psdot(-0.19,1.55)
\psline[linestyle=dashed]{<->}(-0.19,1.8)(1.2,1.8)
\put (0.4,1.9){$a^2$}
\psarc(2.4,0.8){0.9}{0}{180}
\psarc[linestyle=dashed](2.4,0.8){0.9}{180}{360}
\psdot(2.89,1.55)
\qline(1.5,0.8)(3.3,0.8)
\qline(2.89,1.55)(2.89,0.8)
\psdot(2.89,0.8)
\psline[linestyle=dashed]{<->}(1.5,0.75)(2.89,0.75)
\psline[linestyle=dashed]{<->}(2.9,0.75)(3.3,0.75)
\put(2.3,0.6){$x=a^2$}
\put(3.1,0.6){1}
\put(2.66, 1.1){$\sqrt x$}
\endpspicture}

\def\EuclidParallelFigure{%
\pspicture(2.7, 1.2)
\qline(0.0,0.15)(2.47,0.15)
\qline(2.53,0.15)(2.7,0.15)
\qline(0.0,0.7)(2.7,0.7)
\qline(0.0,0.9)(2.473,0.163)
\qline(2.5275,0.14)(2.7,0.085)
\pscircle(2.5,0.15){0.03}
\psdot(0.67,0.7)
\put(0.67,0.78){$p$}
\psdot(1.29,0.515)
\put(1.28,0.39) {$a$}
\psdot(0.9,0.15)
\put(0.88,0.04) {$q$}
\psline[linestyle=dashed](0.67,0.7)(0.9,0.15)
\psdot(1.49,0.7)
\put(1.5,0.78) {$r$}
\qline(0.9,0.15)(1.49,0.7)
\put(-0.15,0.12) {$L$}
\put(-0.15,0.68) {$K$}
\put(-0.15,0.88) {$M$}
\endpspicture}

\def\EuclidParallelRawFigure{%
\pspicture(2.7, 1.2)
\pspolygon[fillstyle=solid,fillcolor=yellow](0.67,0.7)(0.78,0.425)(1.49,0.7)
\pspolygon[fillstyle=solid,fillcolor=yellow](0.07,0.15)(0.78,0.425)(0.9,0.15)
\qline(-0.15,0.15)(2.47,0.15)
\qline(2.53,0.15)(2.7,0.15)
\qline(-0.15,0.7)(2.7,0.7)
\qline(0.0,0.9)(2.473,0.163)
\qline(2.5275,0.14)(2.7,0.085)
\pscircle(2.5,0.15){0.03}
\psdot(0.67,0.7)
\put(0.67,0.78){$p$}
\psdot(1.29,0.515)
\put(1.28,0.39) {$a$}
\psdot(0.9,0.15)
\put(0.88,0.04) {$q$}
\psdot(0.07,0.15)
\put(0.07,0.04){$s$}
\qline(0.67,0.7)(0.9,0.15)
\psdot(1.49,0.7)
\put(1.5,0.78) {$r$}
\qline(1.49,0.7)(0.07,0.15)  %rs
\psdot(0.78,0.425)
\put(0.65,0.3){$t$}
\qline(0.9,0.15)(1.49,0.7)  %qr
\put(-0.3,0.12) {$L$}
\put(-0.3,0.68) {$K$}
\put(-0.15,0.88) {$M$}
\endpspicture}

\def\StrongParallelRawFigure{%
\pspicture(2.6, 1.2)
\pspolygon[fillstyle=solid,fillcolor=yellow](0.67,0.7)(0.78,0.425)(1.49,0.7)
\pspolygon[fillstyle=solid,fillcolor=yellow](0.07,0.15)(0.78,0.425)(0.9,0.15)
\qline(-0.15,0.15)(2.47,0.15)
\qline(2.53,0.15)(2.7,0.15)
\qline(-0.15,0.7)(2.7,0.7)
\qline(0.0,0.9)(2.473,0.163)
\qline(2.5275,0.14)(2.7,0.085)
\pscircle(2.5,0.15){0.03}
\psdot(0.67,0.7)
\put(0.67,0.78){$p$}
\psdot(1.7,0.393)
\put(1.73,0.42) {$a$}
\psdot(0.9,0.15)
\put(0.88,0.04) {$q$}
\psdot(0.07,0.15)
\put(0.07,0.04){$s$}
\qline(0.67,0.7)(0.9,0.15)
\psdot(1.49,0.7)
\put(1.5,0.78) {$r$}
\qline(1.49,0.7)(0.07,0.15)  %rs
\psdot(0.78,0.425)
\put(0.65,0.3){$t$}
\put(-0.3,0.12) {$L$}
\put(-0.3,0.68) {$K$}
\put(-0.15,0.88) {$M$}
\endpspicture}

\def\AlternateInteriorAnglesFigure{%
\pspolygon[fillstyle=solid,fillcolor=lightgray](0.7,0.7)(1.3,0.7)(0.8,0.425)
\pspolygon[fillstyle=solid,fillcolor=lightgray](0.9,0.15)(0.3,0.15)(0.8,0.425)
\pspicture(2, 0.9)
\qline(0.0,0.15)(2,0.15)
\qline(0.0,0.7)(2,0.7)
\psdot(0.7,0.7)
\put(0.67,0.78){$p$}
\psdot(0.9,0.15)
\put(0.88,0.04) {$q$}
\psdot(0.3,0.15)
\put(0.3,0.04){$s$}
\psdot(0.8,0.425)
\put(0.85,0.36){$t$}
\qline(0.3,0.15)(1.3,0.7)
\qline(0.7,0.7)(0.9,0.15)  %pq
\psdot(1.3,0.7)
\put(1.3,0.78) {$r$}
\put(-0.15,0.12) {$L$}
\put(-0.15,0.68) {$K$}
\endpspicture}

\def\AlternateInteriorAnglesSevenFigure{%
\pspicture(2, 0.9)
\pspolygon[fillstyle=solid,fillcolor=lightgray](0.7,0.7)(1.3,0.7)(0.75,0.55625)
\pspolygon[fillstyle=solid,fillcolor=lightgray](0.9,0.15)(0.3,0.15)(0.85,0.2875)
\qline(0.0,0.15)(2,0.15)
\qline(0.0,0.7)(2,0.7)
\psdot(0.7,0.7)
\put(0.67,0.78){$p$}
\psdot(0.9,0.15)
\put(0.88,0.04) {$q$}
\psdot(0.3,0.15)
\put(0.3,0.04){$s$}
\psdot(0.8,0.425)
\put(0.85,0.36){$t$}
\psdot(0.75,0.55625)
\put(0.64,0.53){$x$}
\psdot(0.85,0.2875)
\put(0.9,0.26){$y$}
\psline[linestyle=dashed](0.3,0.15)(1.3,0.7)
\qline(0.7,0.7)(0.9,0.15)  %pq
\psdot(1.3,0.7)
\put(1.3,0.78) {$r$}
\put(-0.15,0.12) {$L$}
\put(-0.15,0.68) {$K$}
\qline(0.3,0.15)(0.85,0.2875) %sy
\qline(0.75,0.55625)(1.3,0.7)  %xr
\endpspicture}

\def\PlayfairFigure{%
\pspicture(0,0.15)(2.2, 1.05)
\qline(0.0,0.15)(2.0,0.15)
\qline(0.0,0.7)(2.0,0.7)
\qline(0.0,0.9)(2.0,0.5)
\psdot(1,0.7)
\put(1,0.78){$p$}
\put(-0.15,0.12) {$L$}
\put(-0.15,0.68) {$K$}
\put(-0.15,0.88) {$M$}
\endpspicture}

\def\PlaneSeparationFigure{%
\pspicture(2.2,1)(0,-0.3)
\qline(0.0,0.15)(1.22,0.15)
\qline(1.27,0.15)(1.7,0.15)
\put(1.75,0.1){$L$}
\qline(0.9,-0.3)(0.9,1.0)
\qline(0.1,0.15)(1.7,0.75)
\qline(0.9,0.75)(1.7,0.75)
\qline(0.9,-0.3)(1.215,0.125)
\pscircle[linewidth=1pt](1.237,0.15){0.03}
\qline(1.2565,0.175)(1.7,0.75)
\psline[linestyle=dashed](0.9,0.6)(1.7,0.75)
\psdot(0.9,0.60)
\put(0.78,0.6){$a$}
\psdot(1.7,0.75)
\put(1.75,0.72){$b$}
\psdot(0.9,0.75)
\put(0.78,0.75){$w$}
\psdot(0.1,0.15)
\put(0.0,0.05){$r$}
\psdot(0.9,0.15)
\put(0.78,0.05){$x$}
\psdot(0.9,0.45)
\put(0.95,0.38){$m$}
\psdot(0.9,-0.3)
\put(0.78,-0.28){$c$}
\put(1.28,0.05){$z$}
\endpspicture}

\def\VerifyingPaschFigure{%
\pspicture(0,-0.1)(2,1.4)
\psset{unit=2cm}
\psdot(0,0)
\put(-0.1,-0.16){$a=(0,0)$}
\psdot(1,0)
\put(0.95,-0.14){$p$}
\psdot(2,0)
\put(2,-0.14){$c$}
\qline(0,0)(2,0)
\psdot(0.5,2)
\put(0.5,2.07){$b$}
\qline(0.5,2)(2,0)
\psdot(1.25,1)
\put(1.32,1.0){$q$}
\qline(0,0)(0.805,0.64)
\qline(0.85,0.68)(1.25,1)
\qline(1,0)(0.835,0.635)
\qline(0.815,0.69)(0.5,2)
\pscircle(0.825,0.66){0.04}
\put(0.68,0.66){$z$}
\psset{unit=3cm}
\endpspicture}

\def\InnerOuterPaschFigure{%
\pspicture(0,0.2)(3.5,1.55)
\psdot(0,0.6)
\put(0,0.47){$a$}
\pscircle(1,0.6){0.03}
\put(1,0.47){$x$}
\qline(0,0.6)(0.97,0.6)
\psdot(1.3,0.6)
\put(1.35,0.55){$q$}
\qline(1.03,0.6)(1.3,0.6)
\psdot(1.1,1.4)
\put(1.1,1.45){$c$}
\qline(0,0.6)(1.1,1.4)  % ac
\qline(1.1,1.4)(1.375,0.3) % cb
\psdot(1.375,0.3)
\put(1.43,0.25){$b$}
\psdot(0.522,0.977)
\put(0.43,1.03){$p$}
\qline(1.375,0.3)(1.015,0.584)  %bx
\qline(0.52,0.98)(0.979,0.62)  %px

\psdot(2,0.6)
\put(2,0.47){$a$}
\psdot(3,0.6)
\put(3,0.47){$p$}
\qline(2,0.6)(3.3,0.6)   %ac
\psdot(3.3,0.6)
\put(3.35,0.55){$c$}
\psdot(3.1,1.4)
\put(3.1,1.45){$q$}
\qline(2,0.6)(2.5,0.96) % ax
\qline(3.1,1.4)(2.543,0.987)  % qx
\qline(3.1,1.4)(3.375,0.3) % qb
\psdot(3.375,0.3)
\put(3.43,0.25){$b$}
\pscircle(2.522,0.977){0.03}
\put(2.43,1.03){$x$}
\qline(3.375,0.3)(2.543,0.96)   %bx
\endpspicture}

\def\MakePointFigure{%
\pspicture(3,2.2)
\qline(0,0.2)(2.7,0.2)
\put(2.8,0.17){$L$}
\put(2.8,1.67){$K$}
\qline(0.5,0)(0.5,2.0)
\psdot(0.5,1.7)
\put(0.41,1.73){$z$}
\psdot(0.5,0.2)
\put(0.41,0.07){$0$}
\psdot(2.0,0.2)
\put(2.0,0.07){$y$}
\psdot(1.5,0.2)
\put(1.55,0.07){$x$}
\psarc{->}(0.5,0.2){1.5}{0}{90}
\qline(1.5,0)(1.5,1.7)
\qline(1.5,1.7)(1.5,2.0)
\put(1.48,2.04){$U$}
\qline(0,1.7)(1.5,1.7)
\psdot(1.5,1.7) 
\put(1.55,1.6){$p$}
\qline(1.5,1.7)(2.7,1.7)
\endpspicture}

\def\PerpFigure{%
\pspicture(-0.5,-1.3)(2,2.5)
\qline(-0.5,0.6)(1.75,0.6)
\put(1.8,0.57){$L$}
\psdot(0.625,1.0)
\put(0.635,1.05){$x$}
\psdot(0.25,0.6)
\put(0.205,0.48){$a$}
\psdot(0.5,0.6)
\put(0.5,0.48){$b$}
\pscircle[linestyle=dashed](0.5,0.6){0.545}
\put(0.19,1.13){$Q$}
\psdot(1.043,0.6)
\put(1.06,0.48){$c$}
\pscircle(0.625,1.0){1.0185}
\put(1.7,1.2){$C$}
\psdot(-0.305,0.6)
\put(-0.38,0.48){$p$}
\psdot(1.556,0.6)
\put(1.556,0.48){$q$}
\psarc(-0.305,0.6){1.861}{270}{90}
\put(-0.47,2.4){$K$}
\psarc(1.556,0.6){1.861}{90}{270}
\put(1.65,2.4){$R$}
\psdot(0.625,2.21)
\put(0.75,2.2){$e$}
\psdot(0.625,-1.01)
\put(0.75,-1.04){$d$}
\qline(0.625,2.3)(0.625,-1.1)
\endpspicture}

\def\CenterImpliesPlayfairFigure{%
\pspicture(0.3,-0.5)(3,1.3)
\qline(0.4,0)(2.97,0)
\put(0.2,-0.05){$L$}
\qline(0.4,1)(2.9,1)
\put(2.93,0.95){$K$}
\psdot(1,1)
\put(0.89,0.89){$p$}
\pscircle(3,0){0.03}
\put(3,-0.1){$e$}
\psline(0.75,1.125)(2.98,0.015)
\put(0.6,1.125){$M$}
\psdot(2,0.5)
\put(2.05,0.53){$r$}
\psdot(1,0)
\put(1.05,-0.1){$w$}
\qline(1,-0.5)(1,1.1)
\put(0.95,1.15){$J$}
\qline(1,0.5)(1.403,1.305)
\psdot(1.403,1.305)
\put(1.3,1.305){$t$}
\psdot(1,0.5)
\put(0.84,0.5){$m$}
\psdot(1.2015,0.905)
\put(1.27,0.9){$x$}
\psdot(1,-0.5)
\put(1.05,-0.42){$z$}
\psline[linestyle=dashed](1,-0.5)(2.97,-0.01)
\psline[linestyle=dashed](1,0.5)(2.97,0.01)
\psline[linestyle=dashed](1.403,1.305)(2.97,0.03)
\endpspicture}

\def\CenterImpliesStrongParallelFigure{%
\pspicture(0,-0.5)(3,1.7)
\qline(0.22,0)(2.97,0)
\put(0.1,-0.05){$L$}
\qline(0.24,1)(2.9,1)
\put(0.1,0.95){$K$}
\psdot(1,1)
\put(1.02,0.89){$p$}
\pscircle(3,0){0.03}
\put(3,-0.1){$e$}
\psline(0.75,1.125)(2.98,0.015)
\put(0.6,1.125){$M$}
\psdot(2,0.5)
\put(2.05,0.53){$r$}
\psdot(1,0)
\put(1,-0.1){$w$}
\qline(1,0)(1,1)
\qline(0.4,-0.2)(1.225,1.45)
\put(0.45,-0.2){$J$}
\psdot(0.775,0.55)
\put(0.64,0.5){$x$}
\psdot(1.225,1.45)
\put(1.14,1.45){$y$}
\psdot(1,0.5)
\put(1.03,0.57){$m$}
\pscircle[linestyle=dotted](1,1){0.5} 
\psdot(0.78,-0.55)
\put(0.66,-0.55){$z$}
\qline(0.78,-0.55)(0.78,0.55)
\psline[linestyle=dashed](0.78,-0.55)(2.97,-0.01)
\psline[linestyle=dashed](0.78,0.55)(2.97,0.01)
\psline[linestyle=dashed](1.225,1.45)(2.97,0.03)
\endpspicture}

\def\RotationHelperFigure{%
\pspicture(-2,-0.87)(2,0.65)
\qline(-1.6,0)(1.6,0)  %J
\put(1.65,-0.04){$J$}
\psdot(1.3,0)
\put(1.33,-0.08){$m$}
\psdot(-1.3,0)
\put(-1.33,0.05){$n$}
\psdot(0,0)
\put(-0.04,0.05){$b$}
\psdot(0.9,0)
\put(0.93,-0.08){$p$}
\pscircle(0.9,0.45){0.03}
\put(0.88,0.53){$e$}
\qline(-1.6,-0.8)(0.88,0.44) %qe
\qline(0.92,0.46)(1.6,0.8) %er
\put(1.65,0.77){$L$}
\psdot(1.3,0.65)
\put(1.28,0.71){$r$}
\psdot(-1.3,-0.65)
\put(-1.35,-0.75){$q$}
\qline(0.9,-0.85)(0.9,0.42) % K
\put(0.85,-0.95){$K$}
\qline(1.3,0.65)(1.3,-0.85) % M 
\put(1.25,-0.95){$M$}
\qline(-1.6,-0.65)(1.6,-0.65)   % R
\put(1.65,-0.69){$R$}
\psdot(0.9,-0.65)  %g
\put(0.93,-0.75){$g$}
\psdot(1.3,-0.65) %f
\put(1.33,-0.75){$f$}
\qline(-1.3,-0.65)(-1.3,0)  %qn
\endpspicture}

\def\EuclidPerpFigure{%
\pspicture(-1,-0.1)(2,1.4)
\pscircle(-1,0){0.03}  %e
\put(-1,-0.09){$e$}
\qline(-0.97,0)(1.7,0) %L
\put(1.75,-0.03){$L$}
\psdot(1,0)   %f
\put(1,-0.1){$f$}
\psdot(1.5,0)  %a
\put(1.5,-0.09){$a$}
\qline(1,0)(1,1)  %pf
\psdot(1,1)   %p
\put(0.97,1.07){$p$}
\qline(1,1)(1.5,0) %pa
\qline(-0.98,0.021)(1.7,1.35)  %M
\put(1.75,1.32){$M$}
\endpspicture}

\def\EuclidRawPerpFigure{%
\pspicture(-1,-0.1)(2,1.4)
\pspolygon[fillstyle=solid,fillcolor=yellow](0.5,1)(1,1)(1,0.5)  %cpm
\pspolygon[fillstyle=solid,fillcolor=yellow](1.5,0)(1,0)(1,0.5)  %afm
\psline(1.5,0)(1.5,1)  % qa
\pscircle(-1,0){0.03}  %e
\put(-1,-0.09){$e$}
\qline(-0.97,0)(1.7,0) %L
\put(1.75,-0.03){$L$}
\psdot(1.5,1)  %q
\put(1.47,1.07){$q$}
\psdot(1,1)   %p
\put(0.97,1.07){$p$}
\qline(1,0)(1.5,1) %fq
\pspolygon[linecolor=red](1.25,0.5)(1,0)(1,1)(1.5,0)(1.1667,0.3333)(1,0.5)(1,1)  %xfpamp
\psdot(1,0.5)   %m
\put(1.03,0.53){$m$}
\psdot(1,0)   %f
\put(1,-0.1){$f$}
\psdot(1.5,0)  %a
\put(1.5,-0.09){$a$}
\qline(-0.98,0.021)(1.7,1.35)  %M
\put(1.75,1.32){$M$}
\psdot(0.5,1.0)  %a
\psdot(0.66,0.84)  %z
\put(0.63,0.72){$z$}
\psdot(0.5,1)    %c
\put(0.47,1.07){$c$}
\psdot(1.25,0.5)   %x
\put(1.3,0.5){$x$}
\psdot(1.1667,0.3333)   %b
\put(1.21,0.325){$b$}
\qline(0.5,1)(1.5,1) %cq

\endpspicture}

\def\EuclidPerpFigureTwo{%
\pspicture(-1,-0.1)(2,1.5)
\pspolygon[fillstyle=solid,fillcolor=yellow](0.5,1)(1,1)(1,0.5)  %cpm
\pspolygon[fillstyle=solid,fillcolor=yellow](1.5,0)(1,0)(1,0.5)  %afm
\pscircle(-1,0){0.03}  %e
\put(-1,-0.09){$e$}
\qline(-0.97,0)(1.7,0) %L
\put(1.75,-0.03){$L$}
\psdot(1,0)   %f
\put(1,-0.1){$f$}
\psdot(1.5,0)  %a
\put(1.5,-0.09){$a$}
\qline(1,0)(1,1)  %pf
\psdot(1,1)   %p
\put(0.97,1.07){$p$}
\psdot(1,0.5)   %m
\put(1.03,0.53){$m$}
\qline(1,1)(1.5,0) %pa
\qline(-0.98,0.021)(1.7,1.35)  %M
\put(1.75,1.32){$M$}
\psdot(0.5,1.0)  %a
\psline[linecolor=red](1.5,0)(0.5,1)   %ac
\psdot(0.66,0.84)  %z
\put(0.63,0.73){$z$}
\psdot(0.5,1)    %c
\put(0.47,1.07){$c$}
\endpspicture}

\def\OneSixteenFigure{%
\pspicture(0,0)(2,0.75)
\psline[linestyle=dashed, linecolor=red](1,0.4)(1.8,0)  %ad
\psline(0,0)(1.8,0)  % bd
\qline(0,0)(0.8,0.8)  %ba
\qline(0.8,0.8)(1.2,0)  %ac
\psdot(0.8,0.8)   %a
\put(0.76,0.68){$a$}
\psdot(1.2,0)  %c 
\put(1.17,-0.09){$c$}
\psdot(1.8,0)  %d
\put(1.75,-0.11){$d$}
\psline[linecolor=lightgray](1.2,0)(2 ,0.8) %cf
\psline[linecolor=lightgray](0,0)(2,0.8)  %bf
\psdot(1,0.4)  %e
\put(0.95,0.3){$e$}
\psdot(1.4,0.2)  %x
\put(1.37,0.24){$x$} 
\psdot(2,0.8)
\put(2,0.67){$f$}
\psdot(0,0)  %b
\put(-0.04,-0.11){$b$}
\endpspicture}

\def\MiddleParallelFigure{%
\pspicture(-0.2,-0.1)(2.3,1)
\qline(-0.2,0)(2.1,0)   %L
\put(2.13,-0.03){$L$}
\qline(0,0)(0,0.8)   %J
\put(-0.03,0.85){$J$}
\qline(1,0)(1,0.47)   %K, lower half
\qline(1,0.53)(1,0.8)   % K, upper half
\put(0.97,0.85){$K$}
\pscircle(1,0.5){0.03}  % r
\qline(2,0)(2,0.8)   %N
\put(1.97,0.85){$N$}
\qline(-0.2, 0.3)(0.97,0.49)  %M, left half
\qline(1.03,0.5)(2.1,0.66)  %M, right half
\put(2.13,0.64){$M$}
\psdot(0,0)  %j 
\put(-0.03,-0.12){$j$}
\psdot(1,0)  %k
\put(0.97,-0.12){$k$}
\psdot(2,0) %n 
\put(1.97,-0.12){$n$}
\psdot(0,0.33)
\put(0.04,0.26){$e$}
\psdot(2,0.645)
\put(1.91,0.54){$f$}
\endpspicture}

\begin{document}

\title{Constructive Geometry and the Parallel Postulate} 
\author{Michael Beeson}        % Enter your name between curly braces

%\address{Department of Computer Science, San Jos\/e State University,\\
%San Jos\'e, CA 95192, USA \\
%E-mail: beesonpublic@gmail.com\\
%www.MichaelBeeson.com/Research}
\maketitle 

\begin{abstract}
 Euclidean geometry, as presented by Euclid,  consists of straightedge-and-compass constructions 
and rigorous reasoning about the results of those constructions.  We show that Euclidean geometry 
can be developed using only intuitionistic logic. This involves finding ``uniform'' constructions where
normally a case distinction is used.  For example, in finding a perpendicular to line $L$ through point 
$p$, one usually uses two different constructions, ``erecting'' a perpendicular when $p$ is on $L$,
and ``dropping'' a perpendicular when $p$ is not on $L$, but in constructive geometry, it must
be done without a case distinction.  Classically, the models of 
Euclidean (straightedge-and-compass) geometry are planes over Euclidean fields.  We prove a similar theorem for constructive
Euclidean geometry, by showing how to define addition and multiplication without a case distinction
about the sign of the arguments.  With intuitionistic logic,
there are two possible definitions of Euclidean fields, which turn out to correspond to different
versions of the parallel postulate.

 We consider three versions of Euclid's parallel postulate.  The two most important are  
Euclid's own formulation in his Postulate 5, which says that under certain conditions two lines meet,
and Playfair's axiom (dating from 1795), which says there cannot be 
two distinct parallels to line $L$ through the same point $p$.  These differ in that Euclid 5
makes an existence assertion, while Playfair's axiom does not.  The third variant, which
we call the {\em strong parallel postulate},  isolates the existence assertion from the 
geometry:  it amounts to Playfair's axiom plus the principle that two distinct lines
that are not parallel do intersect.
The first main result of this paper is that Euclid 5 suffices to define coordinates,
addition, multiplication, and square roots geometrically.  

We completely settle the
questions about implications between the three versions of the parallel postulate.
The strong parallel postulate easily 
implies Euclid 5,  and Euclid 5 also implies the strong parallel postulate, as a 
corollary of coordinatization and definability of arithmetic.
We show that Playfair does not imply Euclid 5, and we also give some other independence results.
Our independence proofs are given  without discussing the exact choice of the other axioms of geometry;  
all we need is that one can interpret the geometric axioms in Euclidean field theory.   The 
independence proofs use 
Kripke models of Euclidean field theories based on carefully constructed rings of real-valued functions.  
``Field elements'' in these models are real-valued functions. 

\end{abstract}
\vfil\eject

\tableofcontents

\section{Introduction } 
\subsection{The purposes of this paper}
Euclid's geometry, written down about 300 BCE, has been extraordinarily influential in the development of mathematics, and 
prior to the twentieth century was regarded as a paradigmatic example of pure reasoning.%
\footnote{Readers interested in the historical context of Euclid are recommended to read \cite{dehn}, where Max Dehn puts forward
the hypothesis that Euclid's rigor was a reaction to the first ``foundational crisis'', the Pythagorean discovery of 
the irrationality of $\sqrt 2$.}
In this paper, we re-examine Euclidean geometry from the viewpoint of constructive mathematics.  
The phrase ``constructive geometry''
suggests, on the one hand,  that ``constructive'' refers to geometrical constructions with straightedge and compass.  
 On the other hand, the word ``constructive'' suggests the use of intuitionistic logic.%
 \footnote{ In intuitionistic logic, one is not allowed to prove that something exists
 by assuming it does not exist and reasoning to a contradiction; instead, one has to show 
 how to construct it.  That in turn leads to rejecting some proofs by cases,  unless we have 
 a way to decide which case holds.  These are considerations that apply to all mathematical reasoning;
 in this paper, we apply them to geometrical reasoning in particular. In this paper, ``constructive''
 and ``intuitionistic'' are synonymous; any differences between those terms are not relevant here.}
  We investigate the connections between these two meanings of the word.   
 Our method is to focus on the body of mathematics in Euclid's {\em Elements}, and to examine
what in Euclid is constructive, in the sense of ``constructive mathematics''.   
Our first aim was to formulate a suitable formal theory that would be faithful to both the ideas of Euclid and 
the constructive approach of 
Errett Bishop.  We presented a first version of such a theory in \cite{beeson-kobe}, based 
on Hilbert's axioms, and a second version in \cite{beeson2015b}, based on Tarski's axioms.
These axiomatizations helped us to see that there is a coherent body of informal mathematics
deserving the name ``elementary constructive geometry,''  or \ECG.  The first purpose of 
this paper is to reveal constructive geometry and show both what it has in common with 
traditional geometry, and what separates it from that subject.%
\footnote{The phrase ``classical
geometry'' is synonymous with ``traditional geometry''; it means
geometry with ordinary, non-constructive logic allowed.  Incidentally, that geometry is also 
classical in the sense of being old. }

The second purpose of this paper is to demonstrate that constructive geometry is 
sufficient to develop the propositions of Euclid's {\em Elements},
and more than that, it is possible to connect it, just as in classical geometry, to the 
theory of Euclidean fields.   A Euclidean field is an ordered field $\F$ in which non-negative
elements have square  roots.  Classically, the models of elementary Euclidean geometry are 
planes $\F^2$ over a Euclidean field $\F$.  We show that something similar is true for 
\ECG.   This requires an axiomatization of \ECG, and also a constructive axiomatization of 
Euclidean field theory; but just as in classical geometry, the result is true for any 
reasonable axiomatization. 

The third purpose of this paper,  represented in its title, is to investigate the 
famous Postulate 5 of Euclid.  There are several ways of 
formulating the parallel postulate, which are classically equivalent, but 
not all constructively equivalent.  These versions are called ``Euclid 5'', the 
``strong parallel postulate'', and ``Playfair's axiom.''   The first two each assert that 
two lines actually meet, under slightly different hypotheses about the interior
angles made by a transversal.   Playfair's axiom simply says there cannot be two distinct
parallels to line $L$ through point $p$ not on $L$;  so Playfair's axiom seems to be 
constructively weaker, as it makes no existential assertion.  

The obvious question is:  are these different versions of the parallel postulate
 constructively equivalent, or not?    
 It turns out that Euclid 5 and the strong parallel postulate are equivalent in constructive 
 geometry  (contrary to  \cite{beeson2012}, where it is mistakenly stated that Euclid 5 does not 
imply the strong parallel postulate).   The proof depends on showing that coordinatization 
and multiplication can be defined geometrically using only Euclid 5, so it is somewhat lengthy,
but conceptually straightforward. 

On the other hand, we show that Playfair's axiom does not imply Euclid 5 (or the strong parallel 
axiom).   This is done in two steps:  First, we define a Playfair ring, which
is something like a Euclidean field, but with a weaker requirement about the existence of 
reciprocals.   We show that if $\F$ is a Playfair ring, then $\F^2$ is a model of geometry
with Playfair's axiom.  If Playfair implies Euclid 5, then the axioms of Playfair rings
would imply the Euclidean field axioms.   The second step is to show that the axioms for Playfair rings
are constructively weaker than those for Euclidean fields.  This is proved using 
Kripke models, where the field elements are interpreted as certain real-valued functions.

Because these formal independence results
are proved in the context of Euclidean field theory, they apply to any axiomatizations of geometry 
that have planes over Euclidean fields as models.
 For example, they apply to the Hilbert-style 
axiomatization in \cite{beeson-kobe}, and to the Tarski-style axiomatization in \cite{beeson2015b}.
It is therefore not necessary in this paper to settle upon a particular axiomatization of geometry.
All that we require is that certain basic theorems of Euclid be provable.

I would like to thank Jeremy Avigad for encouraging 
 me to investigate the logical relations between the different parallel postulates.

\subsection{What is constructive geometry?} \label{section:whatisconstructivegeometry}
In constructive mathematics,  if one proves something exists, one has to show how to construct it. 
 In Euclid's geometry, the means of construction are not arbitrary computer programs, but ruler and compass.  Therefore it is natural to look for 
quantifier-free axioms,  with function symbols for the basic ruler-and-compass constructions.  The terms of such a theory
correspond to ruler-and-compass constructions.  These constructions should depend continuously on 
parameters.    We can see that 
dramatically in computer animations of Euclidean constructions, in which one can select some of the original points and 
drag them, and the entire construction ``follows along.''
We expect that if one constructively proves that points forming a certain 
configuration exist,  that the construction can be done ``uniformly'', i.e., by a single construction 
depending continuously on parameters.  

To illustrate what we mean by a uniform construction, we consider an important example.
There are two well-known classical constructions for constructing a perpendicular to line $L$
through point $p$:  one of them is called ``dropping a perpendicular'', and works when $p$ is not on $L$.
The other is called ``erecting a perpendicular'', and works when $p$ is on $L$.   Classically we 
may argue by cases and conclude that for every $p$ and $L$, there exists a perpendicular to $L$
through $p$.  But constructively, we are not allowed to argue by cases.  If we want to prove
that for every $p$ and $L$, there exists a perpendicular to $L$ through $p$, then we must give 
a single, ``uniform'' ruler-and-compass construction that works for any $p$, whether or not $p$ is 
on $L$.

Readers new to constructive mathematics may not understand why an argument by cases is 
not allowed.  Let me explain.  The statement, ``for every $p$ and $L$ there exists a 
perpendicular to $L$ through $p$''  means that, when given $p$ and $L$, we can construct the
perpendicular.  The crux of the matter is what it might mean to be ``given'' a point and a line.
In order to be able to decide algorithmically (let alone with ruler and compass) whether $p$ 
lies on $L$ or not,  we would have to make very drastic assumptions about what it means to be 
``given'' a point or a line.    For example, if we were to assume
 that every point has rational coordinates relative to some 
lines chosen as the $x$ and $y$ axes, then we could compute whether $p$ lies on $L$ or not;
but that would require bringing number theory into geometry.   Anyway, we wouldn't be able to 
construct an equilateral triangle.  We want to allow an open-ended concept of ``point'',
permitting at least the interpretation in which points are pairs of real numbers; and then 
there is no algorithm for deciding if $p$ lies on $L$ or not.  That is why an argument by 
cases is not allowed in constructive geometry.

The other type of argument that is famously not allowed in constructive mathematics is proof
by contradiction.  There are some common points of confusion about this restriction.  The main thing one is not 
allowed to do is to prove an existential statement by contradiction.  For example,  we are not 
allowed to prove that there exists a perpendicular to $L$ through $x$ by assuming there is none,
and reaching a contradiction.   From the constructive point of view, that proof of course 
proves {\em something},  but that something is weaker than existence.  We write it 
$\neg \neg \exists$,  and constructively, you cannot cancel the two negation signs.

However, you {\em are}  allowed to prove inequality or order relations between points
by contradiction.  Consider for the moment two points $x$ and $y$ on a line.  If we 
derive a contradiction from the assumption $x\neq y$ then $x=y$.   This is taken as 
an axiom of constructive geometry: $\neg\, x\neq y \implies x=y$.  In the style of Euclid:
{\em Things that are not unequal are equal.} It has its intuitive grounding in the idea that 
$x=y$ does not make any existential statement.  Similarly,  we are allowed to prove $x < y$
by contradiction; that is, we take $\neg \neg\, x < y \implies x < y$ as an axiom of constructive
geometry.  (That can also be written $\neg\, y \le x \implies x < y$.)
Without this principle, Euclidean geometry would be complicated and nuanced, since
arguments using it occur in Euclid and are not always avoidable. 
 Since order on a line is 
not a primitive relation, we take instead the corresponding axiom for the betweenness relation,
namely $ \neg \neg\, \B(a,b,c) \implies \B(a,b,c)$.  
These two axioms are called the ``stability'' of equality and betweenness.%
\footnote{In adopting the stability of betweenness as an axiom, we should perhaps offer
another justification than the pragmatic necessity for Euclid.  The principle reduces to 
a question, given two points $s$ and $t$, can we find two circles with centers $s$ and $t$
that separate the two points?  The obvious candidate for the radius is half the distance
between $s$ and $t$, so the question boils down to Markov's principle for numbers: if that 
radius is not not positive, must it be positive?  The stability of betweenness thus has 
the same philosophical status as Markov's principle:  it is self-justifying,  not
provable from other constructive axioms, and leads to no trouble in constructive mathematics,
while simplifying many proofs.}

The theorems of Euclidean geometry all have a fairly simple logical form:  Given some 
points, lines, and circles bearing certain relations,  then there exist some further points
bearing certain relations to each other and the original points.   This logical simplicity
implies that (although this may not be obvious at first consideration) 
if we allow the stability axioms, then essentially the only differences between classical and 
constructive geometry are the two requirements:
\begin{itemize}
\item You may not prove existence statements by contradiction; you must provide a construction.
\item The construction you provide must be uniform; that is, it must be proved to work without an argument by cases.
\end{itemize}

Sometimes, when doing constructive mathematics, one may use a mental picture in which one imagines
a point $p$ as having a not-quite-yet-determined location.  For example, think of a point $p$ which 
is very close to line $L$.  We may turn up our microscope and we still can't see whether $p$ is 
or is not on $L$.  We think ``we don't know whether $p$ is on $L$ or not.''   Our construction 
of a perpendicular must be visualized to work on such points $p$.  Of course, this is just a 
mental picture and is not used in actual proofs.  It can be thought of as a way of 
conceptualizing ``we don't have an algorithm for determining whether $p$ is on $L$ or not.''

We illustrate these principles with a second example.
Consider the problem of finding the reflection of point $p$ in line $L$.
Once we know how to construct a perpendicular to $L$ through $p$, it is still not trivial to 
find the reflection of $p$ in $L$.  Of course, if $p$ is on $L$, then it is its own reflection,
and if $p$ is not on $L$, then we can just drop a perpendicular to $L$, meeting $L$ at the foot $f$,
and extend the segment {\em pf} an equal length on the other side of $f$ to get the reflection.  But what about 
the case when we don't know whether $p$ is or is not on $L$?   Of course, that sentence technically 
makes no sense;  but it illustrates the point that we are not allowed to argue by 
cases.  The solution to this problem may not be immediately obvious; but in the body of this 
paper, we will exhibit a uniform construction of the reflection of $p$ in $L$.

\subsection{Is Euclid constructive?}
Our interest in constructive geometry was first awakened by  computer 
animations of Euclid's constructions; we noticed that the proposition in Euclid~I.2
does not depend continuously on parameters, but instead depends on an argument by cases.%
\footnote{Should any reader not possess a copy of Euclid, we recommend \cite{euclid2007} and \cite{thebones}; or 
for those wishing a scholarly commentary as well, \cite{euclid1956}. There are also several 
versions of Euclid available on the Internet for free.}
Proposition I.2 says that
for every points $a$, $b$, and $c$,
there is a point $d$ such that $ad=bc$ (in the sense of segment congruence).
The computer animation of Euclid's construction of $d$ shows that as $b$ spirals inwards to $a$,
$d$ makes big circles around $a$, and so does not depend continuously on $b$ as $b$ goes to $a$.
Classically, we can argue by cases, taking $d=a$ if $b=a$, but Euclid's proof of I.2 is 
non-constructive.%
\footnote{Prop.~I.2 can be views as showing that a ``collapsible compass'' can simulate a 
``rigid compass.''.  
A ``rigid compass''  can be used to transfer a given distance, say $bc$,  from one location in the plane to another, so it enables construction of a circle ``by center and radius.''   A ``collapsible compass'' 
can only be used to construct a circle ``by center and point'', i.e., to draw a circle with a given center $a$ passing through another given point $b$.  The compass ``collapses''
as soon as you pick it up.  There is a word in Dutch, {\em passer}, for a rigid compass, which was used in navigation 
in the seventeenth century.  But there seems to be no single word in English that distinguishes either of the two types 
of compass from the other. 
}

After discovering the non-constructivity of the proof of I.2, 
we then read Euclid with an eye to the constructivity of the proofs.  
There is no other essentially non-constructive proof in Euclid's {\em Elements} I--IV.%
\footnote{In fact, it is only Euclid's proof that is non-constructive,
not the proposition itself.  A constructive proof of Euclid~I.2 is 
given in Lemma~\ref{lemma:I.2}; see also \S5.1 of \cite{beeson2015b}.}

\subsection{Criteria for a constructive axiomatization of geometry}  
In order to meet the second and third purposes of this paper (namely, to relate its
models to planes over Euclidean fields, and to show that Playfair does not imply 
Euclid 5, we need to make use of some formal theory of constructive 
geometry.  There are, of course, as many ways to formalize constructive geometry as 
there are to formalize classical geometry; see page 4 of \cite{beeson-edinburgh}, where we 
enumerate eleven possible choices, more or less independent,  so there are several thousand
different ways to formalize geometry.  We want an axiomatization that is 
\smallskip

\begin{itemize}
\item is faithful to Euclid's {\em Elements}, i.e., it can express Euclid's 
propositions and proofs in a natural way;
\item has terms for geometrical constructions, and
 those terms describe constructions continuous in parameters;
\item uses intuitionistic logic 
\end{itemize}
\smallskip

In addition we want the axioms of our theory to be
\begin{itemize}
\item  quantifier-free.
\item disjunction free.
\end{itemize}
\smallskip

The last two criteria make the proof theory of such an axiomatization
easy:  they permit the applications of cut-elimination that we want to make, in order to 
connect existence proofs to ruler-and-compass constructions.  They are also 
very reasonable criteria.     First,
disjunction does not occur in Euclid.  Euclid did not use first-order logic;  the logical structure of his theorems is 
discussed in \S\ref{section:logicalformofEuclid}, but for now we just note that he never states or uses a disjunctive proposition, i.e., one with ``or''
in the conclusion.%
\footnote{In Prop. I.26, he mentions ``or'', but in the hypothesis, not the conclusion, so the proposition 
can be expressed as two propositions without ``or''.  Prop. III.4 is an example in which 
Euclid chooses a negative formulation, when a modern mathematician might choose a disjunction: if two 
chords of a circle bisect each each other, then one or the other passes through the center.}

We have developed two rather different axiomatizations of \ECG,  each of which has all 
these properties.  The first one, 
presented in \cite{beeson-kobe}, was based on Hilbert's axioms for geometry.  In \cite{beeson2012} 
we mentioned  a version more closely related to Tarski's concise and elegant axiomatization of 
classical Euclidean geometry, but space did not permit a full description there.  That
theory has now been fully described and studied in \cite{beeson2015b}.  The results in 
this paper apply to either axiomatization, and we believe, to all reasonable variations of 
these axiomatizations.   Therefore, in this paper we do not commit to any particular axiomatization
of geometry; all we require is that certain basic theorems of Euclid should be provable and 
the logical rules should be intuitionistic.

\subsection{Basic metamathematical results}
The results of our metamathematical analysis of constructive geometry are:
\begin{itemize}
\item  Things proved to exist in (a formal version of) \ECG\ can be constructed with ruler and compass.
\item  The models of (a formal version of) \ECG\ are exactly planes $\F^2$ over Euclidean fields $\F$,
with two different choices of the axiom about multiplicative inverses in $\F$ corresponding to 
two different versions of the parallel axiom.
\item These   different versions of the parallel axiom are formally inequivalent.
\end{itemize}

We now describe these results in more detail.
Once a particular axiomatization of \ECG\ is chosen (meeting the criteria listed above),
then one can apply cut-elimination to obtain several fundamental properties of such a theory,
as discussed in \cite{beeson-kobe} and \cite{beeson2015b}.   For example, we have the result alluded to above, that existential theorems are ``realized'' by ruler-and-compass constructions.  For a second example, propositions
involving a negative hypothesis and a disjunctive conclusion cannot be proved. Examples include
$ a < b \implies a < x \lor x < b$ and  $x \neq 0 \implies 0 < x \lor x < 0$, which are discussed
more fully in \S\ref{section:order} and Theorem~\ref{theorem:kripke1}.

Next there is Tarski's ``representation problem''.  We want to characterize the models of \ECG.   Classically,  the 
models of Euclidean geometry are well-known to be planes over Euclidean fields (ordered fields in which nonnegative elements 
have square roots).   But constructively, 
what is a Euclidean field?  It turns out that there are two natural ways of defining that, and that planes over
such fields correspond to two different versions of the parallel postulate.  The proof of such 
a representation theorem  (classically or constructively) involves both the ``arithmetization of geometry'' 
and the ``geometrization of arithmetic.''   The geometrization of arithmetic (originally due to Descartes) shows how to add, subtract, 
multiply, divide, and take square roots of segments, and how to introduce coordinates in an arbitrary plane.  The 
arithmetization of geometry (also due to Descartes) is the verification (using coordinates 
and algebraic calculations) that the axioms of 
geometry hold in $\F^2$,  where $\F$ is a Euclidean field.  Both of these offer some constructive difficulties; for 
example, Descartes's constructions apply only to positive segments,  but we need to define arithmetic on signed numbers, without 
making case distinctions whether the numbers are zero or not.   These difficulties are overcome and the representation 
theorem proved in Theorem~\ref{theorem:representation}.

Finally, we  
settle the relations between three 
versions of the parallel postulate (which are explained below).
Two of these  assert the existence of a point of intersection,  while the third (Playfair's axiom)
does not.  The two that assert existence are equivalent,  and they are not implied by Playfair's axiom.  This independence proof (and others) are
 obtained using Kripke models of Euclidean field theory, so they do not depend on the details of the axiomatization 
of constructive geometry, but only on the fact that constructive geometry can be interpreted in Euclidean field theory.

\subsection{Postulates {\em vs.} axioms in Euclid}

 Euclid presents his readers with both ``postulates'' and ``axioms''.  Modern mathematicians often treat these
words as synonyms.  For Euclid and his contemporaries, however, they had quite different meanings.  Here
is the difference, as explained by Pambuccian \cite{pambuccian2006}, p.~130:   %page 12 in my manuscript copy
\begin{quote}
\small
For Proclus, who relates a view held by Geminus, a 
postulate
prescribes that we construct or provide some simple or easily grasped
object for the exhibition of a character, while an axiom asserts some 
inherent
attribute that is known at once to one's auditors. And just as a problem
differs from a theorem, so a postulate differs from an axiom, even though
both of them are undemonstrated; the one is assumed because it is easy to
construct, the other accepted because it is easy to know. That is, 
postulates
ask for the production, the {\em poesis} of something not yet 
given$\ldots\,$, whereas axioms refer to the {\em gnosis} of a given, to 
insight into the validity of certain relationships that hold between given notions.
\end{quote}

\subsection{The parallel postulate}
\begin{figure}[ht]
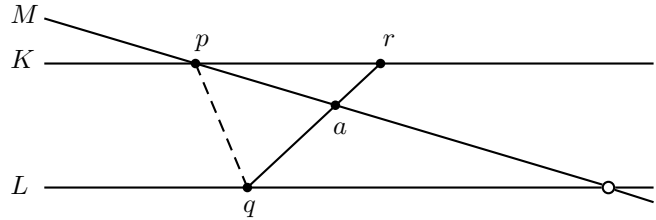
   
\center{
\EuclidParallelFigure
}
\caption{Euclid 5:  $M$ and $L$ must meet on the right side, provided $\B(q,a,r)$ and $pq$ makes
alternate interior angles equal with $K$ and $L$.
\label{figure:EuclidParallelFigureOne}
}
\end{figure}

Euclid's famous ``parallel postulate'' (Postulate 5, henceforth referred to as Euclid 5) states that if two lines $L$ and $M$ are traversed by another line $T$, forming 
interior angles on one side of $T$ adding up to less than two right angles, 
then $L$ and $M$ will intersect on that side of $T$.  The version of Euclid 5 given
in Fig.~\ref{figure:EuclidParallelFigureOne} shows how to state this without referring 
to ``addition of angles.'' 

\noindent
Euclid 5 makes an assertion about the existence of the point of intersection of $K$ and $L$,
and hence it can be viewed as 
a construction method for producing certain triangles.   In view of the remarks of Geminus and Proclus,
it seems likely that Euclid viewed his postulate in this way, or he would have called it an axiom.

\begin{figure}[ht]
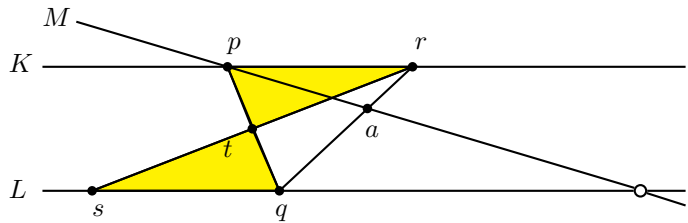
 
\center{ \EuclidParallelRawFigure}
\caption{Euclid 5:  $M$ and $L$ must meet on the right side, provided $\B(q,a,r)$ and $pt=qt$ and $rt=st$.}
\label{figure:EuclidParallelRawFigureOne}
\end{figure}

Fig.~\ref{figure:EuclidParallelRawFigureOne} shows how to eliminate all mention of angles 
in the formulation of Euclid 5.

\begin{figure}[ht]
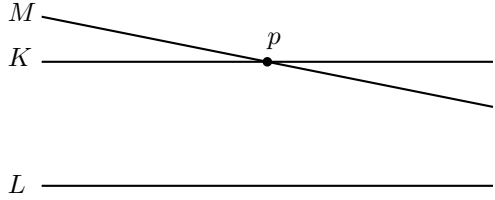
    
 \center{\PlayfairFigure}
\caption{Playfair:  if $K$ and $L$ are parallel, $M$ and $L$ can't fail to meet. }\label{figure:PlayfairFigureOne}
\end{figure}

 In 1795, Playfair introduced the version that is 
usually used today,  which is an axiom rather than a postulate: 
Given a line $L$, and a point $p$ not on $L$, there 
cannot be two distinct lines through $p$ parallel to $L$.
(Parallel lines are by definition lines in the same plane that do not meet.)   
Unlike Euclid 5,
Playfair's axiom does not assert the existence of any specific point.  
See Fig.~\ref{figure:PlayfairFigureOne}.

\FloatBarrier

To compare different versions of the parallel postulate, we begin with 
a fundamental observation.  Given point $p$ not on line $L$,  and a line $T$
through $p$ that meets $L$, there is (without using any parallel postulate at all)
exactly one line $M$ through $p$   that makes
the sum of the interior angles on one 
side of the transversal $T$ equal to two right angles (or makes alternate interior angles equal).
This line $M$ is parallel to $L$,  since if $M$ meets $L$ in $q$
then there is a triangle with two right angles,
contradicting Euclid I.17 (and Euclid does not use the parallel postulate until I.29).
 
With classical logic, Playfair's axiom implies Euclid 5, as follows. Suppose  given point $p$  not on 
line $L$,  and line $T$ through $p$ meeting $L$, and line $K$ through $p$ such that the 
angles formed on one side of $T$ by $K$ and $L$ are together less than two right angles.
Let $M$ be the line through $p$ that does make angles together equal to two right angles with $T$.
As remarked above, $M$ is parallel to $L$,
  Hence, by Playfair,  $K$ cannot be parallel 
to $L$.  Classically, then, it must meet $L$, which is the conclusion of Euclid 5.  But this 
proof by contradiction does not provide an explicit point of intersection.

That raises the question whether Playfair's axiom implies Euclid 5 using only intuitionistic (constructive)
logic.   This question is resolved in this paper.
The answer is in the negative:  Playfair does not imply Euclid 5.
 
\begin{figure}[h]
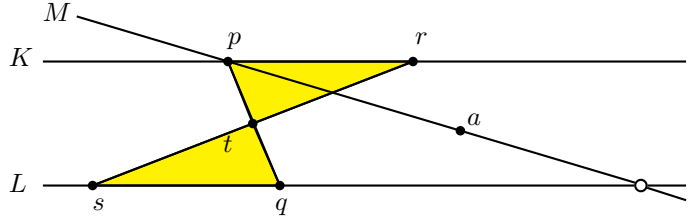
   
\center{\StrongParallelRawFigure}
\caption{Strong Parallel Postulate:  $M$ and $L$ must meet (somewhere) provided $a$ is not on $K$ and and $p$ is not on $L$ and $pt=qt$ and $rt=st$ sand $r\neq p$.}
\label{figure:StrongParallelRawFigureOne}
\end{figure}

  The following is the {\em strong parallel postulate}.
\begin{quote}
If $p$ is not on $L$,  $q$ is on $L$, line $K$ passes through $p$,  and $pq$ makes
alternate interior angles equal with $K$ and $L$, then any line $M$ through $p$ that 
is distinct from $K$ meets $L$.
  See Fig.~\ref{figure:StrongParallelRawFigureOne}.
\end{quote}
This is like Euclid 5, except that the hypothesis is weaker (in that it does not 
require us to specify on which side of $T$ the interior angles are less than two right angles),
and the conclusion is also weaker (not saying on which side of $q$ the two lines meet).
Technically, the hypothesis $B(q,a,r)$ in Fig.~\ref{figure:EuclidParallelFigureOne} is replaced
by simply requiring $a$ to not lie on $K$.

Both the strong parallel postulate and Euclid 5 make existence assertions.  It is not 
hard to show that the strong parallel postulate implies Euclid 5 (hence the name).
Much less obviously, the strong parallel postulate turns out to be equivalent to Euclid 5.  This is proved by developing coordinates, and 
the geometric definitions of addition and multiplication using only Euclid 5, instead of the 
strong parallel postulate.  Once coordinates and arithmetic can be defined, then the 
strong parallel postulate boils down to the property that nonzero points on the $x$-axis have
multiplicative inverses, and Euclid 5 says that positive elements have multiplicative inverses;
but $1/x = x/\vert x \vert^2$, so they are equivalent.  The delicate question is then 
whether Euclid 5 suffices for constructive (case-free) definitions of coordinates and 
arithmetic.  We will show that it does.

In the presence of Playfair's axiom,  the strong parallel postulate is equivalent to 
\begin{quote} If $M$ and $L$ are distinct, non-parallel lines,
then $M$ meets $L$.
\end{quote}
This equivalence (proved in Lemma~\ref{lemma:linestability}) separates the parallel postulate into two parts: a negative statement 
about parallelism (Playfair)  and an existential assertion about the result of a construction
(the two lines meet in a point that can be found with a ruler,  albeit perhaps a {\em very long} ruler).

\subsection{Triangle circumscription as a parallel postulate}
In classical geometry, there are many propositions equivalent to the parallel postulate; for 
example, the triangle circumscription principle, which says that given three non-collinear points $a$, $b$, and $c$,
there is a fourth point equidistant from $a$, $b$, and $c$ (the center of a circle passing through those
three points).  Szmielew used the triangle circumscription principle as an axiom in 
her lectures, although in \cite{schwabhauser}, a different formulation of the parallel
axiom due to Tarski is used.    In Theorem~\ref{theorem:trianglecircumscription} below, we show that the triangle 
circumscription principle is constructively equivalent to the strong parallel postulate,  
 and in \cite{beeson2015b}
we follow Szmielew in adopting this as an axiom.  There is a natural ruler and compass 
construction of the center point of the circle through $a$, $b$, and $c$, as the intersection point of the perpendicular bisectors
of $ab$ and $bc$.

\subsection{Past work on constructive geometry}
L.~E.~J.~Brouwer, the founder of intuitionism,
 apparently held the axiomatic approach to mathematics in low esteem; at least,
he never took that approach in his papers. It is therefore a bit surprising that his 
most famous student, Arend Heyting, wrote his dissertation on an axiomatic (and intuitionistic)
treatment of projective geometry (published two years later as \cite{heyting1927} and 
again 34 years later as \cite{heyting1959}, which is probably the 
same as the easier-to-locate \cite{heyting1959b}).
 Heyting took a relation of ``apartness'' as 
primitive.  This is a ``positive notion of inequality''.  The essence of the notion can be 
captured without a new symbol in the following axiom about order on the line:
$$ a < b \implies x < b \lor a < x$$
Heyting's work was so influential that every subsequent paper about intuitionistic or constructive 
geometry has taken an apartness relation as primitive; it was still being discussed in 1990 and 1996
(see \cite{vanDalen1990,vanDalen1996}).   A somewhat different tack was taken by von Plato  
 \cite{vonPlato1995, vonPlato1998}, who worked on axioms for affine geometry (no congruence).
 
Lombard and Vesley \cite{lombard-vesley}
gave a constructive theory of geometry, perhaps the first to apply to 
Euclidean geometry.  They followed Heyting in taking apartness as primitive, and they wished 
to avoid having equality as a primitive, which led them to the six-place relation ``the sum 
of the lengths of $ab$ and $cd$ is more than the length $pq$.''  They were able to recover 
betweenness and congruence from this relation and give a realizability interpretation.
 
 The unpublished dissertation \cite{kijne}
is about constructive geometry in the sense of geometric constructions, but explicitly non-constructive in that it makes use of ``decision operations'', or ``branching operations''.   
As \cite{lombard-vesley} points out, geometers sometimes use the word ``constructive'' to mean
that the axioms are quantifier-free, i.e., function symbols are used instead of existential quantifiers.
That is neither necessary nor sufficient for a constructive theory in the sense of ``constructive
mathematics''.   Nevertheless our constructive geometry does have quantifier-free axiomatizations.
The first quantifier-free axiomatization of geometry may have been \cite{suppes}; but that was 
a theory with classical logic.

The use of apartness stemmed from a reluctance to apply classical logic to the equality of points.
In our constructive geometry, we assume the stability of equality, and we assume Markov's 
principle $\neg\, x \le y \implies y < x$ (but expressed using betweenness).  This   allows
one to prove the equality or inequality of two points by contradiction.   The strong parallel 
postulate, as we show, allows one to prove that two lines intersect by contradiction.   What is not 
allowed is an argument by cases,  where the cases concern ordering or equality relations between points.  Also, of course, existential assertions are supposed to be proved by explicit construction,
which, coupled with the prohibition on arguments by cases, requires a {\em uniform} construction,
continuous in parameters.
That contrasts with geometry based on apartness, which is designed to {\em allow} certain arguments
by ``overlapping'' cases.  No theory with apartness can have the property that points it proves 
to exist depend continuously on parameters,  which is an important feature of constructive geometry
as developed here.

\section{Euclid's reasoning considered constructively}
In the late twentieth century, contemporaneously with the flowering of computer science, there was a new surge of vigor in algorithmic, or constructive, mathematics,  beginning with Bishop's book \cite{bishop}.  
In algorithmic mathematics, one tries to reduce every ``existence theorem'' to an assertion that a certain algorithm has a certain 
result.  It was discovered by Brouwer that if one restricts the laws of logic suitably (to ``intuitionistic'' logic), 
then one only obtains algorithmic existence theorems, so there is a fundamental connection between methods of proof, and 
the existence of algorithms to construct the things that have been proved to exist.  Brouwer thought it necessary 
to do more than just restrict logic; he also wanted to state some additional principles.  Bishop renounced additional principles
and worked by choosing his definitions very carefully, but using only a restricted form of logic.  Results obtained by Bishop's methods
 are classically valid as well as constructively. 

What happens if we examine Euclid's {\em Elements} from this point of view? 
It turns out that the required changes are few and minor.  Euclid's proof of Prop.~ I.2 
is non-constructive (but the theorem itself has a different and constructive proof, given 
in Lemma~\ref{lemma:I.2}), and  
and the parallel axiom needs a more explicit formulation.
 Euclid is essentially constructive as it stands.  In \S\ref{section:development},
we justify this conclusion in more detail.

\subsection{Order on a line from the constructive viewpoint}\label{section:order}
In this section we explain how a constructivist views the relations $x < y$ and $x \le y$ on the real line.   This 
section can be skipped by readers familiar with constructive order, but it will be very helpful in understanding 
constructive geometry.  Order on a line can be thought of as one-dimensional geometry, so it 
makes sense for the reader new to constructivity to start with the one-dimensional case.
 The theory of order 
 translates directly into the betweenness relation in geometry, since for positive $x$ and $y$ we have that 
$x < y$ is the same as ``$x$ is between 0 and $y$''.   In constructive mathematics, the real numbers are given by constructive 
sequences of rational approximations.  For example, we could take Bishop's definition 
of a  real number as a (constructively given) sequence $x_n$ of rationals such that for every $n \ge 1$, 
$$ \vert x_n - x_m \vert < \frac 1 n + \frac 1 m .$$
This guarantees that the limit $x$ of such a sequence satisfies $\vert x-x_n \vert \le 1/n$. 
Two such sequences are considered equal if $x_n - y_n \le 2/n$;  it is important that different 
sequences can represent equal (or ``the same'', if you prefer) real numbers.  Now observe that 
if $x < y$,  we will eventually become aware of that fact by computing fine enough approximations to $x$ and $y$,
so that $x_n + \frac 1 n < y_n - \frac 1 n$.  Then we will have an explicit positive lower bound on $y-x$, which is what is required to 
assert $x < y$.%

The ``recursive reals'' are defined this way, but specifying that the sequences $x_n$ are to be computed by 
a Turing machine (or other formal definition of ``computer program'').  In more detail:
We write
$\{e\}(n)$ for the result, if any, of the $e$-th Turing machine at input $n$.   Rational numbers are coded as certain integers,
and modulo this coding we can speak of recursive functions from $\N$ to $\Q$.  A ``recursive real number'' $x$ is given by an index of a Turing 
machine $e$, whose output $\{e\}(n)$ at input $n$ codes a rational number, such that $\vert \{e\}(n) - x \vert \le 1/n$ for each $n \in \N$.  We can, if we wish, avoid assuming that the real numbers
are given in advance, by considering the Turing machine index $e$ to {\em be} the real number $x$.

It may be helpful to keep the ``model'' of recursive reals in 
mind, but in the spirit of ``tables, chairs, and beer mugs'',  modern constructivists often prefer not to make this 
commitment, but refer in the abstract  to ``constructive sequences.''  Thus the ``classical model'' is always a possible 
interpretation.

Given two (definitions of) numbers $x$ and $y$, we may compute approximations to $x$ and $y$ for many years and still be uncertain whether $x=y$ or not.  Hence, we are not allowed
to constructively assert the trichotomy law $x <y \lor x = y \lor y < x$, since we have no way to make the decision 
in a finite number of steps.  One can prove  that the recursive reals definitely cannot constructively be proved to  
satisfy the trichotomy law, as that would imply a computable solution to the halting problem.  But without a
commitment to a definite definition of ``constructive sequence'',  the most we can say is that ``we cannot assert''
trichotomy.  This phrase ``we cannot assert'' in constructive mathematics is usually code for,  ``it fails 
in the recursive model.''

Note, however, that we may be able to assert the trichotomy law for various subfields 
of the real numbers.  For example, it is valid for the rational numbers, and it is also valid (though less obviously)
for the real algebraic numbers.  In each of these cases, the elements of the field are given by finite objects, that
can be presented to us ``all at once'', unlike real numbers;  but that property is not sufficient for trichotomy 
to hold,  since the recursive reals fail to satisy trichotomy, but a recursive real can be given ``all at once'' 
by handing over a Turing machine to compute the sequence $x_n$ of approximations.

The relation $x \le y$ is equivalent to $\neg (y < x)$, either by definition, or by a simple theorem if one defines 
$x \le y$ in terms of approximating sequences.   It is definitely not equivalent to $x < y \lor x = y$ (see the refutation 
of the trichotomy law in the recursive model given below). But now consider negating $x \le y$.  Could we assert
$\neg x \le y$ implies $y < x$?   Subtracting $y$ we arrive at an equivalent version of the question with only one 
variable:  can we assert $\neg x \le 0$ implies $ 0 < x$?   Since $x \le 0$ is equivalent to $\neg 0 < x$, the 
question is whether we can assert
\smallskip

\axioms
\qquad $ \neg \neg x > 0 \implies x > 0$ &(Markov's principle)
\endaxioms
\smallskip

In other words,  is it legal to prove that a number is positive by contradiction?  One could argue for this principle 
as follows:  Suppose $\neg \neg x > 0$.  Now 
compute the approximations $x_n$ one by one for $n=1,2,\ldots$.  
Note that trichotomy does hold for $x_n$ and $1/n$, both of 
which are rational.  You must find an $n$ such that $x_n > 1/n$,
since otherwise for all $n$, we have $x_n \le 1/n$, which means $x \le 0$, contradicting $\neg \neg x > 0$.
Well, this is a circular argument:  we have used Markov's principle in the justification of Markov's principle.

Shall we settle it by looking at the recursive model?  There it can easily be shown to boil down to this:  if 
a Turing machine cannot fail to halt,  then it halts.  Again one sees no way to prove this, and some may feel 
is intuitively true, while others may not agree.

Historically, the Russian constructivist school adopted Markov's principle, and the Western constructivists did not.
It reminds one of the split between the branches of the Catholic Church,  which also took place along geographical lines.
In any case, as discussed in \S\ref{section:whatisconstructivegeometry} 
and \S\ref{section:casesplitsinEuclid}, it seems appropriate to adopt this principle for a constructive treatment of Euclid.

From the constructive viewpoint, the main difference between $x < y$ and $x \le y$, as applied 
to real numbers, is that $x<y$ involves an 
existential quantifier; it contains the assertion that we can find a rational lower bound on $\vert y-x \vert$, 
while $x \le y$ is defined with a universal quantifier, and contains no hidden assertions.  If we take any 
formula involving inequalities,  and replace $x < y$ with $\neg\, y \le x$,  we obtain a classically equivalent
assertion no longer containing an existential quantifier.   If in addition we replace $A \lor B$ by 
$\neg\,(\neg A \land \neg B)$,  we will have eliminated all hidden claims, and the result will be classically  valid
if and only if it is constructively valid.   To understand constructive mathematics,  one has to learn to 
see the ``hidden claims'' that are made by disjunctions and existential quantifiers,  which can make a formula 
``stronger'' than its classical interpretation.   Of course, if existential quantifiers or disjunctions occur
in the hypotheses, then a stronger hypothesis can make a weaker theorem. 

A given classical theorem might have more than one (even many) classically equivalent versions with different 
constructive meanings.  Therefore,  finding a constructive version of a given theory is often a matter of choosing 
the right definitions and axioms.

Consider the following proposition, which is weaker than trichotomy:  $$ x \le 0 \lor x \ge 0.$$  This   is 
also not constructively valid.  Intuitively,  no matter how long we keep computing approximations $x_n$, if they 
keep coming out zero we will never know which disjunct is correct.  As soon as we stop computing, the very next
term might have told us.   

We now show that both trichotomy and $x\le 0 \lor x \ge 0$  fail in the recusive reals, if disjunction is interpreted 
as computable decidability.  First consider trichotomy.
We will show that  there is no computable test-for-equality function, that is, no computable function $D$ that operates on 
two Turing machine indices $x$ and $y$,  and produces $0$ when $x$ and $y$ are equal recursive real numbers (i.e., have the same limiting value), and 1 when
they are recursive real numbers with different values.  Proof, if we had such a $D$, we could 
 solve the halting problem by applying $D$
to the point $(E(x),0)$, where $\{E(x)\}(n) = 1/n$ if Turing machine $x$ does not halt at input $x$ in fewer than $n$ steps, and $\{E(x)\}(n) = 1/k$ 
otherwise, where $x$ halts in exactly $k$ steps.  Namely, $\{x\}(x)$  halts if and only if the value of $E(x)$ is not zero, if and only if 
$D(Z,E(x)) \neq 0$, where $Z$ is an index of the constant function whose value is the (number coding the) rational number zero.
 
 Now consider the proposition $x \ge 0 \lor x \le 0$.  We can imitate the above construction, but replacing the halting problem by two recursively inseparable r.e. sets $A$ and $B$,
and making the number $E(x)$ be equal to $1/n$ if at the $n$-th stage of computation we see that $x \in A$, and $-1/n$ if 
we see $x \in B$;  so if $x$ is in neither $A$ nor $B$, $E(x)$ will be equal to zero.  Hence $x \ge 0 \lor x \le 0$ fails 
to hold in the recursive reals. 

Finally we consider this proposition: 
\smallskip

\axioms
\qquad $x \neq 0 \implies x < 0 \lor x > 0$ & (two-sides) 
\endaxioms
\smallskip

We call this principle ``two-sides'' since it is closely related to ``a point not on a line is on one side or the other 
of the line''.  (Here the ``line'' could be the $y$-axis.)   Two-sides does hold in the recursive model, since, assuming 
that $x \neq 0$, if we compute $x_n$ for large enough $n$, eventually we will find that $x_n + 1/n < 0$ or $x_n - 1/n > 0$,
as that is what it means in the recursive model for $x$ to be nonzero.  But since $x_n$ and $1/n$ are rational numbers,
we can decide computably which of the disjuncts holds, and that tells us whether $x < 0$ or $x > 0$.  

On the other hand, this verification of two-sides in the recursive model is not a proof that it is constructively valid;
the following lemma shows that it is at least as ``questionable'' as Markov's principle:

\begin{Lemma} two-sides implies Markov's principle (with intuitionistic logic).
\end{Lemma}

\noindent{\em Proof}. 
Suppose $\neg \neg\, x > 0$, the hypothesis of Markov's principle.  Then $x \neq 0$, so by two-sides $x < 0$ 
or $x > 0$; if $x < 0$ then $\neg x > 0$, contradicting $\neg \neg\, x > 0$; hence the disjunct $x < 0$ is impossible. 
Hence $x > 0$, which is the conclusion of Markov's principle.  That completes the proof of the lemma.

If we assume that points on a line correspond to Cauchy sequences of rational numbers, then we can also prove 
the converse:

\begin{Lemma} If real numbers are determined by Cauchy sequences, then
Markov's principle implies two-sides.
\end{Lemma}

\noindent{\em Proof}.  We may suppose that real numbers are given by ``Bishop sequences'' as described above, as 
it is well-known to be equivalent to the Cauchy sequence definition.  A Bishop sequence for $\vert x \vert$ 
is given by $\vert x \vert_n = \max(x_n,-x_n)$.   Suppose $x \neq 0$ (the hypothesis of two-sides).
  Then $\vert x \vert  \neq 0$.  We claim $\vert x \vert > 0$.  By Markov's principle, it suffices to 
derive a contradiction from $\vert x \vert \le 0$.   Suppose $\vert x \vert \le 0$.  Then $x =0$, 
contradicting $\vert x \vert \neq 0$.  Hence by Markov's principle, $\vert x \vert > 0$.   Then 
by definition of $>$, for some $n$ we have $\max(x_n,-x_n) > 1/n$.  But $x_n$ is a rational number, so 
$x_n > 1/n$ or $-x_n > 1/n$.  In the former case we have $x > 0$; in the latter case $x < 0$.  But this
is the conclusion of two-sides.   That completes the proof.

Since Markov's principle is known to be unprovable in the standard intuitionistic formal theories of arithmetic
and arithmetic of finite types (see \cite{troelstra}, pp. 213 {\em ff.}), two-sides is also not provable in these theories.
However, in the context of geometry we do {\em not} assume that points are always given to us as Cauchy sequences
of rationals;  we do not even assume that we can always construct a Cauchy sequence of rationals corresponding to 
a given point.  (There are non-Archimedean models of elementary geometry, for example.)  The lemma therefore 
does not imply that two-sides and Markov's principle are equivalent in geometry; and indeed  that 
is not the case.   Two-sides is not provable in our geometric theory,  even though we adopt Markov's principle
as an axiom. 
 
It is an open (philosophical) question whether 
our geometrical intuitions compel us to accept Markov's principle, or whether they compel us, having done that, to 
also accept two-sides.%
\footnote{Brouwer did not accept two-sides (or Markov's principle, for that matter), as a discussion about the creative subject 
on p.~492 of \cite{brouwer1948c} makes clear, although the principle is not explicitly stated there.
}  
In this paper we take a pragmatic approach:  geometry without Markov's principle will be 
more complicated, but it is possible without undue complexity to consider two-sides as an added principle, 
which can be accepted or not.
We choose not to   accept it, since (as we show here) it is not required if our goal is to prove the theorems in Euclid or 
of the type found in Euclid.%
\footnote{A philosophical argument can be made that Markov's principle in geometry is related
to Hilbert's ``density axiom'', according to which there exists a point $c$ strictly between any 
two distinct points $a$ and $b$.  For, if $\neg \neg a < b$, then $b\neq a$,
so by the density axiom there is a point $c$ between $a$ and $b$, and the circle with center
$b$ passing through $c$ shows that $a < b$.  But this argument may also be circular, for how 
are we to justify the density axiom? The obvious justification is to take $c$ to be the midpoint
of segment $ab$, but constructing the midpoint uniformly without assuming $a \neq b$ is problematic,
as we discuss below.}
\subsection{Logical form of Euclid's propositions and proofs} \label{section:logicalformofEuclid}

One should remember that Euclid did not work in first-order logic.  This is not because, like Hilbert, he used 
set-theoretical concepts that go beyond first order.   It is instead because he does not use any nested quantifiers
or even arbitrary Boolean combinations of formulas.  All Euclid's propositions have the form, given some points bearing 
certain relations to each other, we can construct one or more additional points bearing certain relations to the original 
points and each other.  A modern logician would describe this by saying that Euclid's theorems have the form, 
a conjunction of literals implies another conjunction of literals,  where a literal is an atomic formula or the negation 
of an atomic formula.  One does not even find negation explicitly in Euclid;  it is hidden in the hypothesis that two 
points are distinct.   Often even this wording is not present, but is left implicit.   

In particular there is no disjunction to be found in the conclusion of any proposition in Euclid. 
Note that a disjunction in the {\em hypothesis} of theorem is inessential:
$(P \lor Q) \implies R$ is equivalent to the conjunction of $P\implies R$ and $Q \implies R$.
This kind of eliminable disjunction occurs implicitly in Euclid, because in
 some of his propositions, a complete proof would include 
an argument by cases, and Euclid handles only one case.  For example, we have already 
discussed Prop.~I.2, where Euclid shows that given $b \neq c$, and 
given $a$, it is possible to construct $d$ with $ad =bc$.   Euclid does not mention the case when $a=b$, presumably 
because it was obvious that in that case one can take $d=c$, and similarly for the case $a=c$.    Euclid was already criticized for this sloppiness about case distinctions
thousands of years ago.  Our point here is that the failure to write out a separate proof 
of every case is not relevant to our claim that Euclid is disjunction-free. 

Euclid's proofs have been analyzed in detail by Avigad {\em et.~al.} in \cite{avigad2009}, and they conclude:

\begin{quote}
Euclidean proofs do little more than introduce objects
satisfying lists of atomic (or negation atomic) assertions, and then draw further
atomic (or negation atomic) conclusions from these, in a simple linear fashion.
There are two minor departures from this pattern. Sometimes a Euclidean proof
involves a case split; for example, if $ab$ and $cd$ are unequal segments, then one
is longer than the other, and one can argue that a desired conclusion follows in
either case. The other exception is that Euclid sometimes uses a {\em reductio}; for
example, if the supposition that $ab$ and $cd$ are unequal yields a contradiction
then one can conclude that $ab$ and $cd$ are equal.
\end{quote}
``{\em Reductio}'' refers to {\em reductio ad absurdum}, which means proof by contradiction.

\subsection{Case splits and {\em reductio} inessential in Euclid} \label{section:casesplitsinEuclid}
It is our purpose in this section to argue that Euclid's reasoning can be supported
in \ECG, including the two types of apparently non-constructive reasoning just mentioned.
The reason for this is that the case splits and {\em reductio} arguments in Euclid
 can be made constructive using the ``stability''
of equality and betweenness.    By the stability of equality, we mean 
$$\neg \neg\, x = y \implies x = y.$$
This formula simply codifies the principle that it is legal to prove equality of two points
by contradiction, and we take it as a fundamental principle.  In words: {\em Things that are not unequal are equal.}
Similarly, if $\B(a,b,c)$ means that $b$ is between $a$ and $c$ on a line, we take 
as an axiom the stability of betweenness:
$$ \neg \neg\, \B(a,b,c) \implies \B(a,b,c).$$
While Euclid never explicitly mentions betweenness (which is, as is well-known, the main
flaw in Euclid),  the stability of betweenness and equality together account for all apparent
instances of nonconstructive arguments in Euclid.

A typical example of such an argument in Euclid is Prop. I.6, whose proof begins
\begin{quote}
Let $ABC$ be a triangle having the angle $ABC$ equal to the angle $ACB$.
  I say that the side $AB$ is also equal to the side $AC$.
  For, if $AB$ is unequal to $AC$, one of them is greater.
  Let $AB$ be greater, $\ldots$
\end{quote}
To render this argument in first-order logic, we have to make sense of Euclid's ``common notion'' of 
quantities being greater or less than other quantities.  This is usually done formally using 
the betweenness relation.  In the case at hand, we would lay off segment $AB$ along ray $AC$, 
finding point $D$ on ray $AC$ with $AD = AC$.  Then since $AB$ is unequal to $AC$,
$D \neq C$.   Then ``one of them is greater'' 
becomes the disjunction  $\B(A,D,C) \lor \B(A,C,D)$.
This disjunction is not constructively valid.   But its double negation is valid, 
since if the disjunction were false, we would have 
$ \neg \B(A,D,C)$ and $\neg \B(A,C,D)$, which would contradict $D \neq C$.   
The rest of Euclid's argument shows that each of the disjuncts implies the desired conclusion.
We therefore conclude that the double negation of the desired conclusion is valid.%
\footnote{Those new to intuitionistic logic may need more detail: if $P \implies R$ 
and $Q \implies R$, then $P \lor Q \implies R$.   Now we can double-negate both 
sides of an implication, using the constructively valid law
$\neg\neg\,(S\implies R) \implies (\neg\neg\, S\implies \neg\neg\,R)$,
obtaining $ \neg\neg\,(P \lor Q) \implies \neg\neg\, R$.  Hence, if we 
know $\neg \neg\, (P \lor Q)$, we can conclude $\neg\neg\, R$.}
The conclusion of I.6, however, is negative (has no $\exists$ or $\lor$).  
  Hence the double negation can be pushed inwards, and then
  the stability of the atomic 
sentences can be applied to make the double negations disappear.

 In  Euclid, disjunctions never appear, even implicitly, in the conclusions of propositions.%
 \footnote{An apparent exception to this rule is Prop. I.13, ``If a straight line set up on a straight 
 line make angles, it will make either two right angles or angles equal to two right angles.''
 Here the disjunction is superfluous: we can just say, ``it will make angles equal to two right angles.''
}
The conclusions are simply conjunctions of 
literals.   Sometimes there is an implicit existential quantifier, but if (as will always be the case
in our formalizations of constructive geometry) we can explicitly exhibit terms for the constructed points, the resulting explicit form 
of the proposition will be quantifier-free.   Then the double negation can be pushed inwards
as just illustrated.  In this way, all arguments of the 
form   beginning ``For, if $AB$ is unequal to $AC$, one of them is greater'',  
can be constructivized.

Prop.~I.26 gives an example of the use of the stability of equality: 
 ``$\ldots DE$ is not unequal to $AB$, and is therefore equal to it.''
It that example, and all other examples in Euclid, 
the following general reasoning applies:   Since the conclusion concerns the equality of certain points,  we can 
simply double-negate each step of the argument, and then add one application of the stability of equality
at the end.    In fact, this {\em had} to happen:  the G\"odel double-negation interpretation
implies that classical logic can in principle be eliminated from proofs of theorems of the form
  found in Euclid.  (See   \cite{beeson-kobe} and \cite{beeson2015b} for the details using 
  a Hilbert-style axiomatization and  Tarski-style axiomatization, respectively.)

\subsection{Betweenness}  
The betweenness relation $\B(a,b,c)$ means that points $a$, $b$, and $c$ lie on a common 
line $L$, and $b$ is the middle one of the three.   This is taken as a 
primitive (undefined) notion.  Euclid never mentioned it, which was later viewed as a flaw.
Betweenness was used in early (nineteenth-century) axiomatic studies in geometry,
for example \cite{pasch1882,pasch1926}. It was used by 
Hilbert in his famous book \cite{hilbert1899}.  When Tarski developed his theory of geometry,
he used non-strict betweenness, but used the same letter $\B$.  To avoid confusion, we
use $\T(a,b,c)$ for non-strict betweenness.  Either of these two relations can be defined
in terms of the other, even constructively; it really does not matter which is taken as 
primitive. We take $\B$ as primitive.  Then $\T$ can be defined by 
$$ \T(a,b,c) := \neg(a \neq b \land b \neq c \land \neg \B(a,b,c))$$
In the other direction,  $\B(a,b,c)$ can be defined as 
$$\T(a,b,c) \land a \neq b \land a \neq c.$$

We use the equality sign for segment congruence: $ab = cd$.  Hilbert viewed segments
as sets of points, and congruence as a relation between segments.  Tarski viewed 
segment congruence as a relation between four points.  In this paper, we do not 
work with formal theories of geometry in detail, and what we present will work with 
either Hilbert's or Tarski's theories.  Occasionally we may mention ``null segment'':
that refers to a segment $aa$.   Corresponding to the confusion about $\T$ versus $\B$ 
is a confusion about ``closed segments'' versus ``open segments''.   Is a null segment 
a singleton or an empty set?  These technicalities do not arise if we follow Tarski in 
thinking of statements about segments as abbreviations for statements mentioning points only.
Even in informal geometry, we do not think about sets of points.

The fundamental properties of betweenness include symmetry ($\B(a,b,c)$ is equivalent 
to $\B(c,b,a)$) and identity ($\B(a,b,a)$ is impossible.)  In terms of $\T$ this 
becomes $\T(a,b,a) \implies a=b$.  

Betweenness is closely related to order.  
\begin{Definition}[Segment ordering] $ab < cd$ if there exists a point $e$ with 
$\B(c,e,d)$ and $ab = ce$.  Similarly, $ab \le cd$ if there exists a point $e$
with $\T(c,e,d)$ and $ab = ce$.
\end{Definition}
Both Hilbert and Tarski gave axioms for betweenness, and both (but especially Tarski)
tried to be parsimonious about the axioms, i.e., they tried to make the axioms few and simple,
at the cost of making the proofs of relatively simple-sounding theorems fairly difficult.
For example, the ``outer transitivity of betweenness'' is 
$$ \T(a,b,c) \land \T(b,c,d) \land b \neq c \implies \T(a,c,d) \land \T(a,b,d)$$
and the ``inner transitivity of betweenness'' is
$$ \T(a,b,d) \land \T(b,c,d) \implies \T(a,b,c) \land \T(a,c,d).$$
We also have the ``transitivity of betweenness'' (with no adjective):
$$ \T(a,b,c) \land \T(a,c,d) \implies \T(b,c,d) \land \T(a,b,d).$$
One can use $\B$ instead of $\T$, of course.  
These axiomatic exercises do not concern us here; see \cite{schwabhauser}, chapters I.3 
and Satz 5.1 for details.   After having established the fundamental properties 
of betweenness, it is easy to show that segment ordering satisfies
the usual properties for a linear ordering.   See, for example, p.~42 of \cite{schwabhauser}.

The trichotomy law classically takes the form 
$$  \B(a,b,c) \land \B(a,b,d) \implies \T(b,c,d) \lor \T(b,d,c).$$
Constructively this is not valid:  we need to double-negate the right hand side.
One  constructive formulation that does not involve double negation is
$$ \B(a,b,c) \land \B(a,b,d) \land \neg \T(a,c,d) \implies \T(a,d,c).$$
We refer to this principle as ``comparability.''  For a proof from Tarski's
axioms, see \cite{schwabhauser}, Satz 5.2.  In this paper, we view geometry informally,
and take the view that any reasonable axiomatic theory of constructive geometry must 
either assume or prove these principles.  

To say ``$x$ is the midpoint of $ab$'' means $ax=xb$ and $\T(a,x,b)$.  (So the 
concept applies even if $a=b$, because $\T$ is used instead of $\B$.)
We will give an example of a lemma about betweenness. The example will help illustrate
the relations between informal and formal geometry, which is of interest at this point,
but we also need to refer to this particular lemma in the proof of Lemma~\ref{lemma:Euclid5impliesone-sided}.
The point about informal and formal geometry is that, informally, we feel free to use
any of the transitivity and order principles discussed above, whereas in formal geometry,
some of these principles are difficult to prove from a minimal set of axioms.   The 
following proof is thus an informal one,  although it is a rigorous proof from the 
above principles.
\begin{Lemma}\label{lemma:midpoint-helper} Suppose $a$, $b$, and $t$ all lie on 
a line $L$.  
If $x$ is the midpoint of $ab$, and $f$ is 
the midpoint of $at$, and $\B(a,f,x)$, then  $\B(a,t,b)$. 
\end{Lemma}

\noindent{\em Proof}.  From the definition of midpoint we have 
$ax = xb$ and $af=ft$.
We have $af < ax$, since $\B(a,f,x)$ and $af = af$.  
Since $\B(a,f,x)$ and $\B(a,x,b)$, by transitivity we have $\B(a,f,b)$, so 
$a \neq b$ and $ab$ is not a null segment. Similarly $f \neq a$ and $at$ is not 
a null segment.
We have
\begin{eqnarray}
\T(a,x,b) && \mbox{\qquad from the definition of midpoint} \label{eq:1} \\
\T(a,f,x) &&  \mbox{\qquad by hypothesis} \label{eq:2}\\
\T(f,x,b) &&  \mbox{\qquad by \ref{eq:2}, \ref{eq:1}, and transitivity }\label{eq:3} \\
\T(a,f,b)  &&  \mbox{\qquad by \ref{eq:2}, \ref{eq:1}, and transitivity} \label{eq:4} \\
\T(a,f,t) &&\mbox{\qquad from the definition of midpoint} \label{eq:5}
\end{eqnarray}
Hence $bx < bf = fb$, by the 
definition of $<$.  Now, if $\T(f,b,t)$, then $fb < ft$, so we have 
$$ af < ax = bx < fb < ft = af,$$
which is impossible. Hence $\neg \T(f,b,t)$.

By (\ref{eq:4}), (\ref{eq:5}), and comparability, we have
$\T(f,t,b)$.   By (\ref{eq:5}) and outer transitivity, if $f \neq t$ we 
have $\T(a,t,b)$.  Since $at$ is not a null segment, $a \neq t$, and since 
$\neg \T(f,b,t)$, $t \neq b$.  Hence $\B(a,t,b)$.  That completes the proof.

\subsection{Incidence and betweenness}
We use $on(x,L)$ to express the relation that point $x$ lies on line $L$.  This paper
does not focus on the details of formalization of geometry, but rather on constructivity 
in geometrical proofs and constructions; nevertheless some confusion will be avoided by 
the following discussion.  Some formal theories of geometry treat lines as first-order objects
(that was the case for Hilbert, although he then treated line segments and circles as sets of points).
Tarski's geometry, by contrast, officially does not even mention lines.  The ``official''
development of geometry from Tarski's axioms, \cite{schwabhauser}, does treat lines as sets of points,
in what amounts to a conservative extension of the original theory, but the authors
insist that these mentions of lines are just 
abbreviations.  The point we wish to call attention to now is the connection between 
the incidence relation $on(x,L)$ and the betweenness relation $\B(a,b,c)$.  Namely, if 
$L = \Line(a,b)$ is given by two distinct points $a$ and $b$, then $on(x,L)$ is 
classically equivalent to 
$$\B(a,b,x) \lor \B(b,x,a) \lor \B(x,a,b) \lor x=a \lor x=b.$$
Constructively, $on(x,L)$ is equivalent to the double negation:
$$ on(x,L) \leftrightarrow \neg \neg (\B(a,b,x) \lor \B(b,x,a) \lor \B(x,a,b) \lor x=a \lor x=b)$$
and we have the following lemma:

\begin{Lemma}\label{lemma:stability-incidence} [Stability of incidence]
 $\neg \neg on(x,L)$ implies $on(x,L)$.
\end{Lemma}

\noindent{\em Discussion}.  Since we are not committing in this paper to a fixed axiomatization
of constructive geometry, we enumerate some possibilities.  In Tarski's theory, $on(x,L)$
is defined by the formula above, so the stability is immediate.  In a Hilbert-style 
axiomatization, we would have a choice to define incidence or take it as a primitive relation.
If we define it in terms of betweenness, we need to use double negations as above, so 
stability is immediate.  If we take it as primitive, we will want to include stability of 
incidence as an axiom.

\section{Constructions in geometry}

\subsection{The elementary constructions} \label{section:elementary}
The Euclidean constructions are carried out by constructing lines and circles and marking 
certain intersection points as newly constructed points.  Our aim is to give an account of this 
process with modern precision. 
We use a system of terms to denote 
the geometrical constructions.  These terms can sometimes be ``undefined'', e.g. if two lines are parallel,
their intersection point is undefined.    A model of such a theory can be 
regarded  as a many-sorted algebra with partial functions representing the basic
geometric constructions.  Specifically, the sorts are
{\it Point}, {\it Line}, and  {\it Circle}.
We have constants and variables of each sort. 

 It is possible to define extensions of this theory
with definitions of \Arc\ and {\em Segment}.  These are ``conservative extensions'',  which means that 
no additional theorems about lines and points are proved by reasoning about arcs and segments. 
Angles are treated as triples of points.  While this conservative extension result applies to logical 
theories, a similar result applies to constructions considered algebraically.  If we replace rays and segments by lines,
and angles by pairs of lines,  then some terms may become defined that were not defined before, as lines may
intersect where rays or segments did not, etc.   But any point that was constructible with rays and angles will still be constructible
when rays and angles are replaced by lines.  Therefore it suffices to restrict attention to points, lines, and circles, which
we do from now on.%
\footnote{Of course,  when
drawing diagrams to be viewed by people, we use rays and segments.}

Lines are constructed by drawing a line through two distinct points; the resulting line is 
$\Line(a,b)$.  Circles are constructed ``by center and radius'';  $\CircleThree(a,b,c)$ is the circle   
with center $a$ and radius $bc$.   This term corresponds to a ``rigid compass''.   

We have already discussed the discontinuity of Euclid's proof of (the uniform version of) Euclid~I.2.
More generally, any construction $ext(a,b,c,d)$ that extends segment $ab$ by $cd$ will exhibit a 
discontinuity as $a$ tends to $b$ with $cd$ fixed, since $a$ can spiral in towards $b$, 
causing $ext(a,b,c,d)$ to make (approximately) circles of a fixed size around $b$.  That shows that 
we cannot hope to define the extension of a segment $ab$ by $cd$ without a case distinction 
on whether $a=b$ or not.  In constructive geometry, we can extend non-null segments,
not {\em any} segment. 

We can, however, 
always construct a point $e(x)$ different from $x$, and that fact can be used
to give a constructive proof of Euclid~I.2.
\begin{Lemma}\label{lemma:I.2}[Euclid~I.2] Given points $a$, $c$, and $d$, there
exists a point $d$ such that $ad = bc$.
\end{Lemma}

\noindent{\em Remarks}. There is no restriction that $a \neq b$ here;  so Euclid~I.2 is constructively provable 
after all:  it is only Euclid's proof that is non-constructive, not the theorem. 
\smallskip

\noindent{\em Proof}. First we show how, given $x$, to construct a point $e(x)$ 
not equal to $x$.  Fix 
two unequal points $\alpha$ and $\beta$.  
Then define $e(x)$ to be the extension of the non-null 
segment $\alpha\beta$ by $\alpha x$.  It is easy to show that if $e(x) = x$ then the 
segment $\alpha x$ is congruent to its proper subsegment $\beta x$, which is not the 
case.  Now if we extend 
the segment with endpoints $e(a)$ and $a$  by $bc$, we obtain a point $d$ such that $ad = bc$.
That completes the proof.
\smallskip

It follows that $\CircleThree$ can in principle be eliminated in favor of $\Circle$,
by the definition
$$\CircleThree(a,b,c) = \Circle(a,ext(e(a),a,b,c))$$

 We allow circles of radius zero,  which we call ``degenerate circles''.  That is, we consider
 $\CircleThree(a,b,c)$ to be defined when $b=c$.   It is not that we need such circles, but allowing them 
permits the avoidance of case distinctions.

 Starting with at least three noncollinear points, we can  
  construct more points using the following six   ``elementary constructions'', each of which has return type {\it Point} (that is, constructs a point):
\begin{eqnarray*}
&&\IntersectLines(\Line K, \Line L) \\
&&\IntersectLineCircleOne(\Line L, \Circle C) \\
&&\IntersectLineCircleTwo(\Line L, \Circle C)\\
&&\IntersectCirclesOne(\Circle C, \Circle K) \\ 
&&\IntersectCirclesTwo(\Circle C, \Circle K) 
\end{eqnarray*}

$\IntersectCirclesOne$
and $\IntersectCirclesTwo$ construct the intersection points 
of two circles;  how the two points are distinguished will be explained below (in this same section).
Euclid does not use function symbols for the intersection points of two circles, and in particular 
did not worry about how to distinguish one from the other, although sometimes (as in Euclid I.9, 
the angle bisector theorem%
\footnote{See Heath's commentary on I.9 in \cite{euclid1956}.  In order to fix I.9,  Euclid should have 
proved a stronger version of I.1,  constructing on a given segment $AB$ two equilateral triangles, one on 
each side of $AB$.}%
) his proofs need to be repaired by selecting ``the intersection point on the 
opposite side of the line joining the centers as the given point $x$'',  or ``on the same side.''
To directly support this kind of argument we might consider adding function symbols
\begin{eqnarray*}
&&\IntersectCirclesSame(\Circle\, C, \Circle\, K, \Point\, p) \\
&&\IntersectCirclesOpp(\Circle\, C, \Circle\, K, \Point\,  p)
\end{eqnarray*}
These two constructions give the intersection point of $C$ and $K$ that is on the same side (or the opposite side)
of the line 
joining their centers as $p$,  when $p$ is not on that line.  It turns out that constructions meeting 
this specification can be defined in terms of those already mentioned, and hence it is not necessary to include them as primitive symbols.  Space does not permit the inclusion of these constructions here;
but even without seeing them explicitly, we know they must exist by the metamathematical results 
of \cite{beeson-kobe} and \cite{beeson2015b}, according to which, since we can prove the existence
of an intersection point on a given side of a line, there must be a way to construct it.

  Another way in which Euclid avoids the need for function 
symbols to identify the intersection points is the phrase ``the other intersection point.''  Of
course ``the other intersection point'' is just the reflection of the given intersection point
in the line connecting the two centers; but constructively, we need uniform reflection, because
we cannot assume that the given intersection point does not lie on that center line (it will 
in case the circles are tangent).   See
 Theorem \ref{theorem:otherintersectioncircles} for details.

Constructively, every line comes equipped with the two points that were used to construct it.  Informally we simply say,  let $L = \Line(a,b)$  to recover the two points $a$ and $b$ 
that were used to construct $L$.  In a formal theory with variables for points and lines, 
such as was used in \cite{beeson-kobe}, 
we  need   ``accessor functions'' to recover points on those lines.
Thus, $\pointOnOne(L)$ and $\pointOnTwo(L)$ would construct two points on a line;
those will be $a$ and $b$ where $L = \Line(a,b)$.  In this paper, we work informally,
so there is no need to mention accessor functions.  In \cite{beeson2015b}, we use a formal 
theory with only one sort of variables, for points.  That keeps the formal apparatus to a minimum.

Similarly, since circles can only be constructed from a given center,
 every circle comes equipped with its center, so in a formal many-sorted theory,
 we would need an ``accessor function'' 
$center(C)$ for the center of $C$.

 There are three issues to decide:  
\begin{itemize}
  \item  when there are two intersection points,  which one is denoted by which term? 
 \item  When the indicated lines and/or circles do not intersect, what do we do about the term(s) for their intersection point(s)?
   \item  In degenerate situations, such as $\Line(p,p)$, or the intersection points of two coincident circles, what do we do? 
\end{itemize}  
We summarize the answers to these questions, with details to follow:
\begin{itemize}
\item The two intersection points   $\IntersectLineCircleOne(L,C)$ and \goodbreak
\noindent
$\IntersectLineCircleTwo(L,C)$ occur on $L$ in the same order as the two points $p$ and $q$
such that $L= \Line(p,q)$.
\item The two intersection points $p =\IntersectCirclesOne(C,K)$ and \goodbreak\noindent
$q=\IntersectCirclesTwo(C,K)$
are distinguished by requiring that $abp$ is a right turn and $abq$ is a left turn, where $a$ and $b$
are the centers of $C$ and $K$ respectively.
\item  When the indicated lines and/or circles do not intersect, or when two circles
coincide,  the term for their
intersection points will be ``undefined''.
\item  In degenerate 
situations such as $\Line(p,p)$,  we will be guided by continuity.  Thus $\Line(p,p)$ is 
undefined, but $\Circle(a,a)$ is defined.
\end{itemize}
\medskip

That summary leaves several details to spell out.  First, there is the matter of points
occurring ``in the same order'' on a line.  We define $\SameOrder(p,q,s,t)$ to say 
 that $p \neq q$ and $s$ and $t$ both lie on $\Line(p,q)$, and then to be a double-negated
disjunction giving the many possible combinations of betweenness relations between $p$, $q$, $s$, and $t$
that are allowed.  The case $s=t$ is allowed and $\SameOrder(p,q,s,s)$ is considered true.  That is,
$\SameOrder(p,q,s,t)$ says that 
the non-strict order of $s,t$ is the same as the order of $p,q$.  (This allows for the case
when a line is tangent to a circle and the two intersection points coincide.)
 
We note that lines are treated {\em intensionally}.
That is, $\Line(a,b)$ is not equal to $\Line(b,a)$,  even though the two lines are extensionally equal, i.e., the same points
lie on them.  The reason for this is that if the line is rotated about the midpoint of segment $ab$ 
by $180^\circ$, the two intersection points of the line with a circle must change places.  
 That is,  a line ``comes equipped'' with two points that were used to construct it;  to ``be given'' a line 
involves being given two points on that line, and order matters.  It thus makes sense to define the two intersection 
points of $\Line(A,B)$ and a circle to be such that their ordering on the line is the same as that of $A$ and $B$.   If one does not want to think about ``lines with direction'',  then just work with 
a points-only theory in the style of Tarski, as developed in \cite{beeson2015b}.  

Second, there is the question of ``undefined terms''. 
The operations of these ``algebras of constructions'' do not have to be defined 
on all values of their arguments.    The same issue, of course, arises in many other algebraic contexts, for example, division 
is not defined when the denominator is zero, and $\sqrt x$ is not defined (when doing real arithmetic) when $x$ is negative.
The symbol for ``defined'' is $\defined$.   In dealing with undefined terms, we adopt the convention of ``strictness''.
For example, if we write $\B(a,t,b)$ for some term $t$, then that implies $t\defined$.  One cannot make assertions about 
undefined objects.  There are precise rules for a ``logic of partial terms'', which 
can be found in \cite{beeson-book}, and geometric theories incorporating these rules are 
discussed in \cite{beeson-kobe} and the last part of \cite{beeson2015b}.  
  For this paper it does not matter--any one of several ways of handling 
``undefined terms'' will work.

 Third, we mentioned the concepts ``right turn'' and ``left turn'' to distinguish the 
 intersection points of two circles. 
There is a definition of $\Right$ and $\Left$ in 
field theory, well known to all computer graphics programmers:  $abc$ is a right or left turn,
provided $a \neq b$,  
according to the (non-strict) sign of the cross product of the vectors $bc$ and $ba$. (The 
two-dimensional cross product is a scalar, defined by $(a,b) \times (c,d) = ad-bc$.)
 This definition of $\Right$ and $\Left$
can be used within a geometrical theory, as soon as coordinates can be introduced based on 
the geometrical axioms.%
\footnote{The relation of continuity to ruler and compass constructions, and the problem 
of distinguishing the intersection points of circles, have been 
considered before by Kortenkamp \cite{kortenkamp}.  His interest, like mine, originally
arose out of programming dynamic geometry software.  The cited dissertation contains 
many interesting examples, on some of which, for example his Fig.~6.9, I do not agree about the most desirable behavior. Kortenkamp's approach is algebraic and computational 
rather than axiomatic, so is not directly related to this paper.
} 
 We refer the
interested reader to  \cite{beeson2015b} for a full discussion.  We emphasize that the reason
to define $\Right$ and $\Left$ is to be able to give axioms guaranteeing the continuity of 
the terms for intersections of circles, and for intersections of circles and lines.   Without such 
axioms, those function symbols might be interpreted by wildly discontinuous functions. 
Regarding the borderline cases:  we
use the non-strict sign of the cross product, so $abb$ is 
both a left and a right turn if $a \neq b$.  We do not even need to ask the question whether
$aab$ is a right or left turn, since we do not ask for the intersection points of concentric 
circles, but officially $aab$ is neither a left turn nor a right turn.

For the metamathematical results in this paper, the introduction of coordinates and the 
definitions of $\Right$ and $\Left$ are minor issues.  We wish to clarify this point.
What we need for the results of this paper is that every Euclidean field $\F$ provides a 
model $\F^2$ for our geometrical theories.  Which version of the parallel axiom is satisfied 
in $\F^2$ depends on which axioms about multiplicative inverses (or division)
 we take in the axioms of Euclidean 
fields.  Thus our independence results do not 
become harder if we change the (non-parallel) axioms, as long as Euclidean planes satisfy the axioms,
and the axioms hold in planes over Playfair fields; in other words the coordinates of any points
asserted to exist should be explicitly bounded in terms of the given points.  We give two 
different proofs of the independence of Euclid 5 from Playfair, and both have this property.

The ability to define coordinates and arithmetic, on the other hand, says that a geometrical theory 
is ``strong enough''.  This in turn implies that {\em every } model is a plane over a Euclidean field
(at least, if we have line-circle and circle-circle continuity as well as coordinates).  It is 
for this ``representation theorem'' that coordinates and arithmetic are required.

\subsection{Circles and lines eliminable when constructing points}

It is intuitively plausible that one does not really need to draw the lines and circles in a geometric diagram;
one just uses the circles and lines to get their points of intersection with other circles and lines.  
Indeed, one can formulate an axiomatic theory that has only variables for points.  One first 
gives a constructive proof of Euclid~I.2, providing a term
 $E(a,b,c)$ for a point $d$ such that $ad=bc$.  Then one 
 can define $\CircleThree(a,b,c) = \Circle(a,E(a,b,c))$. 
One has,  informally, 
``overloaded'' variants of the function symbols mentioned above:
\begin{eqnarray*}
&& \IntersectLines(\Point a, \Point c, \Point c, \Point d)\\
&&\qquad  = \IntersectLines(\Line(a,b), \Line(c,d)) \\
&& \IntersectLineCircleOne(\Point a, \Point b, \Point c, \Point d)\\
&&\qquad = \IntersectLineCircleOne(\Line(a,b), \Circle(c,d)) \\
&&\IntersectLineCircleOne(\Point a, \Point b, \Circle c) \\
&&\qquad = \IntersectLineCircleOne(\Line(a,b), c) 
\end{eqnarray*}
As these three examples illustrate, 
one can regard circles and lines as mere intermediaries;  points are ultimately constructed
from other points.   Such a points-only axiomatization, in the spirit of Tarski,  is
elegant, and convenient for metamathematical studies.   Below we give points-only formulations
of all three versions of the parallel postulate.  In \cite{beeson2015b} we give a 
precise axiomatization of constructive geometry using a points-only language.

\subsection{Describing constructions by terms} \label{section:midpointscript}
It is customary to describe constructions by a sequence of elementary construction steps. In this section, for 
the benefit of readers not expert in logic, we show how terms in a logical system correspond to traditional
descriptions.   For example, 
we might describe bisecting a segment this way:

\begin{verbatim}
Midpoint(a,b)
{  C = Circle(a,b);  // center at A, passing through B
   K = Circle(b,a);
   P = IntersectCircles1(C,K);
   Q = IntersectCircles2(C,K);
   L = Line(P,Q);   // the perpendicular bisector of AB
   J = Line(a,b);
   m = IntersectLines(L,J);  // the desired midpoint
   return m;
}
\end{verbatim}
This description (which actually is in a precisely-defined language used for computer animation of constructions)
looks quite similar to textbook descriptions of constructions (see for example,
\cite{hartshorne}).  We call such a description a
``construction script'' or just a ``script.''  For theoretical purposes, as in this paper, it is 
easier to describe constructions officially by terms as above.  When the midpoint construction is described as a term, 
it looks like this:

\begin{eqnarray*}
 && \IntersectLines(\Line(\IntersectCirclesOne(\Circle(a,b),\Circle(b,a)), \\
  &&  \IntersectCirclesTwo(\Circle(a,b),\Circle(b,a))), \Line(a,b))
\end{eqnarray*}

Conversely,  given a term,  one can ``unpack'' it, introducing new variable names for subterms,  and even, if 
desired, collapsing duplicate subterms (e.g.,  $\Circle(A,B)$ need not be drawn twice). 
Since we have chosen readable names for our function symbols, such as $\IntersectLines$ instead
of $f$ or even $i \ell$,  our terms sometimes grow typographically cumbersome, and for human readability we will 
often use scripts to describe constructions; but we do not define scripts formally; they always
just stand for terms.

\section{Models of the elementary constructions}\label{section:models}
 We will explain four particularly interesting models of the elementary constructions.
 We assume there are three constants $\alpha$, $\beta$, and $\gamma$ of type {\em Point},
 intended to stand for three non-collinear points.

\subsection{The standard plane}
The {\em standard plane} is $\R^2$.
 Points,  lines, and circles (as well as segments, arcs, triangles, squares, etc., in 
extended algebras in which such objects are considered) 
are interpreted as the objects that usually bear those names in the Euclidean plane.  More formally, the interpretation of the type symbol
``$\Point$'' is the set of points, the interpretation of ``$\Line$'' is the set of lines, etc.  In particular we must choose three specific non-collinear
points to serve as the interpretations of $\alpha$, $\beta$, and $\gamma$.  Let us choose $\alpha = (0,1)$, $\beta = (1,0)$, and 
$\gamma = (0,0)$.   The constructor and accessor functions listed above
also have standard and obvious interpretations.  It is when we come to the five operations for intersecting lines and circles that we must 
be more specific.  As discussed above,  the interpretations of the five function symbols for the elementary constructions, 
 such as $\IntersectLineCircleOne$, are partial functions, so that if the required intersection points do not exist, the term 
 simply has no interpretation.  In degenerate situations, terms are defined if and only if they can be defined so as to make the 
 interpreting function continuous; thus $\Line(p,p)$ is undefined and $\Circle(p,p)$ is the zero-radius circle, and the only 
 point on it is $p$.
 
 To distinguish the intersection points of two circles:
  $\IntersectCirclesOne(C,K)$ is the intersection point $p$ such that the angle from $\CircleCenter(C)$ to $\CircleCenter(K)$
to $p$ makes a ``left turn''.   This is defined as in computer graphics, using the sign of the cross product.  Specifically, let $a$ be the center of circle $C$ and $b$ the center of $K$.
 Then the sign of $(a-b) \times (p-b)$ determines whether angle $abp$ is a ``left turn'' or a 
``right turn''.   These definitions can be given in geometry, as soon as one can 
introduce coordinates and define addition and multiplication.
This notion is constructively appealing, because of continuity:
there exists a unique continuous function of $C$ and $K$ that satisfies the stated
handedness condition for $\IntersectCirclesOne$ when $C$ and $K$ have two distinct intersection points, and is defined whenever $C$ and $K$ intersect (at all).   But it is also interesting,
even with classical logic, to introduce function symbols for the intersection points and then
distinguish the intersection points in a continuous way. 

The principle of continuity leads us to make $\IntersectCirclesOne(C,K)$ and \\ $\IntersectCirclesTwo(C,K)$ undefined in the 
``degenerate situation'' when circles $C$ and $K$
coincide, i.e.,  have the same center and radius.  Otherwise, as the center of $C$ passes through the center of $K$, there is a discontinuity.
It makes sense, anyway, to have them undefined when $C$ and $K$ coincide, as the usual formulas for computing them get zero denominators, 
and there is no natural way to select two of the infinitely many intersection points.

\subsection{The recursive model}
The {\em recursive plane} $\Rrec^2$ consists of points in the plane whose coordinates are given by ``recursive reals'', which 
were defined above.  It is a routine exercise to show that the recursive
points in the plane are closed under the Euclidean constructions.   In particular, given approximations to two circles (or to a circle and a line),
we can compute approximations to their ``intersection points'',  even though it may turn out that when we compute still better approximations
to the circles, we see that they do not intersect at all.   

We showed above that in the recursive plane, there is no computable test-for-equality function.  Recall that the proof 
made use of a recursive real number $E(x)$ such that $E(x) > 0$ if $\{x\}(x)$ halts  and $E(x) = 0$ if $\{x\}(x)$ does not halt.
The same example shows that we cannot test computably for incidence of a point on a line;  the number $E(x)$ lies on the $y$-axis 
if and only if Turing machine $x$ does not halt at $x$.

In the recursive model,  we cannot decide in a finite number of steps whether a circle and a line intersect or whether 
two circles intersect;  if they intersect transversally we will find that out after computing them to a sufficient accuracy, but 
if they are tangential,  we won't know that at any finite approximation.   Technically, the circle of center $(0,1)$ and 
radius $1-E(x)$ will meet the $x$-axis if and only if the $x$-th Turing machine does not halt at $x$.

Readers familiar with recursion theory may realize that there are several ways to define computable functions of real numbers.  The model
we have just described is essentially the plane version of the ``effective operations''.  It is a well-known theorem of Tseitin, Kreisel, 
LaCombe, and Shoenfield, known traditionally as KLS (and easily adapted to the plane) that effective operations are continuous.  Of course,
in the case at hand we can check the continuity of the elementary constructions directly.

\subsection{The algebraic model}
The {\em algebraic plane} consists of points in the plane whose coordinates are algebraic.  Since intersection points of 
circles and lines are given by solutions of algebraic equations, the algebraic plane is also closed under these constructions.  
Note that the elements of this plane are ``given'' all at once.  We assume algebraic numbers are given by means of 
a rational interval $(a,b)$ and a square-free polynomial $f \in \Q[x]$ such that $f$ has only one root in $(a,b)$.

In the algebraic plane, there {\em is} a computable test-for-equality function $D$.  
To determine if $(a,b)$ and $f$ determine the same or a different real number than $(p,q)$ and $g$,  first check if the two rational 
intervals overlap.  If not, the two reals are different.  If so, let $(r,s)$ be the intersection.  Now we have to determine if $f$ and 
$g$ have a common zero on $(r,s)$.   There is a simple recursive algorithm to do that:  Say $g$ has degree greater than or equal to that 
of $f$.  Then write $f = gh + r$ with $r$ of lower degree than $g$.  Then $f$ and $g$ have a common zero on $(r,s)$ if and only if 
$f$ and $r$ have a common zero.  Recurse until both polynomials are linear, when the decision is very easy to make.  
Similarly,  we can compute whether two circles or a circle and a line intersect.

Since algebraic numbers can be computed,  the algebraic model is isomorphic to a submodel of the 
recursive plane.

\subsection{The Tarski model}
The {\em Tarksi model} is $K \times K$,  where $K = \Q({}^{\sqrt{\ }})$ is the least subfield of the reals containing the rationals and closed under taking
the square root of positive elements.  This is a submodel of the algebraic plane.   Its points are the points constructible
with ruler and compass.  Here $\alpha$, $\beta$, and $\gamma$ are interpreted as three fixed rational points.   That 
the intersection points of circles can be computed using only the solution of quadratic equations is checked in detail in section
\ref{section:continuityovereuclideanfields}.

\section{Development of constructive geometry} \label{section:development}
In this section we sketch some important ideas of constructive geometry,
with emphasis on where it differs from classical geometry.  These developments
lay the foundations for the geometric definitions of addition and multiplication,
which are the culmination of the theory, as they show that Euclidean field theory 
can be interpreted in geometry.

\subsection{Pasch} \label{section:pasch}
In 1882, Moritz Pasch axiomatized geometry using betweenness, and pointed out 
the necessity for a new axiom that would guarantee the existence of the intersection
points of certain lines.   This new axiom has several common versions,  so sometimes,
if one version is taken as an axiom, other versions are theorems.  The most 
traditional version is this:

\begin{Theorem}\label{theorem:pasch}
[Pasch 1882] Let line $L$ meet side $ab$ of triangle $abc$, and 
let $L$ not meet side $bc$ or any vertex. If $L$ lies in the plane of 
triangle $abc$, then $L$ meets side $ac$.
\end{Theorem}

Without the codicil that $L$ lie in the plane of the triangle, 
this theorem fails in more than two dimensions.  In that form it says something 
about the dimensionality of space as well as something about betweenness.  
But the form we have given requires the concept of ``plane.''
Therefore two other versions,  which hold  without considerations of dimension,
are also often considered:

\begin{figure} [ht]
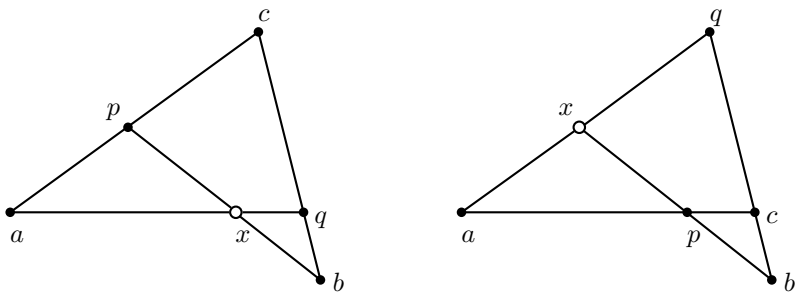

\center{\InnerOuterPaschFigure}
  \caption{Left: inner Pasch.  Right: Outer Pasch. \goodbreak
 Open circles indicate the point 
asserted to exist.}
\end{figure}

Either one of inner Pasch or outer Pasch implies the other, but not easily; and 
either implies Pasch's theorem; and these implications are also valid in 
constructive geometry.   See \cite{schwabhauser}, \cite{beeson2015b},
\cite{tarski-givant} for further discussion.  We will freely use any of these 
principles in this paper.

\subsection{Some Euclidean lemmas}

In this section we collect some lemmas, whose constructive proofs are in Euclid or 
follow easily from Euclid with a bit of attention to constructivity.  
They will be cited in this paper and hence are summarized here without proofs, 
except for a few remarks about the constructive issues. 
They are provable in neutral constructive geometry.
By ``neutral constructive geometry'' we mean geometry with no parallel postulate,
but with line-circle continuity and circle-circle continuity.%
\footnote{The terminology is by no means consistent in the literature; geometers 
do not even agree over whether ``absolute geometry'' means the same thing as ``neutral geometry'',
or whether line-circle continuity is included or not.  Our usage is consistent, 
for example, with \cite{greenberg}.}
 Formal proofs of these theorems
can be given from the axioms listed in \cite{beeson2015b}.
 The existence of constructive proofs of many, but not all, of  these lemmas 
follows from the double-negation interpretation, as shown in \cite{beeson2015b}, 
so it is not really necessary to write out constructive proofs of those lemmas directly.
In this paper, we are not interested in the minimal axiom set required to deduce
such theorems;  any correct axiomatization of constructive geometry should be able
to prove these lemmas.  The only reason to cite \cite{beeson2015b} here is to 
indicate that the details have been checked for at least one formal axiomatization.
The fact that all these theorems can be proved even without line-circle or circle-circle
continuity (due to Gupta, presented in \cite{schwabhauser}, and constructive according to \cite{beeson2015b}) is amazing, but irrelevant here,  since for our purpose of studying constructive geometry
versus classical geometry, we are perfectly willing to use circles.

\begin{Lemma} \label{lemma:euclidI.1}
 Every segment $ab$ with $a \neq b$ has a midpoint (a point $m$ between $a$ and $b$ with $ma = mb$)
and a perpendicular bisector (a line perpendicular to $\Line(a,b)$ at $m$).  Both $m$
and a perpendicular bisector can be constructed (defined) by terms of \ECG.
\end{Lemma}

Note the hypothesis $a\neq b$;  if $a$ and $b$ approach each other, the perpendicular 
bisector of $ab$ is not continuous as $a$ becomes equal to $b$.  See 
Lemma~\ref{lemma:uniformmidpoint}
for a proof that  if $a$ and $b$ 
are restricted to a fixed line $L$,  then one can construct a point $x$ on $L$ 
such that $ax=xb$,  without assuming $a \neq b$.  

\begin{Lemma} \label{lemma:subsegment}  If $\B(a,b,c)$ then $ab$ is not congruent to $ac$.
\end{Lemma}

\begin{Lemma} \label{lemma:exactlytwo}  If line $L$ contains a point strictly inside circle $C$,
then $L$ meets $C$ in exactly two points.
\end{Lemma}

It will be shown in Lemma~\ref{theorem:uniformreflection} that a point can be reflected in a line,
without a case distinction whether the point lies on the line or not.  The following much simpler 
lemma assumes the existence of the reflected point.

 \begin{Lemma}[Reflection in a diameter] \label{lemma:reflectionindiameter}
  Let line $L$ pass through the center of circle $C$
 and let $p$ lie on $C$.   Suppose $q$ is the reflection of $p$ in $L$.  
 Then $q$ lies on $C$.
 \end{Lemma}
 
 \begin{Definition}\label{definition:perpendicular}
Line $K$ is {\bf perpendicular to} line $L$ at $m$ if $m$ is on both $L$ and $K$, and 
there is a point $a$ on $K$ and points
$b$ and $c$ on $L$ with $\B(b,m,c)$ such that $ab=ac$ and $bm = mc$.   $K$ is perpendicular 
to $L$ if there is a point $m$ such that $K$ is perpendicular to $L$ at $m$.
\end{Definition}
One can show that 
if the condition in the definition is true for one point $a$, then it is true for any point $a$ on $K$ not equal to $m$.  It is worth noting that perpendicularity, and hence ``right angle'',
 can be defined without discussing angles in general.

\begin{Lemma}[Stability of perpendicularity]\label{lemma:perpstable}
If two lines are distinct and have a point $m$ in common, and are not not perpendicular at $m$, 
then they are perpendicular at $m$.
\end{Lemma}

\noindent{\em Proof}.  Let $L$ and $K$ be two lines that are not not perpendicular at $m$,
meeting at point $m$.  Let $a$
be a point on $K$ different from $m$ and let $b$ be a point on $L$ different from $m$,
and construct $c$ such that $\B(b,m,c)$ and $bm = mc$.  Then since $L$ and $K$ are 
not not perpendicular, we have $\neg \neg\,(ab=ac)$.  By the stability of equality we have 
$ab=ac$.  Hence $K$ is perpendicular to $L$.

 \begin{Lemma} \label{lemma:perpendicular-symmetric}  If $L$ and $K$ are perpendicular, then $K$ and $L$ are 
perpendicular.
\end{Lemma}

\begin{Lemma} \label{lemma:droppedperpendicularunique}
If $p$ is not on $L$ and $a$ and $b$ are on $L$ and $pa \perp L$ and $pb \perp L$, then 
$a=b$.
\end{Lemma}

\begin{Lemma} \label{lemma:Euclid4} All right angles are congruent.   In other words, 
if $abc$ and $ABC$ are right angles with $ab = AB$ and $bc = BC$ then $ac = AC$.
\end{Lemma}

This is a particularly interesting lemma, since Euclid took it as a postulate; but in 
Hilbert's treatment, and in Tarski's,  it is a theorem, namely Satz 10.12 in \cite{schwabhauser}. 
A constructive proof from Tarski's axioms is given in \cite{beeson2015b}.

The following lemma requires some axiom that implies all the points lie in the same plane,
since it fails in three-dimensional space.   This lemma might even be assumed as the 
upper dimension axiom;  Tarki's geometry has an upper dimension axiom for each $n$, stating 
that space is at most $n$-dimensional. (Only one of these axioms is used at a time.) 
This lemma is equivalent to Tarski's dimension axiom (A9) for the case $n=2$.  
We use the phrase ``in plane geometry'' to mean that some
such axiom is assumed.

 \begin{Lemma} [Uniqueness of erected perpendicular] 
 \label{lemma:uniqueperpendicular}
Suppose $K$ and $M$ are lines perpendicular to another line $L$ at $m$.  Then, in plane geometry,
 $K$ and $M$ coincide (have the same points).   
\end{Lemma}

  \begin{Definition}\label{definition:tangent}
 Line $L$ is {\bf tangent to circle $C$ at $x$} if $L$ meets $C$ at $x$ and only at $x$.
\end{Definition}

\noindent{\em Remark}.  Some books take the property in the following lemma as the definition of {\em tangent} instead.
\medskip

\begin{Lemma} \label{lemma:tangentcircle} (i) If circle $C$ with center $e$ is tangent to $L$ at $x$
then $ex \perp L$. 
\smallskip

 (ii) Given line $L$,  point $x$ on $L$, and two points $a$ and $e$ not on $L$,
with $ex \perp L$ and $e$ lying on a perpendicular bisector of $xa$,  then $\Circle(e,a)$
is tangent to $L$ at $x$.  
\end{Lemma}
 
  \begin{Lemma}\label{lemma:oneparallel}  
 Lines perpendicular to the same line are parallel.
 \end{Lemma}
 
\begin{Definition} [Reflection in a point] The {\bf reflection} of $p$ in $a$ is the \goodbreak
\noindent
(unique) point $t$
such that $\T(p,a,t)$ and $pa = at$.
\end{Definition}

We shall only use reflection in a point when the two points are distinct. In that case it is 
completely unproblematic.  (However, see Corollary~\ref{lemma:uniformreflectioninpoint} for 
a construction that works without a case distinction.)
 The following lemma says that reflection is an isometry:
  
\begin{Lemma} \label{lemma:reflectionpreservescongruence}
Reflection in a point preserves betweenness and congruence.
\end{Lemma}

\begin{Corollary} \label{lemma:chordbisector}
Let $ab$ be a chord of  circle $C$.  Then the center of $C$ lies on the perpendicular bisector of $ab$.
\end{Corollary}

\noindent{\em Remark.}  This is essentially Euclid III.1.
\medskip

\noindent{\em Proof}. Let $m$ be the midpoint of $ab$ and $p$ the center of $C$.  Then $pa = pb$ and $ma = mb$
and $pm = pm$, so triangle $pam$ is congruent to triangle $pbm$.  Hence $mp$ is perpendicular to $ab$ at $m$.
But also the perpendicular bisector of $ab$ is perpendicular to $ab$ at $m$; 
so by Lemma~\ref{lemma:uniqueperpendicular}, $p$ 
is on the perpendicular bisector. That completes the proof. 
 
\begin{Lemma} [Additivity of congruence]  \label{lemma:congruenceadditive}
 Suppose $\B(a,b,c)$ and $\B(A,B,C)$ and $ab=AB$ and $bc = BC$.  Then $ac = AC$.

\end{Lemma} 
 
 \begin{Lemma}[Transitivity of parallel]\label{lemma:transitivityofparallel}
 Suppose line $L$ is parallel to both line $K$ and line $M$, and suppose $M$
 and $K$ do not coincide.  Then, assuming Playfair's axiom,  $M$ and $K$ are parallel.  
 \end{Lemma}

\begin{Definition} \label{definition:rectangle}
 A {\bf rectangle} is a quadrilateral lying in a plane with four right angles.
\end{Definition}

There are related notions, equivalent in Euclidean geometry,
but not in neutral geometry.  

\begin{Definition} \label{definition:saccheri}
A {\bf Saccheri quadrilateral} is a quadrilateral (lying in a plane) with two adjacent right angles
and the opposite sides adjacent to those angles equal.  For example, quadrilateral
$abcd$ is a Saccheri quadrilateral if the angles at $b$ and $c$ are right angles 
and $ab = cd$.   
\end{Definition}

\begin{Definition} \label{definition:parallelogram-rhomboid}
A {\bf rhomboid} is a quadrilateral in which opposite sides are congruent.  (Traditionally
also not all sides are equal, but we do not require that.)   
A {\bf parallelogram} is a quadrilateral with each pair of opposite
sides parallel.  
\end{Definition}

\begin{Lemma} \label{lemma:saccheri}
 In neutral geometry, a Saccheri quadrilateral is a parallelogram.
Playfair's axiom implies that a Saccheri quadrilateral is a rectangle.
\end{Lemma}

\noindent{\em Proof}.  Suppose $abcd$ is a Saccheri quadrilateral.   Let $m$
be the midpoint of the summit $ad$.  It is well known (see e.g. \cite{greenberg},pp. ~177-178,
and the proof there is constructive) 
that the line $K$ joining the midpoint of $bc$ to $m$ is perpendicular
to both the base $bc$ and the summit $ad$.  Hence $ad$ is parallel to $bc$.  
  That completes the proof of the first assertion.

To prove the second claim, we drop a perpendicular $J$ from $m$ to $ab$.  Then $J$ 
is parallel to $bc$, since both are perpendicular to $K$.  By Playfair, $J$
coincides with $ad$, since both are parallel to $bc$ through $m$.   But by 
construction, $J$ is perpendicular to $ab$; since it meets $ab$ at $a$,  the 
angle $bad$ is a right angle.  Similarly, angle $abd$ is a right angle. 
That completes the proof.

\begin{Lemma} \label{lemma:lambert}
Playfair's axiom implies that if a quadrilateral in a plane has three right angles, then 
it is a rectangle.
\end{Lemma}

\begin{Lemma} \label{lemma:midpointprojection} Playfair's axiom implies that orthogonal
projection preserves midpoints.  Here the ``orthogonal projection'' of $p$ onto 
a line $L$ is the foot of the perpendicular from $p$ to $L$.
\end{Lemma}

Assuming the parallel postulate and that space is two-dimensional (for example that
an upper dimension axiom such as Lemma~\ref{lemma:uniqueperpendicular} holds), 
a rhomboid is a parallelogram and vice versa.   Without an upper dimension axiom,
a rhomboid need not lie in a plane, but because ``parallel'' is defined to require
that the lines lie in a plane, a parallelogram must lie in a plane.
Without the parallel postulate, and indeed without an upper dimension axiom,
we can prove the following lemma.

\begin{Lemma}\label{lemma:diagonals}  If the diagonals of a rhomboid meet, then 
the point of intersection bisects both diagonals.   The diagonals of a parallelogram
always meet.  
\end{Lemma}

The hypothesis that a rectangle exists implies Playfair's postulate, so the 
following lemma is still provable in neutral geometry.

\begin{Lemma} \label{lemma:diagonalsofrectangle} The diagonals of a rectangle 
meet and bisect each other.
\end{Lemma}

\subsection{Sides of a line}
The plane separation theorem says that a line separates a plane into two ``sides'' and every
point not on the line is in one side or the other.  This is too strong constructively; it
is essentially the principle $x \neq 0 \implies x < 0 \lor x > 0$.  But we can prove some versions
of the theorem constructively.

Two points $a$ and $b$ not on line $L$ are on opposite sides of $L$ if $a \neq b$ and 
there is a point of $L$ between $a$ and $b$, i.e., the segment $ab$ meets $L$.

\begin{Definition} \label{definition:oppositesides}
\begin{eqnarray*}
\hskip -0.5cm &&\OppositeSide(a,b,L) := \B(a,\IntersectLines(\Line(a,b),L),b) 
\end{eqnarray*}
\end{Definition}

The definition of being on the same side is less straightforward.  One obvious definition of 
$\SameSide(a,b,L)$ would be that segment $ab$ does not meet $L$.  That was Hilbert's definition.
Formally that would be 
$$  \forall x \neg (\B(a,x,b) \land on(x,L)).$$
This proposed definition has two problems.  First,   Hilbert's definition
 works only in the presence of some ``dimension axiom'' that assures that all points lie in a plane.
 The second problem with this definition is the universal quantifier.  In classical 
 geometry we could eliminate it by writing
 $$  a=b \lor \neg \B(a,\IntersectLines(\Line(a,b),L),b)).$$
 But constructively, that implies a case distinction whether $a=b$ or not, which is 
 not possible.  When $a=b$,  the term $\Line(a,b)$ is undefined, so we would not be 
 able to prove $a$ is on the same side of $L$ as itself.  The problem arises only 
 when we try to eliminate the quantifier;  the quantified version is constructively 
 okay.  
 
 Tarski gave a definition that avoids these difficulties: 
$a$ and $b$ are on the same side of $L$ if there is 
some $c$ such that both $a$ and $b$ are on the opposite side of $L$ from $c$.  This 
version has existential instead of universal quantifiers.  Two of those quantifiers 
can be replaced by terms, but the third remains.

\begin{Definition} \label{definition:sameside}
\begin{eqnarray*}
 \SameSide(a,b,L) :=\exists c &&(\B(a,\IntersectLines(\Line(a,c),L,c) \land \\
&& \ \B(b,\IntersectLines(\Line(b,c),L),c)) 
\end{eqnarray*}
\end{Definition}
The point $c$ serves as a ``witness'' that $a$ and $b$ are on the same side of $L$.
We would need $c$ anyway to use Hilbert's definition in 3-space,  where $c$  and the 
two intersection points in Tarski's definition would be
needed as a witness that $a$ and $b$ are in the same plane with $L$.  The definitions
are not so different, after all, and of course they turn out to be equivalent, as
the following theorems show.

\begin{Lemma}[Plane separation theorem] \label{lemma:sameside}
 If $a$ and $b$ are not on $L$ and no point of $L$ is between $a$ and $b$ then $\SameSide(a,b,L)$.
 Indeed a specific point $c$, namely the reflection of $a$ in $L$, is constructible (i.e., given by 
 a term $t(a,L)$)  such that 
 $c$ is on the opposite side of $L$ from both $a$ and $b$.
\end{Lemma}

\noindent{\em Remark}.   This lemma
  obviously fails in three-space, since if $a$ and $b$ are on the same side of $L$ then 
they are coplanar with $L$, but it may be that segment $ab$ does not meet $L$ even though they are not coplanar.
\medskip

\begin{figure}[ht]
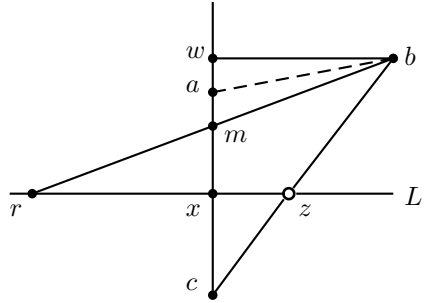
  
\center{  
\PlaneSeparationFigure
}
\caption{If $ab$ does not meet $L$ then $a$ and $b$ are on the same side of $L$ because $bc$ meets $L$. (Lemma~\ref{lemma:sameside}).
\label{figure:PlaneSeparationFigure}
}
\end{figure}

\noindent{\em Proof}. Fig.~\ref{figure:PlaneSeparationFigure} illustrates the proof; the point illustrated
by an open circle must be proved to exist.  Drop a perpendicular from $b$ to $ac$; 
its foot is $w$.  Let $m$ be the midpoint of $wx$.  Reflect $b$ in point $m$ to obtain
$r$.  We do not know yet that $r$ lies on $L$, but we can show triangles $bwm$ and $rxm$
are congruent, so $rx \perp ac$.  Then, by Lemma~\ref{lemma:uniqueperpendicular},
 $r$ must lie on $L$. Line $L$ enters triangle $bmc$ at $x$;  hence by Pasch 
(Theorem~\ref{theorem:pasch}) it must exit the triangle.  By 
hypothesis it does not meet side $ab$.  Hence point $z$, indicated by the diagram in the 
open circle, exists.   Technically, the version of Pasch required here is 
``outer Pasch'', as discussed in \S\ref{section:pasch}.  Hence, the theorem can be proved in any axiomatization of constructive
geometry, as soon as Lemma~\ref{lemma:uniqueperpendicular} and outer Pasch are proved.

\begin{Lemma}[Plane separation theorem, part 2] \label{lemma:sameside2} 
\

If $\SameSide(a,b,L)$ and $\OppositeSide(a,c,L)$ then $\OppositeSide(b,c,L)$.
\end{Lemma}

\noindent{\em Proof}.  This is just a restatement of Pasch's theorem. (We say 
``theorem'' instead of ``axiom'' because Pasch's theorem 
can be derived from special cases, such as ``inner Pasch'' or ``outer Pasch'',
as explained in  \cite{beeson2015b} or \cite{schwabhauser}.)

Hilbert took planes as a primitive concept.  He did not attempt geometry in 
$n$ dimensions.  Tarski instead defined planes; the plane determined by a line $L$
and a point $p$ not on $L$ is the set of points that are either on $L$, or on the
same side of $L$ as $p$, or on the opposite side.   Constructively,  that case 
distinction is unacceptable,  and the quantifier involved in ``same side'' is also 
annoying.  Instead we make the following definition:

\begin{Definition} \label{definition:plane}
  Given line $L$ and point $p$ not on $L$,
let $\bar p$ be the reflection of $p$ in $L$.  Then
a point $x$ {\bf lies in the plane determined by $p$ and $L$} if not not $xp$ meets $L$
or $x\bar p$ meets $L$.
\end{Definition}

Tarski had an ``upper dimension axiom'' (A9), which can be used to specify that
space is at most $n$-dimensional.  The case $n=2$ says that if three points are 
each equidistant from two fixed points $a$ and $b$ then the three points are collinear.
In view of the definition of ``perpendicular'' and theorems about the existence of 
perpendiculars and midpoints, this is equivalent to saying that any two lines perpendicular to 
a given line $L$ coincide.

Here is an application of the plane separation theorem that we need later on:

\begin{Lemma}[Middle Parallel]\label{lemma:middleparallel}
Let $J$, $K$, and $N$ be three weakly parallel lines, all perpendicular to the same line $L$,
and meeting $L$ in points $j$, $k$, and $n$,  in that order (that is, $\T(j,k,n)$); 
and suppose that $j \neq n$.  Let $M$ be another line, meeting $J$ and $N$ in $e$ and 
$f$ respectively.  Then 
$M$ also meets the middle parallel $K$, in a point (non-strictly) between $e$ and $f$.
\end{Lemma}

\begin{figure}[ht]
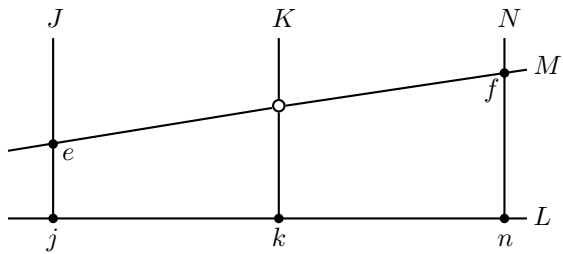

\center{\MiddleParallelFigure}
\caption{The middle parallel meets $M$ if the outer two do}
\label{figure:middleparallel}
\end{figure}

\noindent{\em Proof}. Let $M$ meet $J$ and $N$ at points $e$ and $f$, respectively.
Extend segment $fe$ (by any amount) to a point $e^\prime$.  By the uniform perpendicular
construction, there exists a line $J^\prime$ through $e^\prime$ perpendicular to $L$.
Replacing $J$ by $J^\prime$, we can assume without loss of generality that $J$ and $K$
do not coincide (or in other words, $j \neq k$).  Similarly we may assume $k \neq n$, 
i.e., all three parallels are distinct. 

 Now $fn$ does not meet $K$, so by 
Lemma~\ref{lemma:sameside}, $f$ and $n$ are on the same side of $K$. 
Similarly, $e$ and $j$ are on the same side of $K$.  Since $\B(j,k,n)$, 
$j$ and $n$ are on opposite sides of $K$.  Then by Lemma~\ref{lemma:sameside2},
$e$ and $f$ are on opposite sides of $K$. Hence $ef$ meets $K$, in a point $r$ strictly 
between $e$ and $f$.  We get strict betweenness because we have have assumed $\B(j,k,n)$;
when we allow non-strict betweenness as in the first paragraph, we get only $\T(e,r,f)$,
as stated in the conclusion.  That completes the proof
of the lemma.

  \subsection{Angles} \label{section:angles}
 Euclid took angles as a primitive concept, and assumed that angles can be compared.
 That is, he assumes an angle ordering relation $abc < uvw$. 
 (There is no need to use the notation $\angle abc < \angle uvw$
 since $abc < uvw$ could not have any other meaning.)
 Often Euclid proves two angles are equal by an argument beginning, ``if not, then one of them is 
greater''  and proceeding to a contradiction.  If angles could be compared by 
comparing segments, this would be justifiable in the same way that we justify proving equality
of betweenness of points by contradiction, i.e., by stability of equality, congruence, and betweeness.

 Hilbert, following Euclid, treated angles as a primitive concept.  Tarski
 treated angles as triples of points.  In \cite{beeson2015b}, 
we follow Tarski in treating angles as triples of points.  In this paper, 
in order not to commit to a specific formal treatment, and especially to avoid
a long development of the properties of angle ordering, we just avoid angles
entirely in the formal statements of the various versions of the parallel axioms.  However,
we indicate how the theory of angles can be developed if angles are taken as
a defined concept, following \cite{schwabhauser}.  In order to verify the 
constructivity of Euclid, we need to justify  

 \begin{Definition} [Ordering of angles] \label{definition:angleordering}
 Point $p$ is in  angle $abc$ if $p \neq b$ and there exists a point $x$
 on $\Line(b,p)$ such that $\T(a,x,c)$.  If instead of $\T(a,x,c)$ we use $\B(a,x,c)$,
 then $p$ is said to be in the interior of angle $abc$. 
 Angle $abc <  ABD$ if 
 there is a point $C$ in the interior of angle $ABD$ such that angle $abc$ 
 is congruent to angle $ABC$.   Similarly, $abc \le ABD$ if there is a point 
 $C$ in angle $ABD$ such that angle $abc$ is congruent to angle $ABC$.
 \end{Definition}
 
The theory of angle congruence and angle ordering is
developed in \cite{schwabhauser}, Chapter 11.  To get an idea of what is 
involved, you can try the following exercise:  prove that if $p$ is in the interior of $abc$, 
and $\B(b,a,A)$ then $p$ is also in the interior of $Abc$.  
In making the cited development constructive, 
there are some technical problems, 
turning on the exact version of 
Pasch's axiom that is needed.  These problems are discussed
and solved
in \cite{beeson2015b}.  The following theorem justifies the constructivity
of Euclid's use of angle ordering:

\begin{Theorem}\label{theorem:anglestability}
 [Stability of angle ordering] \label{lemma:stabilityofangleordering}
$$ \neg \neg\, abc < def \implies abc < def $$
and similarly with $\le$ in place of $<$.
\end{Theorem}

\noindent{\em Remark}.  The proof is given 
in \cite{beeson2015b}.
\medskip

It is often instructive to see what a theorem about angles ``really states'' 
when angles are eliminated.  The following famous proposition I.16 of Euclid
is a good example.
 
 \begin{Lemma}[Exterior angle theorem, Euclid I.16] \label{lemma:exteriorangle}
 Suppose $\B(b,c,d)$ and $a$ does not lie on $\Line(b,c)$.
Then  $acd >bac$.  
\end{Lemma}

\begin{figure}[ht]
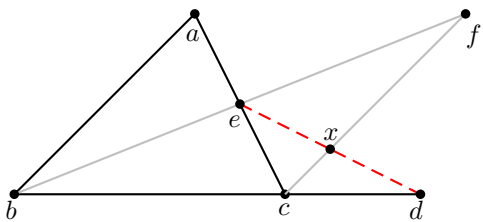

\center{\OneSixteenFigure}
\caption{Construct $L$ parallel to $ab$.  Then $e$ and $x$ exist by Pasch.}
\label{figure:I.16}
\end{figure}

\noindent{\em Proof}.   Euclid's diagram is given by the solid lines 
in Fig.~\ref{figure:I.16}. 
  We show how to complete Euclid's proof.
To prove $acd > bac$, we must construct a point $f$ in the interior of $acd$
such that $acf = bac$.  Euclid knew what angle ordering means:  he
 constructs $f$.  But he does not prove that 
$f$ lies in the interior of $acd$.  We supply the proof:  By  inner Pasch,
applied to the five-point configuration $becdc$, there is a point $x$ such that both
$\B(c,x,f)$ and $\B(e,x,d)$.  Hence $f$ is in the interior of angle $ecd$.
But $ecd$ is the same angle as $acd$, so $f$ is in the interior of angle $acd$ as well,
according to the exercise mentioned above.
That completes the proof.

The immediate corollary (Euclid I.17) is that any two angles of a triangle, taken together,
are less than two right angles.  In particular, no triangle contains two right angles,
although in neutral geometry, nothing prevents the angle sum of (all three angles of)
 a triangle from being 
more than two right angles.   

\begin{Lemma} \label{lemma:legsmallerhypotenuse}
In a right triangle, the hypotenuse is greater than either leg.
\end{Lemma}

\noindent{\em Proof}.  First prove Euclid I.18 and I.19.
Then apply them as indicated in Exercise~22, p.~198 of \cite{greenberg}.

Euclid's Prop.~I.18 says that ``In any triangle the greater side subtends the greater 
angle.''  Euclid's proof from I.16 is perfectly constructive.  The following lemma then
 follows quickly from his Postulate~4 (all right angles are equal);  
 
\begin{Lemma} \label{lemma:rightanglecontains}  An angle contained in a right angle
is not a right angle.  That is, if 
  $abc$ is a right angle, and $p$ and $q$ both lie between $a$ and $c$, 
then $pbq$ is not a right angle.
\end{Lemma}

\noindent{\em Proof}. By Euclid~I.18, $abc$ is greater than $abq$, which in turn is 
greater than $pbq$.  We need the transitivity of angle ordering to conclude $abc > pbq$,
and by Lemma~\ref{lemma:Euclid4} (which is Euclid 4), if $pbq$ is a right angle, then 
$abc = pbq$, contradiction.  That completes the proof of the lemma, modulo the basic
properties of angle ordering, which are developed in \cite{schwabhauser}, Chapter 11,
with attention to constructivity in \cite{beeson2015b}.

\subsection{Uniform perpendicular }
In this section we show that some fundamental constructions needed to define (signed) arithmetic can be defined
from the elementary constructions.

A very fundamental fact in constructive geometry is the existence of a 
 {\bf uniform perpendicular} $\Perp(x,L)$,
which is a line perpendicular to line $L$ passing through $x$, constructed
without a case distinction as to 
whether $x$ lies on $L$ or not.  In classical treatments of geometry,  this case distinction is made
and a different construction is used for each case.  

There are several ways to construct a uniform perpendicular.  In this paper we 
give a construction that does not require the parallel postulate, but it does use
circles (and line-circle continuity).  In \cite{beeson2015b},  we give a different
construction that avoids the use of circles entirely, but it does require the parallel axiom.
At present it is an open problem whether a uniform perpendicular can be constructed
without using either circles or the parallel axiom.  

\begin{figure}[h]
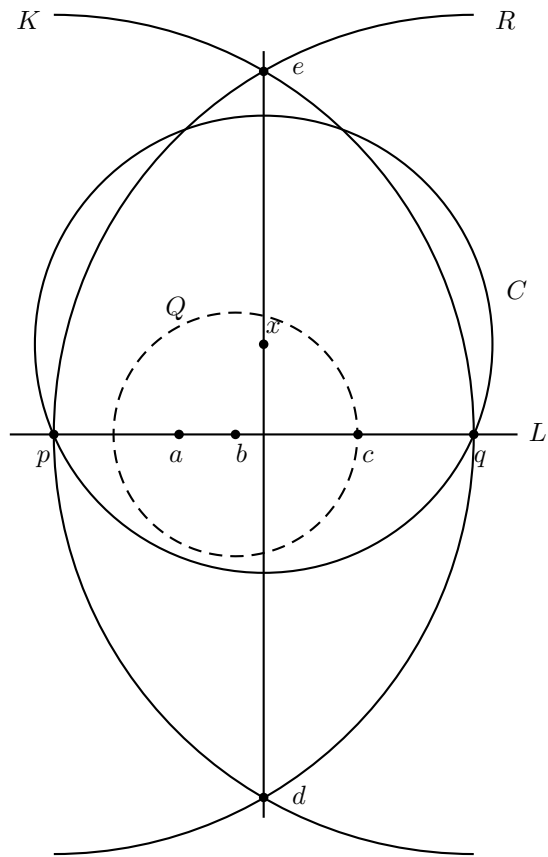
  
\center{
\PerpFigure
\caption{ $M=\Perp(x,L))$ is constructed perpendicular to $L$ without a case distinction whether $x$ is on $L$ or not.  Note $bc=xa$ so the radius $ac$ of $C$ is long enough to meet $L$ twice.\label{figure:PerpFigure}}
}
\end{figure}

Here is the construction that uses circles:

\begin{Definition} \label{definition:perpdef}
The term $\Perp(x,L)$, or expressed in terms of points only, $\Perp(x,a,b)$,
is given by the following construction script, where $L = \Line(a,b)$. (See Fig.~\ref{figure:PerpFigure}.)
\begin{verbatim}
Line Perp(Point x, Point a, Point b)
{ Q = Circle3(b,x,a)
  c = IntersectLineCircle2(L,Q)
  C = Circle3(x,a,c)
  p = IntersectLineCircle1(L,C)
  q = IntersectLineCircle2(L,C)
  K = Circle(p,q)
  R = Circle(q,p)
  d = IntersectCircles1(K,R)
  e = IntersectCircles2(K,R)
  return Line(d,e)
}
\end{verbatim}
\end{Definition}

\begin{Lemma} [Uniform construction of perpendiculars]  \label{lemma:uniformperpendicular}
  Let $\Perp(x,L)$ be defined as in Definition~\ref{definition:perpdef}.
For every point $x$ and line $L$, $\Perp(x,L)$ is defined, and
$\Perp(x,L)$ is perpendicular to $L$, and $x$ lies on $\Perp(x,L)$.
\end{Lemma} 

\noindent{\em Remark}.  The important point here is that no case distinction is required as to whether 
$x$ is or is not on $L$; that is, there is no distinction in this construction between dropped perpendiculars
and erected perpendiculars.  That is good, since constructively we cannot decide if a given point $x$ is 
or is not on a given line, yet we still need to be able to project $x$ onto the line.
\medskip

\noindent{\em Proof}.  We plan to
construct two points $p$ and $q$ on line $L$ such that the perpendicular bisector of $pq$
will be the desired line;  thus we must have $p \neq q$ even if $x$ is on $L$, and whether 
or not $x=a$ or $x=b$.  We can get such points $p$ and $q$ if we can draw a  
  circle with center $x$,  whose radius $r$ is ``large enough.''   It will suffice if we 
take $r$ to be the length of $ax$ plus the length of $ab$.  Since $a \neq b$, that will 
be a nonzero radius, and since it is more than $ax$, the circle with that radius will 
meet $L$ twice, in points $p$ and $q$.  Then we just bisect segment $pq$ by an ordinary 
method,  drawing circles with centers $p$ and $q$ that pass through each other's centers. 

Formally, we go through the script defining $\Perp(x,L)$ and show
that the term on the right of each line is defined.  
  In line 1, 
$\CircleThree(b,x,a)$ is defined because $\CircleThree$ is everywhere defined.
Since the center of that circle is $b$ and $b$ lies on $L$, $b$ is inside that circle, so 
by the line-circle continuity axiom, $c$ in line 4 is defined.  Circle $C$ in line 5
is defined since $\CircleThree$ is everywhere defined.   To show that $p$ and $q$ in lines 5 and 6 are defined and distinct,
 we must 
show that $L$ has a point strictly inside $C$.  That point is $a$, as we shall prove.  
According to the definition of $\IntersectLineCircleTwo$, 
 the intersection points of $L$ and $Q$ occur in the same order on $L= \Line(a,b)$
as $a$ and $b$.  That implies that $c$ occurs on the opposite side of $b$ from $a$,
or in other words,  $\B(a,b,c)$.  By the definition of $<$ we then have $bc < ac$.
Since $c$ lies on $Q$ by construction, we have $bc = ax$.  Hence $ax < ac$.  
Then since the radius of $C$ is $ac$, and its center is $x$, 
 $a$ is strictly inside $C$, as claimed. 

 Hence points $p$ 
and $q$ in lines 4 and 5 of the script   are defined.  By Lemma~\ref{lemma:exactlytwo},
 $p \neq q$.
Then $K$ and $R$ in the next lines of the script are 
distinct circles, so $d$ and $e$ are defined.
 $\Line(d,e)$ is the
perpendicular bisector of segment $pq$, which exists by  Lemma~\ref{lemma:euclidI.1};  the proof of that lemma 
provides the explicit term mentioned in the script for constructing $\Line(d,e)$. 
That completes the proof of the lemma.

\subsection{Uniform projection and parallel}
  The construction $\project(p,L)$ takes a point $p$ and line $L$ and produces a point $q$ on $L$ such that 
$p$ lies on the perpendicular to $L$ at $q$.   The well-known Euclidean construction for the projection applies only if $p$ is known 
not to be on $L$.  To define $\project$ using that construction, we would require a test-for-incidence that allows us to test whether 
point $p$ is on line $L$ or not.  However, we can define uniform projection using the 
uniform perpendicular.  Namely, the projection of $p$ on $L$ is given by 
a ruler-and-compass construction term $\project(p,L)$, giving the foot of the uniform 
perpendicular.

\begin{Lemma} \label{lemma:uniformparallelinformal}
 We can define a construction $\para$ 
such that, for any line $L$ and any point $p$ (which may or may not be on $L$),  $\para(p,L)$ passes through $p$, and if $p$ is not on $L$ 
then $\para(p,L)$ is parallel to $L$, while if $p$ is on $L$, then $\para(p,L)$ has the same points as $L$.
\end{Lemma}

\noindent{\em Proof}. The definition of $\para$ is 
$$ \para(p,L) =  \Perp(p,\Perp(p,L)). $$
In words:  First find the perpendicular to $L$ passing through $P$.  Then erect the perpendicular to that line at $P$.
It is a theorem of neutral geometry that two lines with a common 
perpendicular are parallel.  Since $\Perp$ is everywhere defined, regardless of whether $p$ is or is not on $L$,
the same is true of $\para$.  That completes the proof of the lemma.
\medskip

The following gives an explicit construction for extending a non-null segment on a given end by a (possibly null) segment:
\begin{Definition} \label{definition:Extend}
$$  \Extend(a,b,c,d) = \IntersectLineCircleTwo(\Line(a,b),\CircleThree(b,c,d))$$
\end{Definition}

\begin{Lemma}[Segment extension] \label{lemma:extension}
Let $x = \Extend(a,b,c,d)$.  Then 
$$\B(a,b,x) \  \mbox{ and } \  bx=cd.$$
\end{Lemma}

\noindent{\em Proof}.
Let $y$ be the other intersection point,
namely $$y = \IntersectLineCircleOne(\Line(a,b),\CircleThree(b,d,c).$$ 
Then $b$ is between the two intersection points $x$ and $y$.  Since 
the two intersection points occur on $\Line(a,b)$ in the same 
order as $a$ and $b$, $y$ is on the same side of $b$ as $a$.  Then $x$ 
is on the other side of $b$ from $a$, that is, $\B(a,b,x)$, which was
to be shown.  The congruence $bx = cd$ is part of the line-circle continuity axiom.
That completes the proof.
\medskip

\noindent{\em Remark}. In Tarski's axiom system, there is an axiom guaranteeing
segment extension.  Line-circle continuity is regarded as a more complex axiom,
only to be added when absolutely needed.   

\subsection{Midpoints}
By definition, $m$ is the midpoint of $ab$ if $\T(a,m,b)$ and $am=mb$.  The case $a=b$ is 
allowed in the definition.  The question arises 
whether for every $a,b$ there exists a midpoint $m$ of $a,b$.  When $a \neq b$, we can 
use the  Euclidean construction (repairing Euclid's proof using Pasch's axiom); and when 
$a=b$ we can take $m=a$; but constructively, we are not allowed to conclude that 
for every $a,b$ there exists a midpoint of $ab$.  Indeed, the theorem is not 
constructively valid, as we will now show.  Let $p$ and $q$ be two fixed distinct points.
 Suppose $m(a,b)$ gives the midpoint of 
$a,b$ for $a \neq b$.  Then, for $a \neq b$, we can erect the perpendicular $L(a,b)$ 
to $ab$ at $m(a,b)$, and define $v(a,b) = \IntersectLineCircleOne(L(a,b),\CircleThree(m(a,b),p,q)$.
Then when $a$ approaches $b$ by spiraling in to $b$,  the point $v(a,b)$ does not converge
to a limit.  Hence $m$ cannot be extended to a continuous function defined when $a=b$.
But any ruler-and-compass construction is continuous.  Hence there is no uniform 
construction of the midpoint of $ab$, defined whether or not $a=b$.

We can avoid the problem of the spiraling approach by restricting $a$ and $b$ to lie 
on a fixed line $L$.  Given $a$ and $b$ on $L$, we wish to construct the midpoint of $ab$.
 Our plan is to enlarge the (possibly null) 
``segment'' $ab$ by extending it the same amount in both directions on $L$, and then
use the usual midpoint construction to bisect the resulting non-null segment.  The problem 
is that we do not know the order of $a$ and $b$ on line $L$.  The solution to this problem is,
intuitively, to move $b$ to the right by $ab$ plus one, and move $a$ to the left the same amount;
after that, $ab$ is a non-null segment with the same midpoint as the original $ab$.

\begin{Lemma} \label{lemma:uniformmidpoint} [Uniform midpoint]
There is a construction $m(a,b,L)$ such that for all $a$ and 
$b$ on line $L$,  $m(a,b,L)$ is a point $x$ on $L$ such that $ax = xb$.  
\end{Lemma}

\noindent{\em Remark}.  Then $\Perp(x,L)$ is the perpendicular bisector of $ab$ when $a \neq b$, and 
is perpendicular to $L$ in any case.  
\medskip

\noindent{\em Proof}.  Since every line is given by two points, we have $L = \Line(p,q)$
for some   points $p$ and $q$ with $p \neq q$.  Think of $p$ as left of $q$, and think of 
the length of $pq$ as ``one''.

Here is a construction script for the midpoint of $ab$, given that $a$ and $b$ are on $\Line(p,q)$:
\medskip

\begin{verbatim}
Point m(Point a,Point b, Point p, Point q){
L = Line(p,q)
u = IntersectLineCircle2(L,Circle3(b,p,q))
B = IntersectLineCircle2(L,Circle3(u,a,b))
v = IntersectLineCircle1(L,Circle3(a,p,q))
A = IntersectLineCircle1(L,Circle3(v,a,b))
x = Midpoint(A,B)
return x
}
\end{verbatim}
\medskip

Recall the predicate $\SameOrder(p,q,s,t)$ from Section~\ref{section:elementary},
which says that $p,q$ occur on line $L$ in the same order as $s,t$ do.  Let us say
``$s$ is left of $t$'' to mean $\SameOrder(p,q,s,t)$.  Let us say ``$s$ is to the left
of $t$ by more than $ef$'' to mean that $s$ is to the left of $t$ and there is a proper
subsegment of $st$ congruent to $ef$.
Because of the rules for distinguishing the two intersection points of a line and circle,
 $u$ is to the right of $b$ and $ub =pq$ and $B$ is to the right of and $Bu=ab$,
 so $B$ is to the right of $b$ by more than $ab$.  Similarly,
  $v$ is to the left of $a$ and $av = pq$ and $A$ is to the left of $v$ and  $Av=ab$,
  so $A$ is to the left of $a$ by more than $ab$. 
We claim $A \neq B$.  Suppose $A=B$.  Then $a$ is to the right of $B$ by more than $ab$,
and since $B$ is to the right of $b$ by more than $ab$, it follows that $a$ is to the right 
of $b$ by more than twice $ab$, which is impossible.  Hence, as claimed, $A \neq B$.
Then the Euclidean construction of the midpoint $x$ of $AB$ in the last line is defined.
 
We must show $ax=bx$. 
By the stability of betweenness and equality, we may argue by cases according as $a=b$, 
or $a$ is to the left of $b$ or to the right.  We omit the details.
That completes the proof.   

  \section{The parallel postulate} \label{section:parallel}
We will introduce three versions of the parallel postulate:  Euclid's own version (Euclid 5),
which says that under certain conditions, two lines will meet in a point;  a version we call the 
strong parallel postulate,  which weakens Euclid's conditions for the two lines to meet;  and 
Playfair's version, which says two parallels to the same line through the same point must coincide
(and makes no assertion at all about the existence of an intersection point).  We then 
start to consider the relations of implication between these versions 
(relative to neutral constructive geometry)
and finally draw some consequences of the strong parallel postulate.  
Since our independence results will be obtained using the theory of Euclidean fields, 
the exact axiomatization of neutral constructive geometry will not be important in this paper; but 
  a specific axiomatization is given in  \cite{beeson2015b}.

\subsection{Alternate interior angles}
In connection with parallel lines, the following terminology is traditional:  If $L$ and $K$
are two lines, and $p$ lies on $L$ and $q$ lies on $K$ but not on $L$, then $pq$   
 is a {\it transversal}
of $L$ and $K$.   Then $\Line(p,q)$ makes four angles with $L$ and four angles with $K$; certain 
pairs of them are called {\em alternate interior angles} or {\em  interior angles on the same 
side of the transversal}. 
In the formal theories of constructive geometry given in \cite{beeson-kobe, beeson2015b},
  there is no primitive concept of ``angle'';  angles are treated as triples of points,
and angle congruence is a defined concept; that is also true of Tarski's classical theories
of geometry. In this paper, we use angles informally, but we avoid them 
in the statements of theorems.
 In the next section we will discuss Euclid's parallel 
postulate, and we will need to express it without mentioning angles directly.  We therefore 
discuss how to do that now.
Angle congruence is
defined in terms of triangle congruence, which is in turned defined by the SSS criterion: the 
three sides are pairwise congruent. 

\begin{figure}[ht]
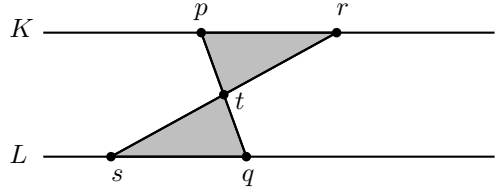
 
\center{
\AlternateInteriorAnglesFigure
}
\caption{Transversal $pq$ makes alternate interior angles equal with $L$ and $K$,  if $pt=tq$ and 
$rt=st$.}
\label{figure:AlternateInteriorAnglesFigure}
\end{figure}

We will now unwind these definitions to express the 
concept of alternate interior angles being equal, without mentioning angles.  
The following definition is  illustrated in 
Fig.~\ref{figure:AlternateInteriorAnglesFigure}.  
    
\FloatBarrier

\begin{Definition} \label{definition:alternateinteriorangles}
 ``{\bf $pq$  makes alternate interior angles equal with $K$ and $L$}'',
written $AIE(p,q,L,K)$, means
\begin{eqnarray*}
&&  \on(p,K) \land \on(q,L) \land \neg \on (q,K) \land p \neq q   \\
&& \land\ \exists r,s,t\, (on(r,K) \land \on(s,L) \land \B(r,t,s) \land \B(p,t,q) \land r \neq p \land q \neq s \\
&&    \land\ pt = qt \land rt = st)
\end{eqnarray*}
This can also be expressed as ``{\bf $K$ and $L$ make alternate interior angles equal with $pq$}''.
\end{Definition}

\begin{figure}[h]
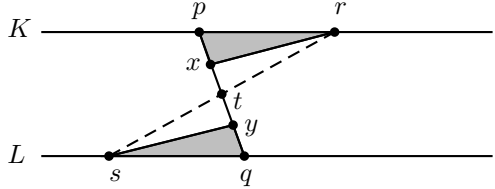
 
\center{\AlternateInteriorAnglesSevenFigure}  
\caption{Another way to say that transversal $pq$ makes alternate interior angles equal with $L$ and $K$.     The shaded triangles are congruent.}
\label{figure:AlternateInteriorAnglesSevenFigure}
 \end{figure}

The reason for requiring $\B(r,t,s)$ is that $r$ and $s$ need to be on opposite sides of $\Line(p,q)$
in order that the angles $rpt$ and $sqt$ are {\em alternate} interior angles.
At the cost of introducing two more points into the definition, we could have made it 
easier to verify;  any two congruent triangles with those angles as a vertex could replace
$rpt$ and $sqt$,  as illustrated in Fig.~\ref{figure:AlternateInteriorAnglesSevenFigure}.
As the theory of angle congruence in \cite{schwabhauser} (Chapter 11) shows
(and one can easily check informally),
if we can produce $x$ and $y$ so that Fig.~\ref{figure:AlternateInteriorAnglesSevenFigure} is correct,
then $AIE(p,q,L,K)$ will hold.

A careful reader might have noticed that we omitted the congruence $rp=qs$ in the definition, 
which would be required to express the congruence of the shaded triangles according to the 
definition of triangle congruence.  This was not an accident.  In case the lines $K$ and $L$
do not coincide (as shown in the diagram), 
the triangles are congruent anyway by SAS, since the angles at $t$
are vertical angles.   The condition is meant to be used only when lines $K$ and $L$ do 
not coincide, and the condition $\neg on(q,L)$ guarantees that it fails when they do coincide.

\smallskip    

On the other hand, instead of introducing more points using $\exists$, we can get rid of $\exists$
entirely by constructing specific points $r$, $s$, and $t$, at least when the two lines 
do not coincide.
\begin{eqnarray*}
r &=& \IntersectLineCircleTwo(K,\Circle(p,\alpha,\beta)) \\
t &=& \Midpoint(p,q) \\
s &=& \IntersectLineCircleTwo(\Line(r,t),\Circle(t,r)) 
\end{eqnarray*}
These points will do for $r$, $t$, and $s$ in $AIE(p,q,L,K)$ if any points will do.
 $\Midpoint(p,q)$ is defined in the script because $L$ and $K$ do not coincide.
In this way $AIE(p,q,L,K)$ can be expressed by a quantifier-free formula.
\medskip

That   alternate interior angles   have something to do with Euclid's parallel postulate is hinted at
by the following lemma:

\begin{Lemma} \label{lemma:alternateinterioranglesimpliesparallel}
Let $K$ and $L$ be two lines, with point $p$ on $K$ but not on $L$,
 and let traversal $pq$ make equal alternate interior angles
with $K$ and $L$.  Then $K$ and $L$ are parallel.
\end{Lemma}

\noindent{\em Proof}.  Refer to Fig.~\ref{figure:AlternateInteriorAnglesFigure} for an illustration;
let points $p$, $r$, $s$, $q$, and $t$ be as shown in the figure with triangle $prt$ 
congruent to triangle $qst$.  Suppose, for proof by contradiction, that $K$ and $L$ meet in point $x$.
(Since $\Parallel$ is defined as not meeting, a constructive proof of it consists in deriving 
a contradiction.)   Consider triangle $pqx$.  If $x$ is on $\Ray(p,r)$, then the angles 
of triangle $pqx$ are $pxq$, $qpx$, and $pqx$. The two interior angles on the same side of $pq$ 
make two right angles together,   but in neutral geometry one can prove that two angles of 
a triangle are less than two right angles.  (This is Euclid I.17, proved in two lines from 
Euclid I.16, which is Lemma~\ref{lemma:exteriorangle}.)
This contradiction shows that $K$ and $L$ cannot meet in a point $x$ on $\Ray(p,r)$.
Similarly, they cannot meet in a point $x$ on $\Ray(r,p)$.  But, if they meet in a point 
$x$ at all, then  
$$\neg \neg\, (\B(x,p,r) \lor x=p \lor \B(p,x,r) \lor x=r \lor \B(p,r,x)).$$
In each of those cases, $x$ lies on $\Ray(r,p)$ or on $\Ray(p,r)$;  so we have a 
contradiction, and $K$ cannot meet $L$ at all.  Hence $K$ is parallel to $L$ as claimed.
That completes the proof of the lemma.

\subsection{Euclid's parallel postulate}

Euclid's postulate 5 is 
\begin{quote} 
{\em If a straight line falling on two straight lines make the interior angles on the 
same side less than two right angles, the two straight lines, if produced indefinitely, meet on that side on which are 
the angles less than the two right angles.}  
\end{quote}

\begin{figure}[h]
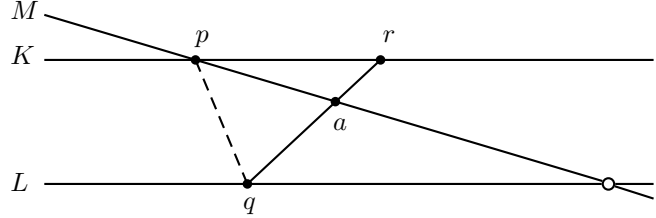
 
\center{\EuclidParallelFigure}  
\caption{Euclid 5:  $M$ and $L$ must meet on the right side, provided $\B(q,a,r)$ and $pq$ makes
alternate interior angles equal with $K$ and $L$.}
\label{figure:EuclidParallelFigure}
\end{figure}

We consider the formal expression of Euclid's parallel axiom.  Let $L$ and $M$ be two straight 
lines, and let $pq$ be the ``straight line falling on'' $L$ and $M$, with $p$ on $M$ and $q$ 
on $L$.  We think that what Euclid meant by ``makes the interior angles on the same side 
less than two right angles'' was that, if $K$ is another line $K$ through $p$,
making the interior angles with $pq$ equal to two right angles, then $M$ would lie in 
the interior of one of those interior angles (see Fig.~\ref{figure:EuclidParallelFigure}).   For $K$
to make the interior angles on the same side of $pq$ equal to two right angles, is the same 
as for $K$ to make the alternate interior angles equal.  In the previous section we discussed
the formal expression of that concept.   
 Once we use this criterion
to express that $L$ and $K$ make interior angles on the same side equal to two right angles,
we can then use three more points to ``witness'' that one ray of line $M$ emanating from $p$
lies in the interior of one of the interior angles made by $M$.   
Fig.~\ref{figure:EuclidParallelFigure}
illustrates the situation.
 
\smallskip    

Here is a formal version of Euclid's parallel axiom, using the formula $AIE$ from the previous section 
to express that $pq$ makes alternate interior angles equal with $K$ and $L$:
\medskip

\axioms
$  on(p,K) \land on(p,M) \land on(a,M) \land on(r,K) \land on(q,L) \land \B(q,a,r)$&    \\
$  \qquad \land \neg on(p,L) \land\ AIE(p,q,K,L) \implies \B(p,a,\IntersectLines(L,M))$ & 
\endaxioms
\smallskip

\noindent
The logical axioms of LPT make it superfluous to state in the conclusion that $\IntersectLines(L,M)$ is defined.  That follows
automatically.  
\smallskip

\begin{figure}[h]
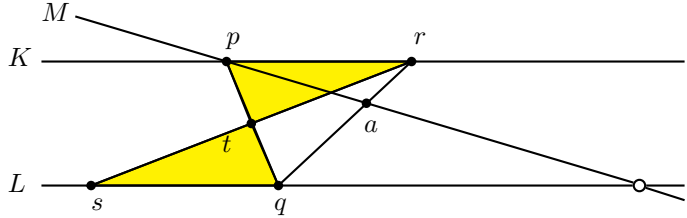
   
\center{\EuclidParallelRawFigure}
\caption{Euclid 5:  $M$ and $L$ must meet on the right side, provided $\B(q,a,r)$ and $pt=qt$ and $rt=st$.}
\label{figure:EuclidParallelRawFigure}
 
\end{figure}

We now write out Euclid 5 in a more primitive syntax, eliminating
the defined concept $AIE$.
The result, illustrated in Fig.~\ref{figure:EuclidParallelRawFigure},
is the official version of Euclid 5;  that is, when we refer to Euclid 5 in subsequent 
sections, this is what we mean.

\smallskip    
 
\medskip

\axioms
$   \B(q,a,r) \land \B(p,t,q) \land \B(s,t,r)  \land pt=qt  \land p \neq r  $&\hskip-1.2cm (Euclid 5)\\
$ \land\, rt=st  \land L = \Line(s,q) \land M = \Line(p,a) \land \neg\, \on(p,L) $ & \\
  $ \implies \B(p,a,\IntersectLines(L,M)) \land \B(s,q,\IntersectLines(\Line(L,M)))$&
\endaxioms
\smallskip

\subsection{The rigorous use of Euclid 5}
Each of  Hilbert and Tarski replaced Euclid's parallel postulate with some different
parallel postulate.  Why did they do that?  We can only guess at the motivations
of Hilbert and Tarski, but one reason may have been the imprecision of Euclid's
version.  Euclid 
mentions angles ``on the same side'' of the transversal,  without ever defining ``side''.

 We choose to adhere to Euclid's formulation (after improving its precision),
 in order to investigate the 
constructive relations between the various forms of the parallel postulate.
Our version of Euclid 5 is valid in 3-space,  since the ``witnessing points'' it requires 
will ensure that the three lines lie in the same plane.   Euclid's own formulation either 
fails or is ambiguous in 3-space.  Two skew lines can make interior angles equal, so unless
``same side'' implies ``in the same plane'',  Euclid 5 will fail.   In addition to this problem,
it can be surprisingly tricky to use Euclid 5 rigorously.  In this section we 
illustrate how that is done. 

% {\em Remark.} We have studied several more variants of Euclid 5. For example, the congruent triangles
% could on the other side of $pq$ and $a$ could be on the other side of $K$.   That 
% change then allows us to relax the hypothesis that $p$ is not on $L$, allowing $K$ to 
% possibly coincide with $L$,  as long as some hypothesis guarantees that $a$ is strictly 
% not between $K$ and $L$.  These versions can all be proved from Euclid 5 as formulated here,
% but we omit this material here.

The following lemma
furnishes an example of 
the use of our version of Euclid 5.  The statement of the lemma is simple 
and fundamental;  its proof is surprisingly intricate and technical.   Part of the point
of exhibiting this proof is to show how much Euclid left unsaid;  these difficulties will arise
with any attempt to formalize Euclid using modern logical languages.   It is worth noting 
that neither Hilbert nor Tarski in published work discussed the propositions of Euclid
(as opposed to the fundamental notions and postulates): congruence 
of angles comes only in chapter~11 (out of 16) of \cite{schwabhauser}, Part I, and circles
are never mentioned.  In particular the difficulties in this section have nothing to do 
with constructivity; they arise from making precise the many things Euclid left out.

\begin{Lemma} \label{lemma:euclidperp}  Let $M$ and $L$ be lines in the same plane,
and let point $p$ be on $M$ but not on $L$.  Suppose $ap \perp M$ and $pf \perp L$,
with $a$ and $f$ on $L$ and $a \neq f$.  Then there exists a point $e$ on $M$ with $\B(e,f,a)$.
\end{Lemma}

\begin{figure}[ht]
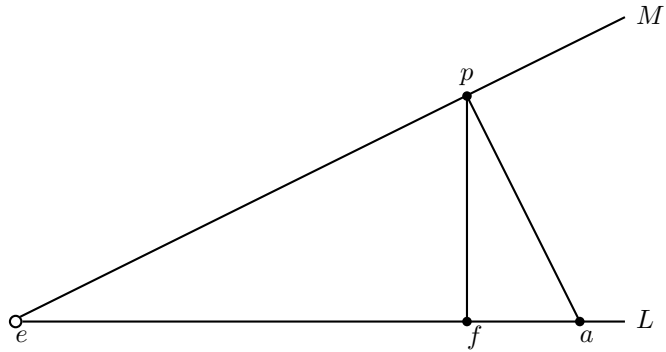

\center{\EuclidPerpFigure}
\caption{Given $p$, $f$, $a$, $M$, and $L$, with $M \perp pa$, to construct $e$ by Euclid 5.}
\label{figure:euclidperpfigure}
\end{figure}

\noindent{\em Remarks}.   The rest of the hypotheses do not 
guarantee that $M$ and $L$ lie in the same plane, and if they are not, 
then the lemma is false.  So the hypothesis that mentions ``plane'' is 
really necessary.   We first give a proof in Euclid's style, and 
then give a more rigorous proof.  We would like to emphasize that the 
difficulties in making this proof rigorous have absolutely nothing to do 
with constructivity.  It is simpler:  Euclid is not rigorous, and neither 
Hilbert nor Tarski has explicitly shown how to work rigorously with Euclid 5.
Each of them replaced Euclid's parallel postulate with a different parallel postulate
instead.
\medskip

\noindent{\em ``Proof'' in Euclid's style}.    The interior 
angles made by traversal $pa$ of $M$ and $L$ are together less than two right angles,
because the angle at $p$ is a right angle, and the angle at $a$ is less than a right 
angle, since {\em pfa} is a triangle with a right angle at $p$.  Therefore $M$
meets $L$ at a point $e$ ``on the side of $pq$ on which angle {\em fap} lies''.  Presumably
that means $\B(e,f,a)$,  although since Euclid did not define ``side'',  so we have to 
say ``presumably.''  That completes the ``proof.''   We never used the hypothesis
that $M$ and $L$ lie in the same plane,  but the theorem fails without it,  so probably 
it is an unstated hypothesis of Euclid 5.
\medskip

\begin{figure}[ht]
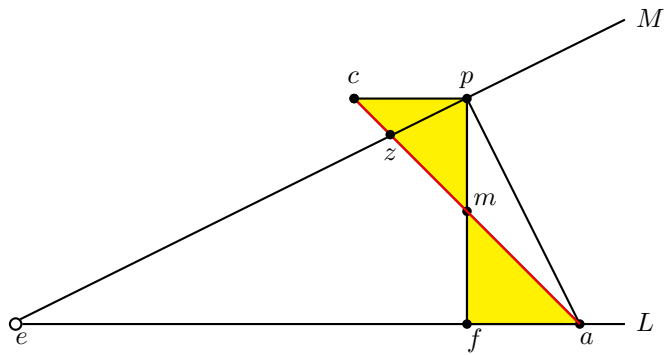

\center{\EuclidPerpFigureTwo}
\caption{Given $p$, $f$, $a$, $M$, and $L$, to construct $e$ by Euclid 5,
first construct $m$ and $c$.  Then 
we need to construct point $z$ using the hypotheses that $M$ lies in a plane with $L$ and $p$ and $M \perp pa$.}
\label{figure:euclidperpfigure2}
\end{figure}

\noindent{\em Proof}.
Now, we give a rigorous proof from our version of Euclid 5.  Let $m$ be the midpoint
of $pf$, and let $c$ be the reflection of $a$ in $m$, so $\B(c,m,a)$ and $cm = am$, 
as shown in Fig.~\ref{figure:euclidperpfigure2}.  
To apply Euclid 5,
we need a point $z$ in that figure; indeed as soon as we have 
a point $z$ on $M$ between $c$ and $a$,  Euclid 5 gives us the existence of 
the desired intersection point $e$ of $M$ and $L$.   Clearly the existence of point $z$
relies on the hypothesis that $M$ and $L$ lie in the same plane.   All that 
remains to fix the Euclidean proof is to prove the existence of $z$.

\begin{figure}[ht]
\center{\EuclidRawPerpFigure}
\caption{Construct $m$, $c$, and $q$.  Then {\em pqaf} is a rectangle, whose diagonals intersect in $x$. Then $b$ exists by inner Pasch, and $z$ by Pasch's theorem applied to triangle $cbq$.}
\label{figure:euclidrawperpfigure}
\end{figure}

\noindent
To construct $z$ turns out to be not quite trivial;  please
refer to Fig.~\ref{figure:euclidrawperpfigure}.
Let $q$ be the reflection of $c$ in $p$.
Then {\em qpfa} has right angles
at $p$ and $f$, and $pq = fa$ since both are equal to $cp$. 
We claim it is a Saccheri quadrilateral; it only remains to prove that it 
lies in a plane.  $L$ and $p$ lie in the plane {\em pfa}, and the points $c$ and $q$ 
constructed by reflection then lie in that same plane.  Then {\em qpfa} is a 
Saccheri quadrilateral, as claimed. 
 By Lemma~\ref{lemma:saccheri},
{\em qpfa} is a parallelogram.  Hence, by Lemma~\ref{lemma:diagonals},
the diagonals intersect.   Let the intersection point be $x$. Then $\B(f,x,q)$ and $\B(p,x,a)$.
Apply inner Pasch to the five-point configuration {\em axpfm}; we obtain point $b$ such 
that $\B(f,b,x)$ and $\B(m,b,a)$.

Now we claim that $M$ does not meet $bq$.  
Suppose that $M$ does meet $bq$ in a point $r$. 
Then angle $xpr$ is a right angle, since $x$ lies on $pa$ and $r$ lies on $M$, and $pa \perp M$.
But this right angle is properly contained in the right angle $fpq$, since both $x$ and $r$ 
lie between $f$ and $q$.  That 
contradicts Lemma~\ref{lemma:rightanglecontains}.   Hence (as claimed)
 $M$ does not meet $bq$.  
 
 $M$ meets  side $cq$ of triangle $cqb$ at $p$, and does not meet side $bq$.
 Moreover $c$ does not lie on $M$, since if it did, $pa$ and {\em pf} would both 
 be perpendicular to $M$ at $p$, contradicting Lemma~\ref{lemma:uniqueperpendicular}.
 Also $a$ does not lie on $M$, since $ap \perp M$. 
Now we finally use the hypothesis 
that $M$ and $L$ lie in the 
same plane: we want to apply Pasch's theorem to line $M$ and triangle $cqb$.
We can do that because $M$ lies in the plane of that triangle. 
By Pasch's
 theorem (Theorem~\ref{theorem:pasch}), $M$ meets $cb$ in a point $z$.
 That is the point we were trying to construct. 

  Now the hypotheses of Euclid 5 are 
fulfilled, with $z$ and $a$ substituted for $a$ and $q$ in Euclid 5.  Hence
the intersection point $e$ of $M$ and $L$ exists, and satisfies $\B(e,f,a)$ 
as claimed.  That completes the proof of the lemma. 
\medskip

The following lemma furnishes another example of the rigorous 
use of Euclid 5.  This lemma will be used in the proof that 
Euclid 5 implies a form of triangle circumscription, which is crucial
for verifying the definability of multiplication.    The proof of the lemma is 
more difficult than it at first appears, and illustrates once again the necessity 
of constructing auxiliary points to justify statements about angles that Euclid 
would have taken for granted. These difficulties have nothing to do with constructivity
but simply arise from treating Euclid with modern rigor.

\begin{Lemma}\label{lemma:euclidexample}
 Let points $a$ and $x$ lie on line $L$, and let
 point $p$ not lie on line $L$.  Let $K$ be perpendicular to $L$ 
at $x$ and $M$ perpendicular to $ap$ at $p$.  Let $J$ be the perpendicular
to $L$ through $p$, meeting $L$ at $f$. Suppose $f$ and $x$ lie on the same ray of $L$
emanating from $f$; that is, $\neg \B(a,f,x)$ and for some point $\ell$ on $L$, we have 
$\ell \neq f$ and $\T(\ell,f,a)$ and $\T(\ell,f,x)$. 
  Then (assuming Euclid 5)  $M$ and $K$ meet.  
\end{Lemma}

\begin{figure}[ht]
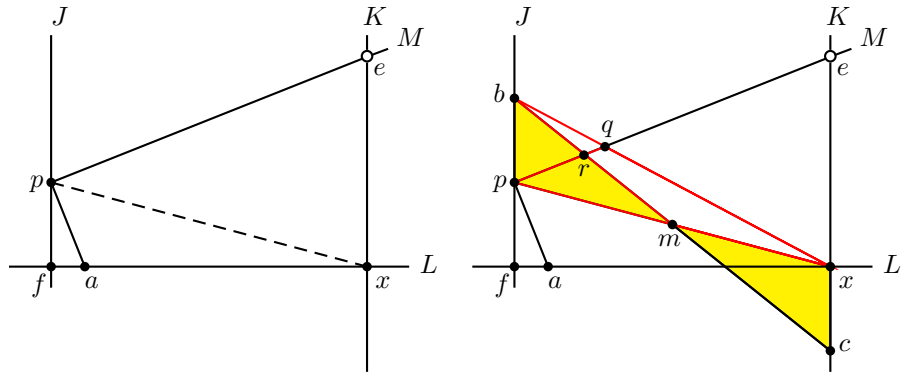

\EuclidExampleFigureOne
\EuclidExampleFigureTwo
\caption{Euclid 5 implies $M$ meets $K$,
whether or not $x=a$. 
Left, attempted traditional proof.  Right, rigorous proof requires point $r$ to exist.
It is found by inner Pasch after constructing $q$ by the plane separation theorem.
}
\label{figure:euclidexample}
\end{figure}

\noindent{\em Remark}.  Thus $a$ can possibly lie on $J$ (in which case $M$ is 
parallel to $L$). Nothing is said 
about the order of $a$ and $x$ on $L$, except that they do not lie on opposite sides of $J$.
The point $\ell$ specifies on which side of $f$ the points $a$ and $x$ lie; it is necessary
since we might have $a=x=f$.  
Fig.~\ref{figure:euclidexample} illustrates the case when $a$ is to the left of $x$.  Point
$\ell$ is not shown; it would be to the left of $f$.
\medskip

\medskip

\noindent{\em Proof}. Refer to Fig.~\ref{figure:euclidexample}.
The figure shows $a$ to the left of $x$, but that is not guaranteed by the 
hypotheses.  We claim that we can assume without loss of generality that $\T(f,a,x)$.
To do so, we replace line $K$ by another line $K^\prime$ perpendicular to $L$ and 
meeting $L$ at a point $x^\prime$ to the right of both $a$ and $x$.  Now, if we prove
the lemma under the extra assumption that $x$ lies to the right of $a$, we can apply
it to show that $M$ meets $L^\prime$ in some point $e^\prime$.  Then $J$, $K$, and $K^\prime$
are three parallel lines, all perpendicular to $L$, and $M$ meets the outer two, 
so according to Lemma~\ref{lemma:middleparallel}, $M$ and $K$ meet, as desired.  That
shows that we may assume $\T(f,a,x)$, as illustrated.

 We wish to apply Euclid 5 to show that $M$ meets $K$,
using the transversal $px$, as shown in the left part of the figure.
  Intuitively, angle $Kxp$ is less than right angle $Kxa$,
and angle $Mpx$ is less than or equal to right angle $Mpa$, so the corresponding interior angles made by 
transversal $px$ are less than two right angles.   

That ``proof'' is insufficiently rigorous.  We must use the
official, points-only version of Euclid 5, which makes no mention of angles.
That requires the construction of a point to witness the fact that traversal $px$
makes corresponding angles with $M$ and $K$ less than two right angles.  The construction
of that point is not accomplished in one step;  
we need to appeal to several theorems, including the plane separation theorem
and the exterior angle theorem. 
The details follow.

First, we claim $M$ does not meet segment $fx$.  Suppose that $M$ meets $fx$ in a 
point $y$.  Then $y \neq p$, since $p$ does not lie on $L$.  Moreover $y \neq a$,
since $pa \perp M$.  Therefore $pay$ is a triangle.
Triangle $pay$ has a right angle at $p$, since $y$ lies on $M$, and $ap \perp M$.
By Lemma~\ref{lemma:exteriorangle}, the remaining angles of $pay$ are each less than a right angle.
Since we are trying to reach a contradiction, we can argue by cases.   If $a=f$ then 
triangle $pay$ has two right angles, contradiction.  Hence $a \neq f$.  If $\B(a,y,x)$
then angle $pay$ is an exterior angle of triangle {\em pfa}, and hence angle $pay$ is greater
than the right angle {\em pfa}, contradiction.  If $\B(f,y,a)$ then angle $pya$ is an exterior 
angle of triangle {\em pfy}, which has a right angle at $f$,  so angle $pya$ is more than a 
right angle, contradiction.  But if $M$ meets $fx$ then,  classically, one of these cases
must hold.  Since all of them are contradictory, we have reached a contradiction, even
constructively, by the stability of betweenness.   We have therefore proved that $M$
does not meet  $fx$.  

By Lemma~\ref{lemma:sameside},  $f$ and $x$ are on the same side of $M$.  Let $b$ be
the reflection of $f$ in $p$.  Then $b$ is on the opposite side of $M$ from $f$.  
By Lemma~\ref{lemma:sameside2}, $b$ is on the opposite side of $M$ from $x$.  Hence 
$bx$ meets $M$ in a point $q$.   Let $m$ be the midpoint of $px$, and let $c$ 
be the reflection of $b$ in $m$.  We claim that $c$ lies on $K$.  We have $bm = mc$
and $pm = mx$, by construction of $c$ and $m$ respectively.  Angles $bmp$ and $xmc$ are
vertical angles, so by SAS, triangles $bpm$ and $cxm$ are congruent.  Hence $cx$ 
and $J$ make alternate interior angles equal with traversal $px$.  Hence $cx$
is parallel to $J$.  Now $cx$ and $K$ are both parallel to $J$ through $x$.  By 
Playfair (and hence by Euclid 5),  $\Line(c,x)$ and $K$ coincide.  Hence $c$ lies 
on $K$.

Now we apply inner Pasch to the five-point configuration $bqxpm$.  The result is a 
point $r$ such that $\B(b,r,m)$ and $\B(p,r,q)$.   This point $r$ is the required
witness that lines $M$ and $K$ make corresponding angles less than two right angles with 
traversal $px$.  That is, we can now apply Euclid 5  with the variables $(L,a,p,r,q,s)$ in 
Euclid 5 replaced by $(K,r,p,b,x,c)$ here.   The conclusion is that $M$ meets $K$
in a point $e$.   That completes the proof.

\subsection{Playfair's axiom}
 Euclid
did not give his Postulate 5 the name ``parallel postulate'' (or any other name).   
A case can be made that it is more of a ``triangle construction postulate.''  Be that 
as it may, we now
  define $\Parallel(L,K)$ for lines $L$ and $K$ to 
mean that the lines do not meet:  
\begin{Definition} \label{definition:parallel}  In plane geometry:
$$ \Parallel(L,K)\quad  \iff \quad \forall x \neg\,(on(x,L) \land on(x,K)).$$   
\end{Definition}

\noindent
This definition assumes we are working in plane geometry;  in an axiomatic setting there
should be some ``dimension axiom'' present, or else the definition needs to require 
as well that $L$ and $K$ lie in the same plane.  That is, skew lines in space are not 
considered parallel.

To express that two lines do not meet requires a universal quantifier over points.
   Compare the given definition to the  formula
$$  \neg\,\IntersectLines(K,L) \defined.$$
If $L$ and $K$ are parallel, then $\IntersectLines(K,L)$ is undefined;
but it will also be undefined if $L$ and $K$ coincide.  
The relation defined by $\neg \IntersectLines(L,K) \defined$ can be referred to 
informally as ``weakly parallel'';  classically it means ``parallel or coincident'',
but constructively, we generally cannot make the decision of two weakly parallel lines 
whether they are coincident or parallel.   Instead, weakly parallel means ``not not 
(parallel or coincident).''  (Consider, for example, the $x$-axis and the line $y=a$;
these lines are parallel if $a \neq 0$ and coincident if $y=0$,  so if we could decide 
whether two weakly parallel lines are coincident or parallel, we could decide $y=0 \lor y \neq 0$.)
The concept ``weakly parallel'' also requires a universal quantifier if expressed in a 
language that does not use the logic of partial terms.

Most modern treatments of geometry formulate the parallel axiom in a different way:
\smallskip

{\em   If two lines $K$ and $M$ are 
parallel to $L$ through point $p$,  then $K=M$.}
\smallskip

 We call this the ``Playfair's axiom'', or for short just ``Playfair'', 
after John Playfair, who published it in 1795, although (according to Greenberg \cite{greenberg}, p. 19) it was referred to by Proclus.  It is also known as ``Hilbert's parallel axiom'', since Hilbert used 
it in his axiomatization of geometry \cite{hilbert1899}.
Since we can always construct one parallel, it would be equivalent to assert that, 
given a line $L$ and a point $p$ not on $L$, there exists exactly one parallel to $L$ through $p$.
\medskip

\begin{figure}[h]
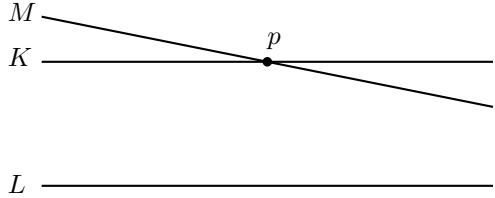
    
 \center{\PlayfairFigure}
\caption{Playfair:  $M$ and $L$ can't fail to meet.}
 \label{figure:PlayfairFigure}
\end{figure}

 \smallskip

The formal expression of Playfair's postulate is as follows:
\medskip

\axioms
$ \Parallel(K,L) \land \Parallel(M,L)$ &(Playfair)\\
\qquad $ \land on(p,K) \land on(p,M)  \implies on(q,K) \iff on(q,M)$&
\endaxioms
\smallskip

There is an ambiguity in this formulation,  concerning $Parallel$.
 If we intend to use this axiom 
only in the presence of an upper dimension axiom, then we can take $Parallel$ to 
mean ``not meeting''.  But if we wish to use Playfair in the absence of a dimension axiom,
then $Parallel$ must mean ``not meeting, and lying in the same plane.''  In that case,
$Parallel$ becomes an existential statement, requiring the existence of five witnessing 
points to the coplanarity of the lines.   But since $\Parallel$ occurs only on the left
of Playfair's axiom,  Playfair's axiom itself is universal, which for an axiom is the 
same as quantifier-free.  In particular, 
Playfair's axiom makes no existential assertion at all.
It is an empirical fact (and formally, 
a consequence of the double-negation interpretation in \cite{beeson2015b}) 
that Playfair's axiom suffices to prove all the consequences 
of the parallel postulate that do not themselves make existential assertions.
 
\begin{Lemma} \label{lemma:playfairimpliesAIE}
  Playfair's axiom implies (in neutral constructive geometry) that  parallel lines make equal 
alternate interior angles with any transversal,  and equal corresponding angles.
\end{Lemma}

\noindent{\em Proof}. Let $L$ and $K$ be  parallel lines, and let $pq$ be a transversal, 
with $p$ on $K$ and $q$ on $L$ and $p \neq q$.  
   Let $s$ be any point on $L$ not equal to $q$.
Let $t$ be the midpoint of $pq$, and extend segment $st$ to point $r$ with $\B(s,t,r)$
and $tr=ts$.   
Let $M = \Line(p,r)$.  Then by Definition \ref{definition:alternateinteriorangles},
$pq$ makes equal alternate interior angles   with $L$ and $M$.
 
If $p$ does not lie on $L$, then by 
Lemma~\ref{lemma:alternateinterioranglesimpliesparallel}, $M$ and $L$ are parallel.
But then $M$ and $K$ are two  lines parallel   to $L$ passing through $p$.  Hence by 
Playfair's axiom, they coincide.   Hence $pq$ makes alternate interior angles with 
$L$ and $K$.   On the other hand, if $p$ does lie on $L$, then $M$ and $K$ coincide 
and we are finished in that case too.  But by the stability of the incidence relation, 
we have $\neg \neg \, (\on(p,L) \lor \neg \on(p,L))$,  so it is constructively legal
to argue by these two cases.   The equality of corresponding angles follows
from the equality of vertical angles.   That completes the proof.

\begin{Lemma} \label{lemma:parallel-perp}
Let $M$ and $K$ be weakly parallel lines, and let $L$ be perpendicular to $M$
and meet $K$.  
Assuming Playfair's axiom,  $L$ is also perpendicular to $K$.
\end{Lemma}

\noindent{\em Proof}.  By the stability of ``perpendicular'' 
(Lemma~\ref{lemma:perpstable}), we may argue by cases according 
as $M$ and $K$ coincide or not.  If they coincide there is nothing to prove.
If $M$ and $K$ do not coincide, then they are parallel.  Let $L$
be perpendicular to $M$ at $p$ and meet $K$ at $q$.  The perpendicular to $L$
at $q$ is parallel to $M$, so by Playfair's axiom, it coincides with $K$.
Therefore $K$ is perpendicular to $L$.  That completes the proof. 

\begin{Lemma} \label{lemma:parallel-transitive} Assuming Playfair's
axiom, weak parallelism is 
transitive.  That is, if $L$ and $M$ are weakly parallel and $M$ and $K$
are weakly parallel, then $L$ and $K$ are weakly parallel.
\end{Lemma}

\noindent{\em Proof}. Suppose $L$ and $M$ are weakly parallel and $M$ and $K$
are weakly parallel, and suppose that 
$\IntersectLines(L,K)\defined$.  We must derive a contradiction to conclude that 
$L$ and $K$ are weakly parallel.  
Let $p$ be the intersection point of $L$ and $K$.
By Lemma~\ref{lemma:stability-incidence}, we may argue by cases to prove a contradiction.

Case 1,  $p$ lies on $M$.  Then $M$ and $L$ coincide, since $M$ and $L$ are weakly parallel;
in that case $L$ and $K$ are weakly parallel since by hypothesis $M$ and $K$ are 
weakly parallel.  That disposes of Case 1.

Case 2, $p$ does not lie on $M$.  Then $M$ and $L$ are two parallels to $K$ through 
$p$, so by Playfair, $M$ and $L$ coincide.  Then $L$ and $K$ are weakly parallel,
 since by hypothesis $M$ and $K$ are weakly parallel.   That completes the proof.
 
\begin{Lemma} \label{lemma:two-perps}  Let lines $M$ and $K$ be weakly parallel,
and let point $p$ be given.  Let lines $U$ and $V$ containing $p$ be perpendicular to $M$ and $K$
respectively.  Then, assuming Playfair,  lines $U$ and $V$ coincide. 
\end{Lemma}

\noindent{\em Proof}. 
Suppose first that $U$ meets $K$. Then by Lemma~\ref{lemma:parallel-perp}, 
$U\perp K$, and hence $U$ and $V$ coincide, as both are perpendiculars to $K$ through $p$.
Similarly, if $V$ meets $M$,  $U$ and $V$ coincide. 

Now we will prove that $U$ and $V$ coincide.  Let $q$ be a point on $U$; we must 
show $q$ is on $V$. By Lemma~\ref{lemma:stability-incidence},  we may prove 
that by contradiction. Suppose, for proof by contradiction, that $q$ is not on $V$.
Then, as shown in the previous paragraph, $U$ does not meet $K$.  

That $U$ does not meet $K$ means that 
$U$ is  
parallel to $K$; since  $K$ is weakly parallel to $M$ by hypothesis,   
  Lemma~\ref{lemma:parallel-transitive} lets us conclude (using Playfair) 
  that   $U$ is weakly
parallel to $M$.  But $U$ meets $M$ and does not coincide with $M$, since $U \perp M$.
That contradiction completes the proof.

\begin{Lemma}[Euclid I.32] \label{lemma:euclidI.32}
 Playfair's  axiom implies that an exterior angle of a triangle is equal to the sum of 
the interior and opposite angles.
\end{Lemma}

\noindent{\em Remark}.  The theorem asserts that the 
exterior angle can be decomposed into two angles, each congruent to one of the opposite interior angles.
Addition of angles has not been defined (either here or in Euclid).  
\medskip

\noindent{\em Proof}.  We have already discussed the closely related Euclid~I.16 in 
Lemma~\ref{lemma:exteriorangle}. Please refer to Fig.~\ref{figure:I.16}.
Euclid's construction provides that $ae=ec$ and $be=ef$, so triangles 
$abe$ and $cfe$ are congruent, making $bac$ and $acf$ alternate interior angles,
and $cf$ parallel to $ba$.  
 
Playfair's axiom implies that parallel lines make equal corresponding angles with 
any transversal, by Lemma~\ref{lemma:playfairimpliesAIE}.
 Hence angle $abc$ is congruent to angle 
$fcd$.  But now, the exterior angle $bcd$ is seen to be composed of two angles, each
 congruent to one of the interior and opposite angles.  That completes the proof.
\smallskip

\noindent{\em Remark}. Euclid~III.20, that an angle inscribed in a semicircle is a right angle, 
follows from I.32 in \ECG\ just as it does in Euclid. 

\subsection{Euclid 5 implies Playfair's axiom}

\begin{Theorem}\label{theorem:Euclid5impliesPlayfair}
Euclid's Postulate 5 implies Playfair's axiom in neutral constructive geometry.
\end{Theorem}

\noindent{\em Proof}.
 Let $M$ and $K$ be 
two parallels to $L$ passing through $p$; we must show that $M$ and $K$ coincide.  Without loss
of generality we can assume that $K$ is perpendicular to the perpendicular $H$ to $L$ passing through $p$.
 Then $H$ makes right angles
with $K$ and with $L$.  Since $M$ does not meet $L$, then by Euclid 5, $M$ does not make less than a right 
angle on either side of $H$.  Thus we have two angles, together equal to two right angles, and 
neither is less than a right angle.  Therefore neither is greater than a right angle; hence each is 
a right angle, and $M$ is perpendicular to $H$ at $p$.  By the uniqueness of the erected 
perpendicular
(Lemma~\ref{lemma:uniqueperpendicular}), 
$M$ coincides with $K$ as desired.  
\medskip

\noindent{\em Remark.}  We defined lines to be parallel if they lie in the same plane 
and do not meet.  Without the coplanar part of the definition of ``parallel'',  Playfair's
axiom would fail in 3-space, while Euclid 5 would hold.   The proof at hand all points lie 
in the same plane;  that enters through the application of  Lemma~\ref{lemma:uniqueperpendicular}
(uniqueness of the erected perpendicular), which is essentially an upper dimension axiom.

\subsection{The strong parallel postulate}

We make two changes in Euclid 5  to get what we call the strong parallel postulate:
\smallskip

\begin{itemize}
\item  We change the hypothesis $\B(q,a,r)$ in Euclid's axiom to 
$\neg on(a,K)$.   In other words,  we require that the two  interior angles on the same side
of the transversal do 
not make exactly two right angles,  instead of requiring that they make less than two right angles.

\item  We change the conclusion to state only that $M$ meets $L$,  without specifying on which 
side of the transversal $pq$ the intersection lies.
 
\end{itemize}
\vskip -0.5cm

\begin{figure}[h]
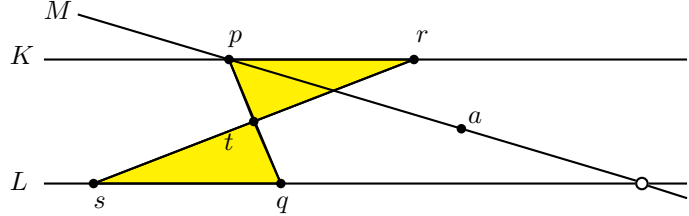
   
\center{\StrongParallelRawFigure}
\caption{Strong Parallel Postulate:  $M$ and $L$ must meet (somewhere) provided $a$ is not on $K$ and $pt=qt$ and $rt=st$.}
\label{figure:StrongParallelRawFigure}
\end{figure}

\FloatBarrier

Here is the strong parallel postulate. 
\smallskip 

\axioms
$ \B(p,t,q) \land \B(s,t,r)     $&  (Strong Parallel Postulate, SPP)\\
\endaxioms
\vskip-0.1cm
\axioms
$ \land\, \neg on(p,L)\land\, \neg\, \on(a,K) \land pt=qt \land rt=st \land p \neq r  $& \\
$\land\, L = \Line(s,q) \land\, M = \Line(p,a) \land K = \Line(p,r) $ & \\
$\qquad \implies on(\IntersectLines(L,M),L) \land on(\IntersectLines(L,M),M)   $ & 
\endaxioms
\smallskip

The strong parallel postulate differs from Euclid's 
version in that we are not required to know in what {\em direction} $M$ passes through $P$; but also the conclusion is weaker, 
in that it does not specify {\em where} $M$ must meet $L$.   In other words, the betweenness hypothesis of Euclid 5 is removed,
and so is the betweenness conclusion.  Since both the hypothesis and conclusion have been 
changed, it is not immediate whether this new postulate is stronger the Euclid 5, 
or equivalent, or possibly even weaker, but it turns out to be stronger--hence the name.

The strong parallel postulate fails in 3-space, because $a$
might not even lie in the same plane with $K$ and $L$.  Hence, the axiom should be used
(in this form) only 
in the presence of other axioms guaranteeing that space is two-dimensional.  We could, of 
course, also require that $a$ lie on some line passing through both $K$ and $L$, but 
since we are only concerned with plane geometry in this paper, we have avoided this complication.

The following lemma separates the strong parallel postulate into two parts:  Playfair, 
which is a negative statement about parallelism,  together with an assertion that two lines
that can't fail to meet, must in fact meet:

\begin{Lemma} \label{lemma:linestability}
[Stability of line intersection]
The strong parallel postulate is equivalent to Playfair's axiom together with 
``if $L$ and $M$ are non-coincident, non-parallel lines, then $M$ meets $L$.''
\end{Lemma}

\noindent{\em Remark}.  The name of the lemma arises from the fact that, expressed 
in symbols, the lemma says $\neg \neg t \defined \implies t \defined$, where 
$t = \IntersectLines(L,M)$.
\medskip

\noindent{\em Proof}.  
Assume the strong parallel postulate, and assume that $L$ and $M$ are non-coincident, 
non-parallel lines.  Our first goal is to prove the existence of a point $p$ on $M$ 
but not on $L$.   Since each line is given by two points, we can find two points $s$ 
and $t$ on $M$.  By the uniform perpendicular construction, we can construct line $J$ 
  through $s$  perpendicular to $L$.  Let $u$  
be the point of intersection of $J$  with $L$.   Let $ab$ be a segment 
long than any of $su$  for example obtained by extending the non-null segment 
$st$ by $su$. Construct a point $q$ on $J$ such that 
$qs=ab$, and a line $L^\prime$ through $q$ parallel to $L$.   Then $s$ lies on $M$ but not 
on $L^\prime$,  and $M$ is not perpendicular to $J$.  Applying the strong parallel postulate,
we can find a point $p$ where $M$ meets $L^\prime$.   That point lies on $M$, but not on $L$,
since $L^\prime$ is parallel to $L$.

Let $J$ be the perpendicular from $p$ to $L$, and let $K$
be the perpendicular to $J$ at $p$.  Let $M$ be a line through 
$p$ that does not coincide with $K$.  Let $a$ be a point on $M$ different from $p$. 
Then $a$ does not lie on $K$, since if it did, $M$ would coincide with $K$ and hence be 
parallel to $L$.  By the strong parallel postulate, $M$ meets $L$ as claimed.
That proves both Playfair and the quoted assertion in the theorem.  

Conversely, assume Playfair and the quoted assertion, and assume the hypotheses of the 
strong parallel axiom.  Then by Playfair, $M$ is not parallel to $L$, and by the 
quoted assertion, $M$ meets $L$, which is the conclusion of the strong parallel axiom.
That completes the proof.

\subsection{The strong parallel postulate proves Euclid 5}
The strong parallel postulate has a weaker conclusion than Euclid  5,  because it does not specify on which side of $P$ the intersection point will lie.
On the other hand, it also has a weaker hypothesis than Euclid 5, so its exact relationship 
to Euclid 5 is not immediately clear.   One direction is settled by the following theorem:

\begin{Theorem} \label{theorem:Euclid5} The strong parallel postulate implies  Euclid's Postulate 5 in 
neutral constructive geometry.
\end{Theorem}
 
\noindent{\em Proof}.  
 Suppose the hypotheses of Euclid~5 hold, as shown in 
Fig.~\ref{figure:EuclidParallelFigure}.  Specifically, let
  $L$ be a line, $p$ a point not on $L$,  $K$ parallel to $L$ through $p$,
$M$ another line through $P$, 
$q$ a point on $L$, $a$  a point on $M$ not on $pq$,  and $r$ the intersection of $qa$ with $K$. Suppose
that the interior angles made by $L$, $M$, and $pq$ make less than two right angles, which formally means that 
$a$ is between $q$ and $r$. 
By the strong parallel postulate,  $M$ does meet $L$ at some point $e$ (indicated by 
the open circle in Fig.~\ref{figure:EuclidParallelFigure}).   
It remains to show that $a$ is between $p$ and $e$.
By Markov's principle for betweenness,
 it suffices to prove that $p$ is not between $e$ and $a$, and $e$ is not between $p$ and $a$.
Suppose first that $p$ is between $e$ and $a$.  Then $e$ is on the opposite
side of line $K$ from $a$.  Segment $eq$ lies on line $L$, and line $L$ does not meet $K$, 
so segment $eq$ does not meet $K$.  Hence by Lemma \ref{lemma:sameside}, 
 $e$ and $q$ are on the same side of line $K$.   Since $e$ is on the opposite side of $K$ from $a$, 
it follows that $a$ and $q$ are on opposite sides of line $K$.  Hence point $r$, the intersection of $aq$ with $K$, must be between 
$a$ and $q$.  But that contradicts the fact that $a$ is between $q$ and $r$.  Hence the assumption that $p$ is between $e$ and $a$
has led to a contradiction.   Now suppose instead that $e$ is between $p$ and $a$.  Then $a$ and $p$ are on opposite sides of $L$.

But $rp$ does not meet $L$, since $rp$ lies on $K$, and  $K$ does not meet $L$.  Hence by 
Lemma~\ref{lemma:sameside}, $r$ and $p$ are on the same side of line $L$.  Hence $a$ and $r$ are on opposite sides of $L$.
But then the intersection  point of $ar$ and $L$, which is $q$, lies between $a$ and $e$, contradicting the fact that $a$ lies 
between $q$ and $r$.   Hence the assumption that $e$ is between $p$ and $a$ has also led to a contradiction.    Hence, as noted 
already, by Markov's principle for betweenness,  $a$ is between $p$ and $e$.  That completes the proof of Euclid's Postulate 5 from the strong parallel postulate.

\subsection{Triangle circumscription}
The {\em triangle circumscription principle} says, classically,  
that any three non-collinear points lie on a circle.  
It is classically equivalent to the parallel postulate.  In this section we show that
it is constructively equivalent to the strong parallel postulate, so we could have taken
it as an axiom in place of the strong parallel postulate. 
The name is appropriate, since
one can prove (even constructively) that a line intersects a circle in at most two points, and from that it is not difficult 
to prove Euclid III.1, that a chord of a circle lies inside the circle.  Hence, once a circle passes through the vertices of a triangle, the whole triangle lies inside the circle.

The triangle circumscription principle 
 is a convenient form of the parallel postulate
for use in a theory with only point variables, since it can be naturally stated without mentioning circles or lines:
given three non-collinear points, there is a fourth point equidistant from all three. 
 According to \cite{tarski-givant}, p.~190, 
Szmielew took it for her parallel axiom in her influential manuscript that formed the basis for Part I of \cite{schwabhauser},
but Schwabh\"auser adopted another form of the parallel axiom for publication.

The triangle circumscription principle makes sense constructively, but 
constructively there are several versions of it to consider.  One of them 
is non-constructive, but the others turn out to be equivalent,
though that is not  {\em a priori} clear.
  These principles are not just a curiosity:  we must study them 
because of their intimate connection to the geometric definition of multiplication.

First, the non-constructive version is this:  Given three points $a$, $b$, and $c$,
there is a circle passing through all three, unless they are distinct and 
collinear.  The problem with this is that if $b$ and $c$ are close together,
and about unit distance from $a$,  the circle through the three points does
not depend continuously on $b$ and $c$ as $c$ approaches $b$.  Hence this 
version is not constructively valid.   

That problem can be fixed drastically by requiring $a$, $b$, and $c$ to be 
three distinct points. That we call the triangle circumscription principle;
but that version is insufficient to define multiplication.   To formulate
a stronger version, we drop the requirement that $a \neq b$, but we require 
instead that $a$ and $b$ lie on a fixed line $L$.

The {\em strong triangle circumscription principle}   
says that if $a$ and $b$ lie on $L$,
and $c$ does not lie on $L$, then there is a point $e$ equidistant from $a$, $b$, and $c$,
and if $a=b$ then no other point of $L$ is on the circle with center $e$ through $a$, i.e., that
circle is tangent to $L$ at $a$.  To prove this principle 
constructively, the point $e$ should be constructed uniformly, that is, without a case distinction
as to whether $a=b$ or $a \neq b$.

The  {\em one-sided triangle circumscription principle} is similar to the strong triangle circumscription
principle, but with an additional hypothesis about the positions of the three points:
$a$ and $b$ do not lie on opposite sides of the perpendicular $W$ from 
$c$ to $L$.  (This corresponds to multiplication
of nonnegative numbers.)

\begin{figure}[ht]
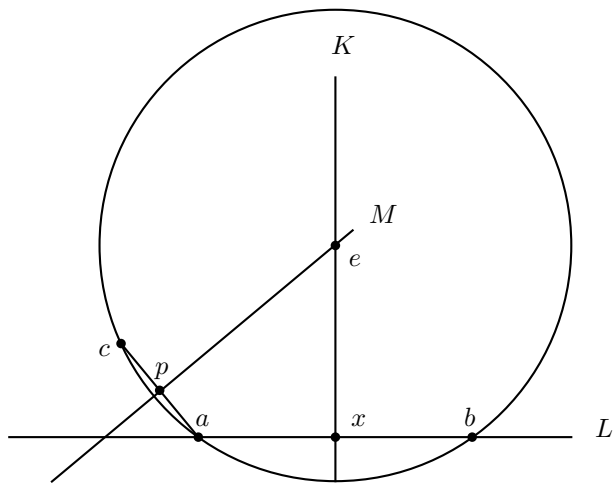
   
\center{\TriangleCircumscriptionFigure}
\caption{The strong parallel postulate implies strong triangle circumscription.
The center $e$ is where non-parallel lines $M$ and $K$ meet.}
\label{figure:TriangleCircumscriptionFigure} 
\end{figure}

\begin{Lemma} \label{lemma:strongparallelimpliescircumscription}
 In neutral constructive geometry,
 the strong parallel postulate implies the strong triangle circumscription principle.
\end{Lemma}

\noindent{\em Proof}. Assume the strong parallel axiom;
we will derive the strong triangle circumscription principle.
Suppose $c$ is not on line $L$, but $a$ and $b$ are on line $L$,
as shown in Fig.~\ref{figure:TriangleCircumscriptionFigure}.
By Lemma~\ref{lemma:uniformmidpoint},
 we find (by a uniform construction) a perpendicular $K$ to $L$, that bisects $ab$ 
 when $a \neq b$, and passes through $a$ when $a=b$.  Let $M$ be the perpendicular
bisector of segment $ac$.  Then $M$ and $K$ are not coincident, since if they 
were, then $ac$ and $L$ would both be perpendicular to $K$ through $a$, so 
by Lemma~\ref{lemma:uniqueperpendicular}, $c$ would lie on $L$, contradiction.
Similarly, $M$ and $K$ are not parallel, since if they were parallel, then 
  $ac$ and $L$ would be perpendicular to parallel lines, and both contain $a$, so 
  by Lemma~\ref{lemma:two-perps}, lines $L$ and $ac$ would coincide, contradiction.
  (Lemma~\ref{lemma:two-perps}
  requires Playfair's axiom, but the strong parallel postulate implies Playfair, so 
  the lemma is applicable.)   Therefore
 $M$ and $K$ are not parallel.   By the strong parallel postulate and
  Lemma~\ref{lemma:linestability},
   lines $M$ and $K$ meet in a point $e$. 
That completes the proof of the lemma.  Note that the proof does not assume anything about the order
of $a$ and $b$ on $L$.

\begin{Lemma} \label{lemma:Euclid5impliesone-sided}
In neutral constructive geometry,  
Euclid 5 implies the one-sided triangle circumscription principle.
\end{Lemma}

As in Fig.~\ref{figure:TriangleCircumscriptionFigure},
we have a line $L$ containing points $a$ and $b$, and a third point $c$ that does
not lie on $L$. The hypothesis of one-sided triangle circumscription tells us that
the perpendicular from $c$ to $L$ does not meet the closed segment $ab$, 
and hence does not meet $ax$,  where as shown $x$ is the midpoint of $ab$,
and $K$ is perpendicular to $L$ at $x$.  $K$ and $x$ are constructed by 
Lemma~\ref{lemma:uniformmidpoint}, so they are defined even if $a=b$.
As in the previous proof, let $M$ be the perpendicular bisector of $ac$.
Now let $J$ be the perpendicular to $L$ passing through $p$.  Let $f$ be the 
foot of $J$ on $L$.  We claim $J$ does not meet the open segment $ax$
(even though we do not know anything about the order of $a$ and $x$ on $L$).
 Under the hypothesis of the one-sided triangle 
circumscription principle,  the relative positions of 
$K$ and $M$ are as shown in  Fig.~\ref{figure:TriangleCircumscriptionFigure}, 
even though we do not know whether $\B(f,a,b)$ or $\B(f,b,a)$ or $a=b$.

Intuitively, what is going on here is this:  the $x$-coordinate of $J$ is $a/2$, 
and the $x$-coordinate of $K$ is $(a+b)/2$, and $a/2 \le (a+b)/2$.  But we cannot 
argue that way, since we need this lemma in order to define multiplication, which 
comes before introducing coordinates.   We must argue geometrically.   

Let $t$ be the foot of the perpendicular from $c$ to $L$.  By Lemma~\ref{lemma:midpointprojection},
 $f$ is the midpoint
of $ta$. (That lemma applies since we are allowed here to use Euclid 5, which implies
Playfair.)   Hence $\T(t,f,a)$.  We have to prove $\neg \B(a,f,x)$.  Suppose $\B(a,f,x)$.
If $\T(t,a,b)$, i.e., $a$ is left of $b$ in the diagram, then by the transitivity 
of betweenness we have $\T(f,a,b)$, so $f$ does not lie between $a$ and $b$,
and hence not between $a$ and $x$.  Hence $\neg\, \T(t,a,b)$.  Hence $\T(t,b,a)$.
That is, $b$ is to the non-strict left of $a$ on $L$, rather than as shown in the diagram.
Now we can apply Lemma~\ref{lemma:midpoint-helper}:  if $x$ is the midpoint of $ab$, and $f$ is 
the midpoint of $at$, and $\T(t,b,a)$, then $\neg\, \B(a,f,x)$.  We conclude $\neg\, \B(a,f,x)$,
contradicting our assumption that $\B(a,f,x)$; that contradiction proves without
assumptions that $\neg\, \B(a,f,x)$. 
 
That verifies the hypothesis that $J$ does not meet $ax$.     
By Lemma~\ref{lemma:euclidexample}, Euclid 5 implies that $M$ and $K$ meet.
That completes the proof.

\begin{Lemma} \label{lemma:circumscription-Playfair}
 The triangle circumscription principle implies Playfair's axiom.
\end{Lemma}

\begin{figure}[ht]
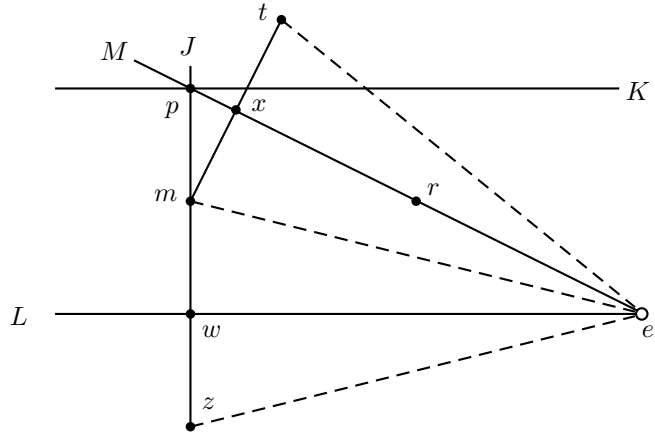
   
\center{\CenterImpliesPlayfairFigure}
\caption{Triangle circumscription classically implies Playfair. Given lines $L$
and $M$, to construct their intersection point as the center $e$ of an appropriate circle. }
\label{figure:CenterImpliesPlayfairFigure}  
\end{figure}

\noindent{\em Remark}. In the next lemma, we will prove that 
triangle circumscription proves the strong parallel postulate, which 
is stronger than this lemma; but we need this lemma to prove that one.
The figure and argument used for this proof will not work constructively
to prove the stronger theorem, because this proof requires a case 
distinction.  The construction (see Fig.~\ref{figure:CenterImpliesPlayfairFigure})
assumes $M$ is not perpendicular to $K$, so the construction 
given does not extend continuously to that case.
\medskip

\noindent{\em Proof}.  Since Playfair's axiom contains no existential 
quantifiers or disjunctions, we can argue classically.  
Assume the triangle circumscription principle; we will 
prove Playfair.  See Fig.~\ref{figure:CenterImpliesPlayfairFigure} for an illustration.
Suppose $L$ is a line and $p$ is a point not on $L$.
Drop a perpendicular $J$ from $p$ to $L$; let $w$ be the foot of this perpendicular,
so $pw \perp L$ at $w$.   
Let $K$ be perpendicular to $pw$ at $p$.  Then 
$K$ is parallel to $L$.   Let $M$ be a line 
through $p$ parallel to $L$ that does not coincide with $K$, as witnessed by 
a point $r$ on $M$ but not on $K$.  We intend to show that $M$ meets $L$.
(Playfair's axiom mentions two lines through $p$ parallel to $L$, and says they 
must coincide;  without loss of generality we can assume one of the lines is $K$.)
Our proof is not constructive, because it requires a case distinction 
as to whether $r$ is collinear with $pw$ or not.  If $r$ is collinear 
with $pw$ then $w$ lies on $M$ and hence $M$ meets $L$, so we are finished.
Therefore (classically) we can assume $r$ is not collinear with $pw$.  
Let $m$ be the midpoint of $pw$.  Since $p$ does not lie on $L$, but $w$ does 
lie on $L$, we have $p \neq w$ and hence $m \neq p$. 
Then $m$ does not lie on $M$, since if it did, $M$ would contain two distinct 
points $m$ and $p$ of $pw$.  Let $t$ be the reflection of $m$ in $M$; since 
$m$ does not lie on $M$, $t \neq m$.  Let $z$ be the reflection of $p$ in $L$;
since $pw \perp L$, $z$ lies on the the line $J$ containing $pw$.  Since $m \neq w$,
we have $m \neq z$. 

 We claim that $t$, $m$, and $z$ are not collinear.
Suppose that they are collinear.  Then $t$ lies on $J$. Since $t$ is 
the reflection of $m$ in $M$, $tm \perp M$.  Since $t$ lies on $J$, we have $J \perp M$.
Then $K$ and $M$ are both perpendicular to $J$ at point $p$.  Hence 
(by Lemma~\ref{lemma:uniqueperpendicular})
$K$ and $M$ coincide.  Then point $r$ (which lies on $M$) lies on $K$, contradicting
the hypothesis that $r$ is not on $K$.  This contradiction proves that $t$, $m$, and $z$
are not collinear.

Since $t$, $m$, and $z$
are not collinear, we can apply the triangle circumscription axiom to obtain 
a circle containing all three; let $e$ be the center of that circle. Technically
the axiom tells us $ey = ex = ez$.    Let $x$ 
be the midpoint of $tm$  (so $x$ lies on $M$).  Then by the definition of 
perpendicular, since $et = em$ we have $pe \perp tm$ at $x$.  But by construction 
of $t$, we have $tm \perp rp$ at $x$.  By the uniqueness of the perpendicular to $tm$ 
at $x$, the lines $pe$ and  $rp$ coincide.  But $\Line(r,p)$ is $M$.  Hence $e$ lies on $M$.
Similarly, 
$we \perp J$ at $w$ because $me=ze$, but also $L \perp J$ at $w$ by construction of $J$.
therefore $L$ coincides with $\Line(w,e)$, so $e$ lies on $L$.  Now we have proved
that $e$ lies on both $L$ and $M$, which was to be shown.  That completes the proof
of the lemma.

\begin{Lemma} \label{lemma:circumscriptionimpliesstrongparallel}
The triangle circumscription principle implies the strong parallel postulate.
\end{Lemma}

\begin{figure}[ht]
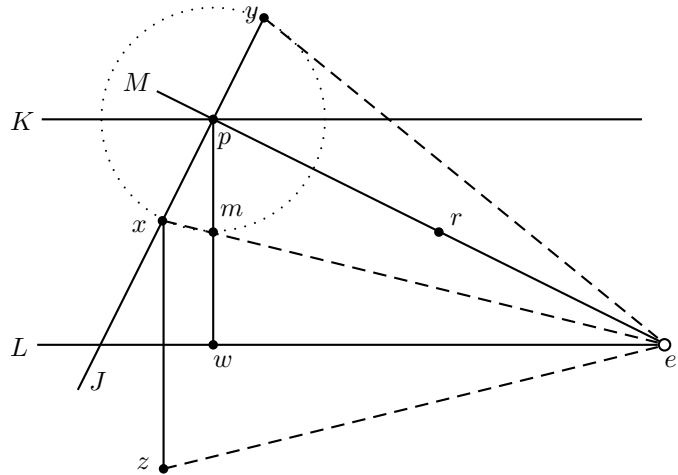
   
\center{\CenterImpliesStrongParallelFigure}
\caption{Triangle circumscription implies the  strong parallel postulate. Given lines $L$
and $M$, to construct their intersection point as the center $e$ of an appropriate circle.
$y$ and $z$ are reflections of $x$ in $M$ and $L$.
}
\label{figure:CenterImpliesStrongParallelFigure}  
\end{figure}

\noindent{\em Proof}.  Assume the triangle circumscription principle; we will 
prove the strong parallel postulate.  See Fig.~\ref{figure:CenterImpliesStrongParallelFigure}.
Suppose $L$ is a line and $p$ is a point not on $L$.  Let $pw$ be perpendicular
to $L$ at point $w$ on $L$ and let $K \perp pw$ at $p$.  Then 
$K$ is parallel to $L$.   Let $M$ be a line 
through $p$ that does not coincide with $K$.  We intend to show that $M$ meets $L$.

 Let $J$ be perpendicular to $M$ 
at $p$.  We need to construct a point $x$ on    $J$ but not on $M$ or $L$.
One way to do this is as follows: Let $m$ be the midpoint of $pw$, and 
choose $x$ on $J$ such that $px = pm$.  Since $p$ does not lie on $L$,
$m \neq p$.  Then $x \neq p$.  Since $x$ lies on $J$ and $J \perp M$, $x$
does not lie on $M$.  We claim $x$ does not lie on $L$ either:  if $x$ lies
on $L$, then $pwx$ is a right triangle,  whose hypotenuse $px$ is equal to half
the leg $pw$.  Now we have a right triangle with hypotenuse less than a leg,
  contradicting Lemma~\ref{lemma:legsmallerhypotenuse}.  That completes
  the proof that $x$ does not lie on $L$.

 Let $y$ be the reflection 
of $x$ in $M$.  Then $x \neq y$ since $x\neq p$.  
 The dotted circle in 
Fig.~\ref{figure:CenterImpliesStrongParallelFigure} illustrates
the fact that $px = pm = py$.  
Let $z$ be the reflection of $x$ in $L$;  since $x$ is 
not on $L$, $z \neq x$.  We claim that $x$, $y$, and $z$ are not collinear.
Once we prove that claim, we can finish the proof as follows:
By the triangle circumscription principle, there is a circle $C$ through $x$, $y$, and $z$.
 Let its center be $e$.  Then since $xy$ is a chord of $C$, its perpendicular bisector
$M$ passes through $e$.  Since $xz$ is a chord of $C$, its perpendicular bisector $L$ passes 
through $C$.  Hence $L$ and $M$ meet in point $e$.   This contradicts the assumption that $M$
is parallel to $L$.  Hence we (will have) have proved that every line $M$ parallel to $p$ through $L$ 
coincides with $K$.  In particular any two such lines coincide with each other.  

It remains to prove that $x$, $y$, and $z$ are not collinear.  Suppose that they are 
collinear.  Then $x$ and $z$ lie
on $J$, which is then perpendicular to both $L$ and $M$.  Since $J$ is perpendicular to $L$
and $K$ is parallel to $L$, $J$ is perpendicular to $K$, by
 Lemma~\ref{lemma:parallel-perp};  that lemma assumes Playfair,  but by Lemma~\ref{lemma:circumscription-Playfair}, 
 Playfair holds here.
  Hence $M$ and $K$ are both perpendicular to $J$ at $p$.  
 By the 
uniqueness of the perpendicular (Lemma~\ref{lemma:uniqueperpendicular}), $M$ and $K$ coincide, contradicting our assumption that $M$ 
and $K$ do not coincide.    That proves that $x$, $y$, and $z$ are not all collinear.
That completes the proof.

\begin{Theorem} \label{theorem:trianglecircumscription}
 The strong triangle circumscription principle, the triangle circumscription 
principle, and the strong parallel postulate are all equivalent in neutral constructive geometry.
\end{Theorem}

\noindent{\em Proof}.  The strong triangle circumscription principle implies the 
triangle circumscription principle, so the previous  lemmas establish a circle of 
implications between the three propositions.   Specifically, the triangle circumscription
principle implies the strong parallel postulate by 
Lemma~\ref{lemma:circumscriptionimpliesstrongparallel}, and the strong parallel 
postulate implies the strong triangle circumscription principle by 
Lemma~\ref{lemma:strongparallelimpliescircumscription}.
That completes the proof.

\section{Rotation and reflection}
In the introduction, we remarked on the meaning of a ``uniform reflection'' construction,
which must construct the reflection of a point $p$ in line $L$, without making a case distinction
according to whether $p$ is on $L$ or not.  In this section, we take up this construction
in detail.  We begin with the definition.

\begin{Definition} \label{definition:reflection}
Point $z$ is the reflection of $a$ in line $L$ provided (i) there is a line $K$
perpendicular to $L$ passing through both $z$ and $a$, and (ii) if $f$ is the 
intersection point of $K$ and $L$ then $af = zf$, and (iii) $\T(a,f,z)$.
\end{Definition}

Thus, if $a$ lies on $L$,  it is its own reflection in $L$,  and if it does not lie on $L$,
we can drop a perpendicular from $a$ to $L$, meeting $L$ at $f$, and extend the non-null 
segment $af$ by $af$ to construct $z$.   But to construct $z$ without making a case 
distinction is more difficult.    

Once we have a uniform perpendicular construction, the problem of reflection in a line
reduces to the problem of reflection in a point.  That is, to reflect $a$ in line $L$,
construct a perpendicular $K$ to $L$ passing through $a$;  let $f$ be the point where 
$K$ meets $L$, and reflect $a$ in the point $f$.  The result will be the reflection of 
$a$ in $L$.   

But reflection in a point, instead of in a line, is also problematic constructively,
even if the two points are known to lie on a fixed line.
Of course if $a \neq p$, we can reflect $a$ in $p$ just by extending $ap$ by $ap$.
And if $a=p$, then the reflection of $a$ in $p$ is just $a$.  But this construction
cannot be made uniform (that is, free of a case distinction), since only non-null
segments can be extended constructively.  Another failed attempt would be to 
find a point $b$ a fixed distance from $a$ with $\B(b,a,p)$ and then extend
the segment $bp$ by $ap$.  But, if we do not know on which side of $p$ lies, we
cannot construct such a point $b$, even if we assume both $a$ and $p$ lie on a known line $L$.

Since attempts to construct the uniform reflection directly appear blocked, we turn 
to another approach.
We use the fact that reflection of $p$ in point $a$ can be obtained by two successive 
rotations of $p$ about $a$, each rotation through a right angle.  It will therefore suffice
to define rotation of point $p$  through a given angle about point $a$ in such a way that no case distinction is 
required as to whether $p=a$ or not.  We call that ``uniform rotation.''  Of course, one way to define rotation is to use 
reflection; so if we had uniform reflection, we could get uniform rotation, and vice-versa.
We get off the ground by proving the existence of uniform rotation another way,
then deriving the existence of uniform reflection.

\subsection{Uniform rotation} \label{section:rotation}
In classical geometry, reflection and rotation can be defined without using the parallel
postulate.  However, we were unable to define uniform reflection or uniform rotation without
using Euclid 5.  In this section, we employ Euclid 5 freely.    

\begin{Lemma} \label{lemma:rotationhelper}
(Assuming Euclid 5.)
Let lines $L =\Line(\ell,b) $ and $J = \Line(d,b)$ meet at $b$, and
suppose that $L$ and $J$ do not coincide, and angle $\ell b d$ is less than a right angle.
    Suppose point $p$
lies on $ J$.  Then 
there exists a point $e$ on $L$ such that $ep$ lies on a line perpendicular to $J$.  
\end{Lemma}

\begin{figure}[ht]
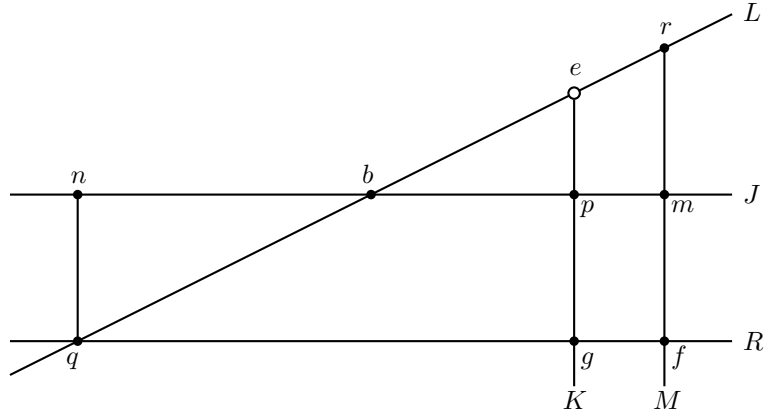

\center{\RotationHelperFigure}
\caption{The perpendicular to $J$ at $p$ meets $L$ in some point $e$.  Euclid 5
is needed to get $r$.}
\label{figure:rotationhelper}
\end{figure}

\noindent{\em Remark}. See Fig.~\ref{figure:rotationhelper}. 
$p$ could be on either side of $b$ or equal to $b$.  Imagine
$p$ moving back and forth along $J$,  passing through $b$.
  $e$ has to be constructed without a case 
distinction as to whether $p=b$ or not, or on which side of $b$ the point $p$ lies.
It would be technically incorrect to say $ep \perp J$, at least if $e=p$, 
since perpendicularity does not 
apply to null segments.
\medskip

\noindent{\em Proof}. Extend
segment $bd$ by $bp$ to point $m$.  
Let $n$ be the reflection of $m$ in $b$.  Then $\B(n,b,m)$ and hence $\B(n,p,m)$. 
Erect line $M$ perpendicular to $J$ at $m$.  Since angle $\ell bm$ is, by hypothesis,
less than a right angle, the interior angles made by traversal $J$ of $L$ and $M$
are together less than two right angles.  By Euclid 5, $M$ and $L$ intersect; 
call the intersection point $r$. 
  
Before proceeding, we discuss this application of Euclid 5.  We have presented it 
with Euclidean (im)precision, but technically we have to use our formulation of Euclid 5,
which requires an auxiliary point to witness that angle $\ell b m$ is less than a right 
angle.  The definition of ``less than'' for angles implies the existence of that point,
which would lie on $L$ somewhere to the right of $b$,
and witness that $\ell bm$ is less than a right angle by being between $m$ and a point
on the perpendicular to $J$ at $b$. 
For simplicity this point is not shown in Fig.~\ref{figure:rotationhelper}.

Let $q$ be the reflection of $r$ in $b$.   Then $q \neq b$ since $\ell \neq b$.
Construct $R$ parallel to $J$ through $q$.  Since we are assuming Euclid 5,
and $J$ is perpendicular to $M$,  also $M\perp R$.  Let 
$f$ be the point of intersection of $M$ and $R$.  There
 exists a line $K$ through $p$ perpendicular to $R$ (since $p$ is not on $R$, we 
 do not even need the uniform perpendicular construction).   Let
$g$ be the intersection point of $K$ and $R$. 

Now triangle $bnq$ is the reflection of triangle $bmr$ in point $b$. By
Lemma~\ref{lemma:reflectionpreservescongruence},
reflection preserves congruence.  Hence  triangle $bnq$ is congruent to triangle $bmr$.
 Therefore
$nq \perp J$.   Therefore each of the quadrilaterals $npgq$, $nmfq$, and $pmfg$ is a 
rectangle.   Therefore $np = qg$ and $nm = qf$ and $pm = gf$.  We already
proved $\B(n,p,m)$; therefore $\B(q,g,f)$.  That is, 
 $K$ meets side $qf$ of triangle $rqf$ at $g$.  $K$ does not meet $rf$ since
$K$ is parallel to $M$, which contains $rf$.  Then by Pasch's theorem (Theorem~\ref{theorem:pasch})
 $K$ meets 
$qr$;  the point of intersection is the desired point $e$, since $ep$ lies on $K$,
which is perpendicular to $J$.  Hence $ep \perp J$ as desired.  That completes
the proof of the lemma.

\begin{Theorem}[Uniform rotation]
(Assuming Euclid 5.) There is a construction $\Rotate$ that rotates
a point $p$ through a given angle $poq$, without a case distinction.  That is,
if $a$ lies on $\Line(o,p)$, then $z = \Rotate(p,o,q,a)$ is a point on $\Line(o,q)$
that is the reflection of $a$ in the line that bisects angle $poq$.
\end{Theorem}

\begin{figure}[ht]
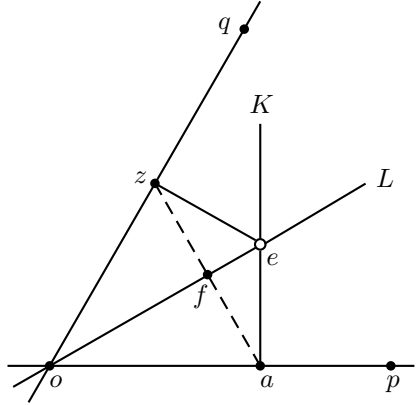
  
 \center{\RotationFigure}
\caption{$z = \Rotate(p,o,q,a)$.  $a$ is allowed to move along $\Line(o,p)$, even through $o$. }
\label{figure:RotationFigure}
\smallskip
\end{figure}

\noindent{\em Proof}.  Let $L$ be the angle bisector of angle $poq$.  (The method 
of Euclid I.9 for bisecting an angle is constructive as it stands; though Euclid 
paid insufficient attention to which of two possible triangles to select, that defect
is easily remedied.) 
 See Fig.~\ref{figure:RotationFigure}, but note that $a$ might also be on the other side of $o$,
or coincide with $o$.  Please imagine an animation in which $a$ moves back and forth on the 
horizontal line through $o$.
 
 The point is that $z$ should be defined even when $a=o$ (in which case it is just $o$, of course), and if $a$ moves along $\Line(o,p)$ through $o$, 
then $z$ moves along $\Line(o,q)$, passing through $a$ as $a$ passes through $o$.
     
The angle at $a$ is a right angle, and the angle at $o$ is half of angle $poq$.
Every angle is less than two right angles, so half of angle $poq$ is less than one right angle.
Let line $K$ be the perpendicular to $\Line(o,p)$ at $a$.  
By Lemma~\ref{lemma:rotationhelper}, $K$ meets $L$ at a point $e$.
Now, by the uniform perpendicular construction, 
there exists a point $z$ on  $\Line(o,q)$ such that $ez \perp \Line(o,q)$.
This point $z$ is the point we had to produce, the value of $\Rotate(p,o,q,a)$.
We claim $oz = oa$;  by the stability of equality, we can prove that by cases. If $a=o$
then $z=o$, so both $oz$ and $oa$ are null sequences, and hence congruent.  If $a \neq o$,
then triangles $zeo$ and $aeo$ are congruent (by SAS), so $oz=oa$.  That completes the proof 
that $oz=za$.  The perpendicular to $L$ through $z$ can be constructed without a 
case distinction by the uniform perpendicular construction.  Let $f$ be 
its intersection with $L$.  Similarly, let $g$ be the point of intersection 
with $L$ of the perpendicular to $L$ through $a$. We claim $f=g$ and $zf = af$.
By the stability of equality, we can prove that by cases.  Case 1, $a=o$.  Then 
$f=g=a=z$ and $af =gf$ because both are null segments.  Case 2, $a \neq o$.  
Then $aog$ and  $aof$ are triangles, and indeed they are right triangles
with equal hypotenuses and one pair of angles equal (the angles at $o$).   Hence
they are congruent triangles.  Hence 
$of=og$, so $f$ and $g$ coincide, and $zf = zf$.  That completes the proof 
by cases (and stability of equality) that $zf = af$ and $\T(z,f,a)$, and $a$ and $z$
both lie on a common perpendicular to $L$.   By Lemma~\ref{lemma:uniqueperpendicular}, 
we have $\B(z,f,a)$.  That is, the perpendicular from $z$ to $L$ also passes through $a$,
and since $zf=af$,  $z$ is the reflection of $a$ in $L$.  
By definition of reflection, that makes
$z$ the reflection of $a$ in $L$.  That completes the proof of the lemma.

\FloatBarrier
 
\subsection{Uniform reflection}
Now we will apply uniform rotation to show how to construct the reflection of a point
in a line uniformly.

\begin{Theorem} [Uniform reflection]\label{theorem:uniformreflection} 
(Assuming Euclid 5.) There is a term $\Reflect(x,L)$ that produces the reflection of $x$ in 
line $L$, without a case distinction whether $x$ is on $L$ or not.  That is, if $z = \Reflect(x,L)$
and $K$ perpendicular to $L$ contains $x$ and meets $L$ at $e$, then $z$ lies on $K$ with $ez=ex$,
and unless $x$ lies on $L$,  $\B(x,e,z)$.
\end{Theorem}

\noindent{\em Proof}.  The idea is that a reflection has the same effect on a single point $x$
as two ninety-degree rotations.  Using the uniform perpendicular construction, let $K$ be perpendicular to $L$,
and let $e$ be the intersection point of $K$ and $L$; then rotate $x$ twice by ninety degrees about $e$
using uniform rotation, obtaining the answer $z$.
 Fig.~\ref{figure:UniformReflectionFigure}
illustrates the construction.  (The  other points in the picture are used in uniform rotation.)

 \begin{figure}[h]   
\center{\UniformReflectionFigure}
\caption{Uniform Reflection}
 \label{figure:UniformReflectionFigure}
\end{figure}

 \subsection{The other intersection point} 
\begin{figure}[h]
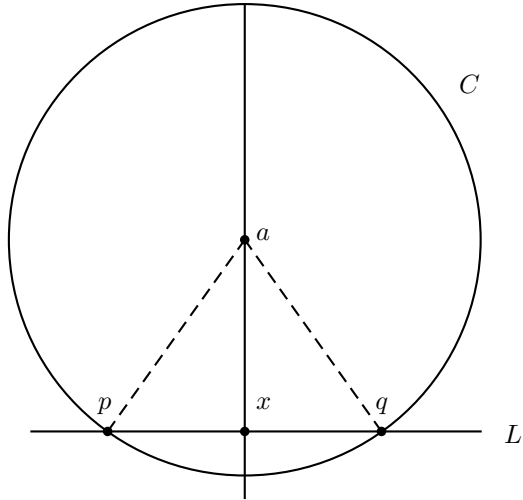
   
\center{\OtherIntersectionPointFigure}
\caption{Definability of ``the other intersection point''}
\label{figure:OtherIntersectionPointFigure}
\end{figure}

Many Euclidean constructions involve constructing one intersection point $p$ of a line $L = \Line(a,b)$ and a circle $C$, 
and then we say ``Let $q$ be the other intersection point of $L$ and $C$''.   The question is, whether in \ECG\ we can 
give a uniform method of constructing $Q$ from $p$, $L$, and $C$.  Yes, we can, as illustrated in 
Fig.~\ref{figure:OtherIntersectionPointFigure}.  The  proof is not as obvious 
as it may seem from the figure, since possibly $L$ could be a diameter of $C$ or tangent to $C$,
and we are not allowed to prove an existential theorem by cases.  Instead, we must give a single
construction that works for all cases.  We are able to do that using uniform perpendiculars
and uniform reflection.  Since uniform reflection requires Euclid 5, so does this construction.

\begin{Theorem} \label{theorem:otherintersectionpoint}
There is a term $Other(p,L,C)$ such that, if point $p$ lies on line $L$ and circle $C$, then
$Other(p,L,C)$ is defined and lies on both $L$ and $C$, and for all $z \neq p$, if $z$ lies on both $L$ and $C$ then $Other(p,L,C) = z$.
\end{Theorem}

\noindent{\em Proof}. Here is the construction script for $Other$:
\medskip

\begin{verbatim}
Other(Point P, Line L, Circle C)
{ a = center(C)
  K = Perp(a,L)
  q = Reflect(p,K);
  return q;
}
\end{verbatim}
\medskip
Since $\Perp$ is always defined, $K$ is defined, whether or not $L$ is a diameter of $C$,
and since $\Reflect$ is uniform reflection, $q$ is defined, whether or not $L$ is tangent to $C$.
But now that we have given a single construction that applies in all cases, we may argue 
by cases that $aq = ap$,  since by the stability of congruence, it suffices to prove
$\neg aq \neq ap$.  If $a$ does not lie on $L$ and $x \neq p$ then $axp$ is a triangle,
and is congruent to $axq$ by SAS; hence in that case $ap=aq$.  If $x=p$ then $q=p$ so 
in that case $ap=aq$.  If $a$ lies on $L$ then $a=x$ so $ap=aq$ because $xp=xq$. 
That completes the proof.
\medskip

We next  show that similarly, one can construct the ``other intersection point'' of two circles.
%We need a preliminary result:
%
%\begin{Lemma}[Construction of the radical axis] \label{lemma:radicalaxis}
%There is a construction that produces, given two non-concentric circles $C$ and $K$,  a line $L = RadicalAxis(C,K)$ such that 
%(i) if $C$ and $K$ intersect in two points, then $L$ passes through both intersection points,
%and (ii) if $C$ and $K$ meet in only one point $x$ then $L$ is tangent to each circle at $x$,
%and (iii) $L$ is defined whether or not $C$ and $K$ intersect.
%\end{Lemma}
%
%\noindent{Remark}. The line we construct is called the {\em radical axis} of $C$ and $K$; it 
%has other interesting properties that are not stated in the lemma. 
%\medskip
%
%\noindent{\em Proof}.  Fix two  non-concentric circles $C$ and $K$ (without any supposition about whether they intersect).
%Let $L$ be the line joining the centers of $C$ and $K$. 
%Let $R$ be another circle that intersects each of $C$ and $K$ in two distinct points, and whose center does not lie on $L$.  Then 
%let $K$ be the line joining the two points of intersection of $R$ and $C$, and let $M$ be the line joining the 
%two points of intersection of $R$ and $K$.  Then

\begin{Theorem}[The other intersection point of two circles] 
\label{theorem:otherintersectioncircles}
There is a term  $Other2(p,C,K)$ such that, if point $p$ lies on circle $C$ and on circle $K$, and $C$ and $K$ 
are not coincident,  then 
$Other2(p,C,K)$ is defined and lies on both $C$ and $K$, and forall $z \neq p$,
if $z$ lies on both $C$ and $K$ then $Other2(p,C,K) = z$.
\end{Theorem}

\noindent{\em Proof}.   
Here is the construction script:
\medskip

\begin{verbatim}
Other2(Point p, Circle C, Circle K)
{ a = center(C);
  b = center(K);
  L = Line(a,b);
  q = Reflect(p,L);
  return q;
}
\end{verbatim}
\medskip

Here is the correctness proof of the script.  We claim $a \neq b$.  To prove that,
assume $a=b$.  
Then $ap=bp$, so for all $z$,  $On(z,C)$ if and only if $za=ap$ if and only if $ab =bp$ if and only if $On(z,K)$.  That is, $C$ and $K$ coincide, contradicting the hypothesis.
Hence $a \neq b$. 
Therefore  $L$ in line 3 of the script is defined.  Since {\it Reflect} is always defined, by Lemma~\ref{theorem:uniformreflection},
$q$ is defined,  and by Lemma~\ref{lemma:reflectionindiameter}, $q$ lies on circles $C$ and $K$.
That completes the proof.

 \section{Geometrization of arithmetic without case distinctions}\label{section:geometrization}

Today we usually think of analytic geometry as coordinatizing a plane and translating 
geometrical relations between points and lines into algebraic equations and inequalities.
But the converse is also possible:  translating algebra into geometry,  and this is 
important for showing 
that the models of the geometry of constructions are planes over Euclidean fields.

In modern geometry books (such as \cite{borsuk-szmielew} or \cite{hartshorne}),
 arithmetic is geometrized as operations on congruence 
classes of segments.  To give a formal definition within geometry we must 
avoid sets.   We operate instead on points on
some fixed line $X = \Line(0,1)$, where $0$ and $1$ are two arbitrarily fixed points.  We 
refer to $X$ as the $x$-axis.  We erect a line $Y$ perpendicular to $X$ at 0, and call it the 
$y$-axis.  We mark off a point $I$ on $Y$ such that $0I = 01$. Then the
following are the key steps to be carried out:
\smallskip

\begin{itemize}
\item Coordinatization.   We need to assign coordinates $X(p)$ and $Y(p)$ (both on the $x$-axis)
to each point $p$.  We need to show that every pair of points on the $x$-axis arises as 
$(X(p),Y(p))$ for some point $p$.
\item Definition of addition.  Points on the $x$-axis can be added, without a case distinction as to sign.
\item Definition of multiplication.  Points on the $x$-axis can be multiplied, without a case distinction as to sign.
\item Definition of square root.  The geometrical definition of the square root of $x$ given by Descartes can be made to work without a case distinction whether $x$ is zero or positive.
\item The laws of addition, multiplication, and square root are satisfied, except for multiplicative inverse.
\item Depending on which parallel axiom is assumed, some form of the axiom of multiplicative inverse holds.
\end{itemize}
 
As far as I can tell, past work on these things has always assumed some discontinuous constructions, 
such as test-for-equality or at least apartness.  The closest thing to a constructive treatment is 
\cite{hartshorne}; the treatment there is discussed further in \S\ref{section:hartshornecoordinates}.
 Since coordinatization and arithmetic themselves are patently 
computable and continuous,  it is  unaesthetic to appeal to discontinuous and 
non-computable ``constructions'' to achieve coordinatization and arithmetization.  
Constructively, it is not only unaesthetic, it is incorrect. 

\subsection{Coordinatization}
With the $x$-axis and $y$-axis, and the points 0,1, and $I$ fixed,  we use the uniform 
perpendicular construction to project any point $p$ onto the $x$-axis, obtaining its
$x$-coordinate $X(p)$.  We project it also onto the $y$-axis, obtaining a point $y(p)$
on the $Y$-axis, and then we rotate that point $90^\circ$ clockwise,  using the uniform 
rotation construction defined above, to obtain $Y(p)$ on the $x$-axis.  That is the 
$y$-coordinate of $p$.  Note that $y(p)$, with lower-case $y$, is on the $y$-axis,
while $Y(p)$, with upper-case $Y$, is on the $x$-axis.  

The first half of coordinatization (the assignment of coordinates $(X(p), Y(p))$ to every point) is thus 
accomplished without any parallel postulate. But to show that every pair of coordinates 
corresponds to a point requires Euclid 5, not to get from $(x,0)$ and $(0,y)$ to $(x,y)$,
but to get from $(y,0)$ to $(0,y)$,  which requires rotation.

\begin{Lemma}[Coordinatization]  \label{lemma:coordinatization} 
Assuming Euclid 5, it is possible to use uniform rotation and 
perpendiculars to define a construction 
 $\MakePoint(x,y)$, which produces a point $p$ such that $x=X(p)$ and $y=Y(p)$
whenever $x$ and $y$ are on the $x$-axis. 
\end{Lemma}

\begin{figure}[ht]
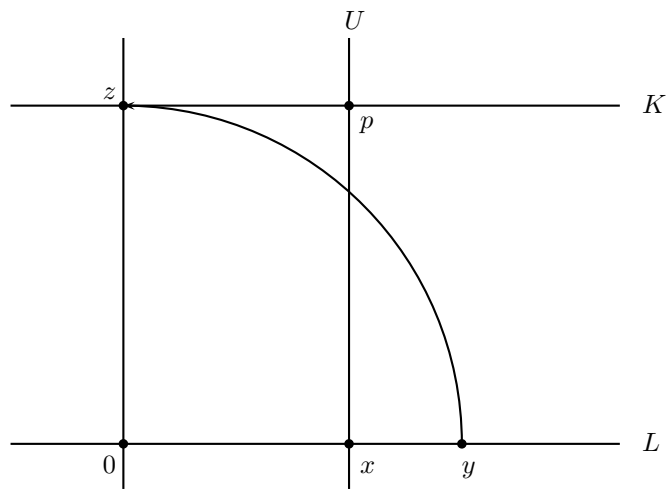
   
\center{\MakePointFigure}
\caption{$\MakePoint$ constructs a point $p$ with specified  coordinates $x$ and $y$. 
$K$ is constructed perpendicular to $U$ but by Playfair is also perpendicular to the $y$-axis.
}
\label{figure:MakePointFigure}  
\end{figure}

\noindent{\em Proof}.  See Fig.~\ref{figure:MakePointFigure}.
To construct the point $p$ from $x$ and $y$, 
first rotate $y$ by $90^\circ$ counter-clockwise,
obtaining a point $z$ on the $y$-axis.  (This is the only place we need Euclid 5.)
 Then
 erect a perpendicular $U$ 
to the $x$-axis at $x$.  Then let $K$ be perpendicular to $U$ through $z$, and let 
$p$ be the foot of that perpendicular on $U$.  Then $K$ and $L$ are both perpendicular 
to $U$.  We claim that $K$ is perpendicular to the $y$-axis.
By the stability of perpendicularity (Lemma~\ref{lemma:perpstable}), we can argue by cases to prove that.
If $U$ is distinct from the $y$-axis and $K$ is distinct from $L$, then $zpx0$ is a 
quadrilateral with 
three right angles, at 0, $x$, and $p$. 
Hence by Lemma~\ref{lemma:lambert}, it is a rectangle. (Here we need
Playfair's axiom, as Lemma~\ref{lemma:lambert} requires it, but we do not need Euclid 5.) Hence $K$ is 
perpendicular to the $y$-axis.   If $U$ coincides with the $y$-axis, then $K$ is 
perpendicular to the $y$-axis, since it is perpendicular to $U$.  If $K$ 
coincides with $L$, then $K$ is perpendicular to the $y$-axis, since $L$ is 
perpendicular to the $y$-axis. These cases are (classically) exhaustive, so by 
the stability of perpendicularity, we conclude that $K$ is perpendicular to the 
$y$-axis.  Then 
 $z$ is the projection 
of $p$ on the $y$-axis. 
See Fig.~\ref{figure:MakePointFigure}.
That completes the proof.   
\smallskip

\noindent{\em Remark}. If instead of constructing $K$ as a perpendicular to $U$ containing $z$,
we construct $K$ as perpendicular to the $y$-axis at $z$, then we have to prove $K$ meets $U$,
which can be done using the strong parallel postulate, but it is not clear how to do it using 
only Euclid 5.  The construction given in the proof avoids not only the strong parallel 
postulate but even Euclid 5, using only Playfair, although we still need Euclid 5 to make 
rotation work at the start of the proof.    It is crucial to avoid at least 
the strong parallel postulate, 
since we will eventually (in Corollary~\ref{lemma:Euclid5suffices}) use this lemma to help prove the strong parallel postulate from 
Euclid 5.

\subsection{Addition}
To perform addition geometrically we suppose given a line $L = \Line(R,S)$ and an ``origin'', 
a point $O$ on $L$  with $S$ between $R$ and $O$.   We need to define a construction $Add(A,B)$,
which also depends, of course, on $S$, $R$, and $O$, 
such that $Add(A,B)$ is a point $C$ on $L$ representing the (signed) sum of $A$ and $B$,
with $O$ considered as origin.       

\begin{Lemma} \label{lemma:informaladdition}
Given line L = Line(R,S), and a point $0$ on $L$ with $S$ between $R$ and $0$,  
we can construct a point $Add(A,B)$ on $L$ representing the signed sum of $A$ and $B$, with $O$ 
considered as origin,  using the elementary constructions and Circle3.
\end{Lemma}

\noindent{\em Remark}. In order to appreciate that this lemma is not trivial, consider the following 
obvious, but incorrect, attempted solution:
$$  Add(A,B):= \IntersectLineCircleTwo(\Line(0,B),Circle3(A,0,B)).$$
This works fine for $B \neq 0$, whatever the ``sign'' of $A$ and $B$, and it even works when $A=0$, but 
when $B=0$ it is undefined.
\medskip

\begin{figure}[h]
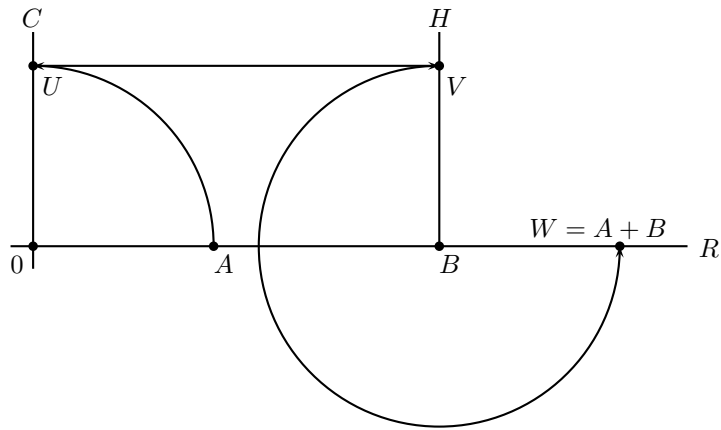
     
\center{\AdditionFigure}
\caption{Signed addition without test-for-equality. \goodbreak
$A$ is rotated to $U$, then projected to $V$, then rotated to $W$.}
\label{figure:AdditionFigure}
\end{figure}

\begin{figure}[h]
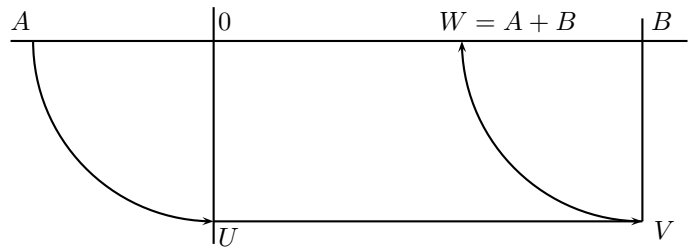
   
\center{\AdditionFigureTwo}
\caption{Signed addition when $A$ is negative}
 \label{figure:AdditionFigureTwo} 
\end{figure}

\noindent{\em Proof}.
With  $\Rotate$ in hand, we can give a construction for $Add(A,B)$ (depending also on $R$, $S$, and $O$).
(The construction is illustrated in Fig.~\ref{figure:AdditionFigure} and Fig.~\ref{figure:AdditionFigureTwo})
First, we replace $R$  with a new point on $L = \Line(R,S)$,  farther away from $O$, so that $O$, $A$, and $B$ are all on the 
same side of $R$, and the new $R$ and $S$ are in the same order on line $L$ as $R$ and $S$ were 
before.  Such points can be constructed using $\Extend$, but we omit the details.
Now erect the perpendicular $K$ to $L$ at $O$, and the perpendicular $H$ to $L$ at $B$.
In the process of erecting these perpendiculars, we will have constructed a points $C$ on $K$   
such that $ROC$ is a right turn.  Then let $D$ be the projection of $C$ on $H$ and let 
\begin{eqnarray*}
U &=& \Rotate(R,O,C,A) \\
V &=& \project(U,H) \\
W &=& \Rotate(D,B,R,V)
\end{eqnarray*}
We set $Add(A,B)= W$.
Then $Add(A,B)$ is defined for all $A,B$.  Suppose $A \neq O$.  Then $UV$ is perpendicular to both $K$ and $H$.  
Then $U$ and $V$ are on the same side of $L$, since if 
$UV$ meets $L$ at a point $X$, then $XU$ and $XO$ are both perpendicular to $K$, which implies $U=O$, 
which implies $A=O$, contradicting $A \neq O$.   It then follows from the property of $\Rotate$ that $B$ and $W$ occur on line $L$ in the same 
order that $O$ and $A$ occur.  Refer to Fig.~\ref{figure:AdditionFigureTwo} for an illustration of the case when $A$ is negative. 
This implies that $Add(A,B)$ represents the algebraic sum of $A$ and $B$, 
since in magnitude $BW = OA$.    

\begin{figure}[h]
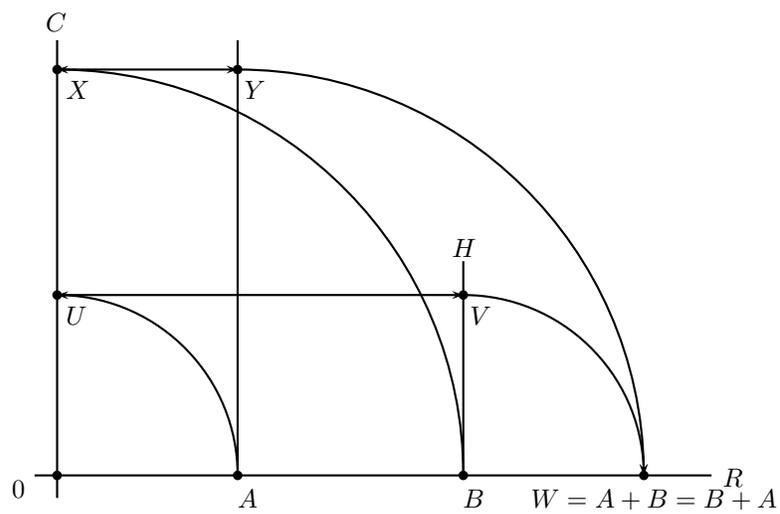
     
\center{\CommutativeAdditionFigure}
\caption{ Commutativity of addition with $A$ and $B$ positive.}
 \label{figure:CommutativeAdditionFigure}
\end{figure}

\begin{figure}[h]
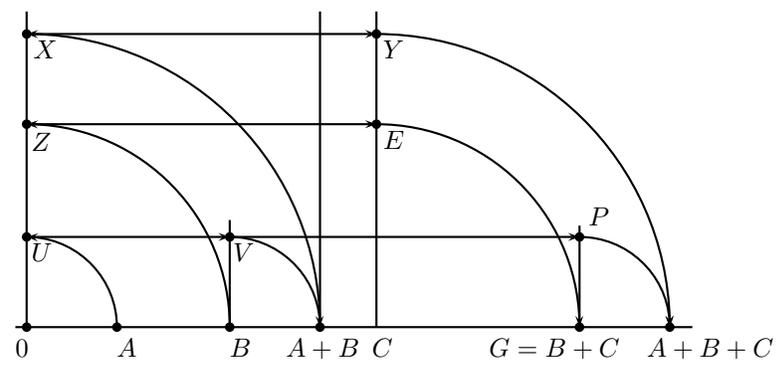
   
\center{\AssociativeAdditionFigure}
\caption{ Associativity of addition, positive arguments.}
\label{figure:AssociativeAdditionFigure}
\end{figure}

We indicate by pictures (Fig.~\ref{figure:CommutativeAdditionFigure} and Fig.~\ref{figure:AssociativeAdditionFigure}) how the commutativity and associativity of addition can 
be proved (without case distinctions).  We note that the definition of uniform 
rotation and the stability of betweenness allows us to argue using case distinctions, 
because double negations can be pushed through and eliminated by stability.

\FloatBarrier

\subsection{Unsigned multiplication}  \label{section:informalmultiplication}
\begin{figure}[ht]
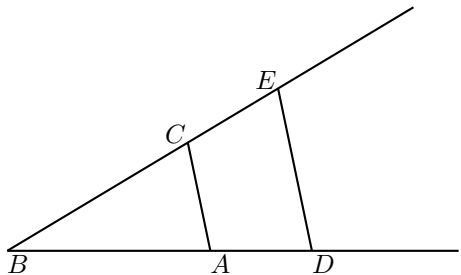
 
\center{\DescartesMultiplicationFigure}
\caption{{\em La Multiplication} according to Descartes}
\label{figure:DescartesMultiplicationFigure}
\end{figure}

The geometrical definitions of multiplication and square root go back (at least) to Descartes.%
\footnote{Although Descartes is usually credited with this, compare Euclid VI.12, which is 
very similar to Descartes's treatment of multiplication and division, and Euclid VI.13, which 
is very similar to Descartes's construction of square roots.}
On the second page of {\em La Geometrie} \cite{descartes}, he gives
constructions for multiplication and square roots.  We reproduce the drawings found 
on page 2 of his book \cite{descartes} in Figures \ref{figure:DescartesMultiplicationFigure} and 
\ref{figure:SquareRootFigure}.

\smallskip

Here is Descartes' explanation of this figure:
\begin{quote}
{\em 
\begin{enumerate}
\item{  Let $AB$ be taken as unity.}
\item{  Let it be required to multiply $BD$ by $BC$.  I have only to join the points $A$ and $C$,
and draw $DE$ parallel to $CA$; then $BE$ is the product of $BD$ and $BC$.}
\item{  If it be required to divide $BE$ by $BD$, I join $E$ and $D$, and draw $AC$ parallel to $DE$; then $BC$ is the result of the division.}

\end{enumerate}
}  % end em
\end{quote}
 
From the point of view of constructive geometry, there is a problem with the construction.  Namely, Descartes has
only told us how to multiply two segments with non-zero lengths, and at least one of whose lengths is not 1 (the length of unity--he needs this
when constructing $AC$ parallel to $DE$), while we want to be able to multiply in general,
without a test-for-equality construction.   

\begin{figure}[ht]
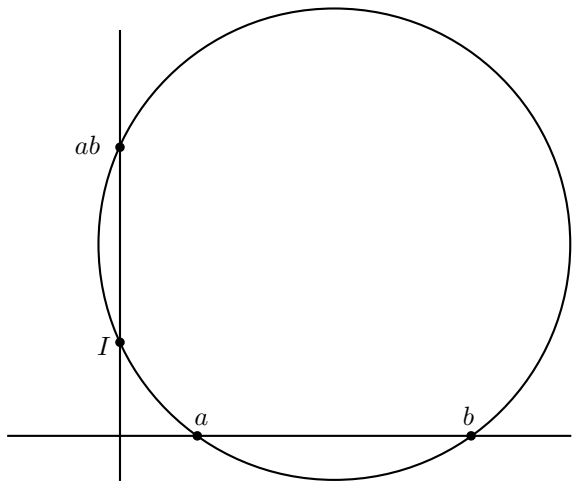
   
\center{\HilbertMultiplicationFigure}
\caption{Multiplication according to Hilbert}
 \label{figure:HilbertMultiplicationFigure}
\end{figure}

Using the $\para$ construction of Lemma~\ref{lemma:uniformparallelinformal} where Descartes calls for ``drawing $DE$ parallel to $CA$'',  we no longer have a problem 
multiplying numbers near 1 or 0.  In \cite{beeson-kobe}, we gave a construction that successfully generalizes
Descartes's multiplication method to signed arguments.  That method uses $\Rotate$.  
However, Hilbert (\cite{hilbert1899}, p. 54) gives another construction, whose result is equivalent to that of Descartes
for positive arguments,  but which directly works correctly for signed arguments as well.  It is illustrated in 
Fig.~\ref{figure:HilbertMultiplicationFigure}.

The construction is as follows:  Start with $\Line(0,1)$.  Erect a  perpendicular to  $\Line(0,1)$ at 0, and 
let $I$ be a point on it such that $0I=01$.  
 Given points $a$ and $b$ on $\Line(0,1)$,  construct a circle passing through $I$, $a$, and $b$.   This must 
 be done uniformly, so that when $a=b$, the circle passes through $I$ and is tangent to $\Line(0,1)$ at $a$.  Assume for the moment that such a circle exists. 
  Then the result of the construction is ``the other point of intersection''  of the circle and $\Line(0,I)$, except that this point 
lies not on $\Line(0,1)$ but on $\Line(0,I)$, so it needs to be rotated clockwise back to $\Line(0,1)$.
This ``other point of intersection'' can be defined by a term,  according to Theorem~\ref{theorem:otherintersectionpoint}. 

The point at issue is,  what version of the parallel postulate is required to prove the 
existence of the required circle?  We have already done the work:  If we restrict
$a$ and $b$ to be nonnegative, then the one-sided triangle circumscription principle 
is exactly what is needed, and Lemma~\ref{lemma:Euclid5impliesone-sided} shows that 
Euclid 5 suffices. 
 
\begin{Theorem} \label{theorem:unsignedmultiplication}
Assuming Euclid 5, 
 Hilbert multiplication as shown in Fig.~\ref{figure:HilbertMultiplicationFigure}
can be defined for all $a$ and $b$ on the non-negative $x$-axis. That is, 
it is defined 
by a term $\HilbertMultiply(a,b)$.
\end{Theorem}

\noindent{\em Proof}.  The proofs of Lemma~\ref{lemma:Euclid5impliesone-sided} and
Theorem~\ref{theorem:otherintersectionpoint} provide the required term.
\medskip

\begin{Lemma} \label{lemma:ringlaws}  Playfair's axiom implies  the commutative ring
laws, with $\HilbertMultiply$ for multiplication and $Add$ for addition, whenever
the terms involved are defined; that is, assuming the required circles and points of 
intersection exist.
\end{Lemma}

\noindent{\em Proof}.  The theorems in question have the form $t=s$, where $t$ and $s$ are
two terms for geometrical constructions.  The formulas $t\defined$ and $s\defined$ are 
existential formulas, specifying that certain ``witnessing points'' exist satisfying certain 
quantifier-free and disjunction-free relations.  
Hence $t\defined \land s\defined \implies t = s$  is equivalent
to a universally quantified negative formula.  Hence, by the double-negation interpretation 
from \cite{beeson2015b}, it suffices to prove this classically.  But that has been 
done many times, for example in \cite{schwabhauser}. Alternately one can check the 
proofs in \cite{schwabhauser} directly to see that they are constructive.
 That completes the proof.

\subsection{An alternative approach to constructive addition} \label{section:hartshornecoordinates}
The geometrization of arithmetic in \cite{hartshorne} can be made constructive,
providing an alternative to the definition given above.  We show how this is done.
Again we fix a line $L$, which we also call the $x$-axis, and points $0$ and $1$ on $L$,
and we construct the $y$-axis perpendicular to $L$ at 0. 
Then $x+1 := ext(0,1,0,x)$ works as the notation indicates for all $x$ on $L$,
 and $x-y := ext(x+1,x,0,y)$ works for $y \ge 0$.  Unary minus, $-x$, is defined
 by uniform reflection in the $y$-axis.   
 We define $\vert x \vert = ext(-1,0,0,x)$.  This is better than 
$\vert x \vert = \sqrt{\vert x^2 \vert}$, since it does not require
  square roots
or multiplication (hence requires no parallel postulate.)

Then, using 
the binary subtraction just defined, we define the negative part of $a$ by 
$$ a^- := \frac 1 2 \big(a - \vert a \vert\big).$$
The factor $\frac 1 2$ represents a midpoint construction, not a multiplication.
The midpoint construction in question is the uniform midpoint of segment $0x$
assuming that $x$ lies on line $L$, as defined in Lemma~\ref{lemma:uniformmidpoint}.  
Thus there is a term with free variable $a$ defining $a^-$.
Then we define $a^+ := -((-a)^-)$, using the unary minus defined above.
Finally we define 
\begin{equation}
 a + b := (a^+ + b^+) - (a^- + b^-). \label{eq:posnegparts}
 \end{equation}
Then we still have to prove the laws of addition and the distributive law.
For that, we can appeal to the double-negation interpretation in \cite{beeson2015b}.%
\footnote{
Note that both 
definitions of addition require the uniform reflection construction, which in turn 
requires a constructive definition of rotation.  Therefore the extra work involved 
in making the definition constructive is essentially the same for both approaches; 
it is just a matter of taste whether one likes 
an approach that is explained with pictures or with equations.  
}

\subsection{From unsigned to signed multiplication}
Hilbert's definition of multiplication, works without case distinctions about the 
signs of $a$ and $b$, but it requires the strong triangle circumscription principle,
while we wish to use only Euclid 5.  We were able to define Hilbert multiplication
for nonnegative arguments using Euclid 5.  Now we define
 a new multiplication $a\cdot b$, which agrees with $\HilbertMultiply(a,b)$
for nonnegative $a$ and $b$ by definition.  Later we will show that 
it agrees for all $a$ and $b$.

 We make use of the positive and negative parts $a^+$ and $a^-$
defined in (\ref{eq:posnegparts}).  
Motivated by the idea that, since $x=x^+-x^-$, we should have 
$$ a \cdot b =  (a^+ - a^-)\cdot(b^+ - b^-),$$
we define
$$ a\cdot b := (a^+b^+ + a^-b^-) - (a^-b^+ + a^+b^-)$$
On the right, multiplication means $\HilbertMultiply$, and is only applied to 
nonnegative elements.  On the right, addition means the geometric $\Add$, and 
subtraction means the geometrically defined binary subtraction.  Hence the right side
abbreviates a term representing a geometrical construction.  That can be taken as 
the definition of multiplication for arbitrary $a$ and $b$ on the $x$-axis.

Multiplication defined in this way satisfies the commutative law:
\begin{eqnarray*}
  a \cdot b  &=& (a^+b^+ + a^-b^-) - (a^-b^+ + a^+b^-) \\
  &=& b \cdot a
\end{eqnarray*}
using the commutativity of $\HilbertMultiply$ and the
commutativity of $Add$ on the right.
Multiplication similarly can be shown to satisfy the associative and distributive laws.
Therefore we may make algebraic calculations using $a \cdot b$ in the usual way.

We need to develop a bit of analytic geometry using this multiplication.
We define the ``length of segment $pq$'', denoted by $\vert pq \vert$,  to be the point 
$r$ on the positive 
$x$-axis such that $0r = pq$.  Officially, $\vert pq \vert $ is constructed
as follows: Let $u = -1$ be the 
reflection of $-1$ in the $y$-axis.  Then $\vert pq \vert = ext(u,0,p,q)$.
Then two segments are congruent if and only if they have the same length;
that follows from the stability of equality and the fact that a segment cannot
be congruent to a proper subsegment.

We use the usual notation for coordinates: instead of   $p = \MakePoint(p_1,p_2)$,
we just write $p = (p_1,p_2)$, and we use subscripts to indicate the coordinates of a point.
In this setting we can formulate a version of
 the Pythagorean theorem:  if $pqr$ is a 
right triangle, with the right angle at $q$, then 
$$\vert pq \vert^2 + \vert qr \vert^2 = \vert pr \vert^2.$$
where $x^2$ means $x \cdot x$.  Note that here only points on the nonnegative $x$-axis
are multiplied.   One can prove this version of the Pythagorean 
theorem, but it is not trivial.  It is Proposition 20.6, p.~178 of \cite{hartshorne}.
(Hartshorne defines multiplication on equivalence classes of segments; the segments
of the form $0p$ with $p$ on the non-negative $x$-axis can be taken as coset representatives,
so Prop. 20.6 does represent the theorem we need.)  The proof is perfectly constructive,
so there is no need to repeat it here.  

We do need a slight improvement:  let $L$ and $K$
be two perpendicular lines meeting at $q$, and let $p$ and $r$ lie on $L$ and $K$ 
respectively (but we do not assume they are different from $q$).  Then the equation of 
the Pythagorean theorem still holds.   Because of the stability of equality, we 
can prove this by cases.  (See \S \ref{section:casesplitsinEuclid} if needed,
to review why stability of equality permits us to argue by cases here.)

 If $p$ and $r$ are different from $q$, it is the Pythagorean
theorem from Prop.~20.6.   If either $p=q$ or $r=q$, then one of the sides has length 0
and the other two sides coincide, so the equation is trivially true.  Hence the 
Pythagorean theorem is constructively true, even for possibly degenerate right triangles.

\begin{Theorem} \label{theorem:signedmultiplication}
[Euclid 5 implies multiplication]
Assume Euclid 5 and line-circle continuity.  Fix a line $L$ for the $x$-axis, and points 0,
1, and $I$ for use in coordinatization.  Then $\HilbertMultiply(a,b)$ is defined
for all $a,b$ on $L$.
\end{Theorem}

\begin{figure}[ht]
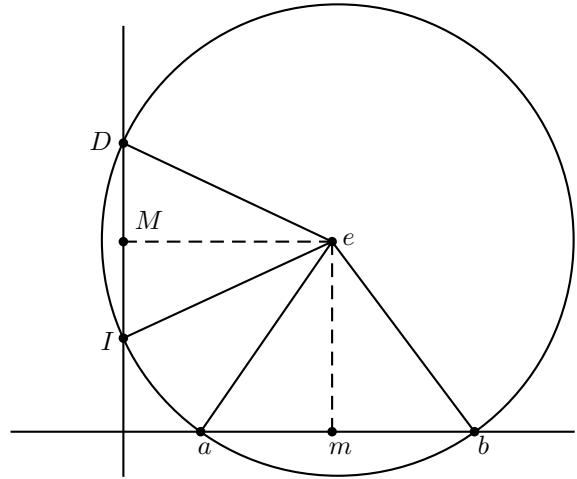

\center{\JustifyMultiplicationFigure}
\caption{$d = a \cdot b$ makes $eI=ea$, so $d=\HilbertMultiply(a,b).$}
\label{figure:JustifyMultiplicationFigure}
\end{figure}

\noindent{\em Proof}. 
By Theorem~\ref{theorem:unsignedmultiplication}, $\HilbertMultiply$ is defined
for nonnegative $a$ and $b$.  It follows that $a \cdot b$ is defined for all $a$ and
$b$. 
We will show that, for all $a$ and $b$, 
$\HilbertMultiply(a,b)$ is equal to $a\cdot b$. 
 Let $d = a \cdot b$.  Then $d$ is a point on the $x$-axis;  let $D$ be the 
rotation of $d$ by ninety degrees counterclockwise, so $D$ lies on the $y$-axis.
Using uniform rotation, we do not need a case distinction on the sign of $d$. 
See Fig.~\ref{figure:JustifyMultiplicationFigure}.  Although the picture shows
$a \le b$, the proof does not assume that.

We must show (i) there is a circle through $a$, $b$, and $I$, and (ii) that circle
passes through $D$.  
While in general we cannot construct a uniform midpoint of an arbitrary (possibly null) 
segment 
$ab$, Lemma~\ref{lemma:uniformmidpoint} says that we can construct a uniform 
midpoint of a (possibly null) segment $ab$ if we know that $a$ and $b$ are restricted to 
lie on a fixed line.  Let $L$ be the fixed line serving as the $x$-axis.
Let $m$ be the uniform midpoint 
of segment $ab$ (constructed by Lemma~\ref{lemma:uniformmidpoint} using $L$), 
and let $M$ be the uniform midpoint of segment $ID$ (constructed
 by Lemma~\ref{lemma:uniformmidpoint} using the $y$-axis for the line).
  Let $M^\prime$ be the 
rotation of $M$ by ninety degrees clockwise; using uniform rotation this can be
constructed without a case distinction about the position of $M$ on the $y$-axis.
 Let  $e = MakePoint(m,M^\prime)$, 
which is defined, assuming Euclid 5, by Lemma~\ref{lemma:coordinatization}.  Then 
$eM0m$ is a rectangle (unless some of its sides coincide, but by stability of equality,
its opposite sides are equal, whether or not they coincide).  Then $eD=eI$ since $e$
lies on the perpendicular bisector of $DI$ (technically, by the stability of equality
we can argue by cases, whether $D=I$ or not).  Similarly $ea=eb$.  To complete the 
proof, it suffices to show that $ea=eI$.  The point $e$ has been constructed uniformly,
with no case distinction about the signs or relative positions of $a$ and $b$;
but now, by the stability of equality, we are allowed to argue by cases for $ea=eI$.
That is, however, not necessary.  

I do not know a direct geometrical proof that $ea=eI$.  We  use analytic 
geometry.  The following algebraic calculations, however, are officially abbreviations 
for geometrical constructions, i.e., the multiplications refer to the multiplication 
$x \cdot y$  defined above, 
and the additions to $\Add$.  Divisions by 2 refer to a uniform midpoint construction.

 Let $R = \vert eI \vert$ and $r = \vert ae \vert$. 
  We make use of the Pythagorean theorem as described above.
We also make use of the distributive law, and of the properties of addition, whose proofs
have been sketched above. 
\begin{eqnarray*}
e &=& \bigg(\frac {a+b} 2, \frac{d+1} 2\bigg) \\
r^2 &=& \bigg(\frac{d+1} 2\bigg)^2 + \bigg(\frac{b-a}2\bigg)^2 \\
&=& \frac{d^2} 4 + \frac 1 4  + \frac{b^2} 4 + \frac {a^2 } 4 
\end{eqnarray*}
\begin{eqnarray*}
R^2 &=& \bigg(\frac {a+b} 2\bigg)^2 + \bigg(\frac {d-1} 2\bigg)^2  \\
&=& \frac {a^2} 4 + \frac {b^2} 4 + \frac {d^2 } 4 + \frac 1 4 \\
&=& r^2
\end{eqnarray*}
Since $R^2 = r^2$ and both $r$ and $R$ are nonnegative, we have $r=R$.
That completes the proof.

\subsection{Reciprocals}

\begin{figure}[h]
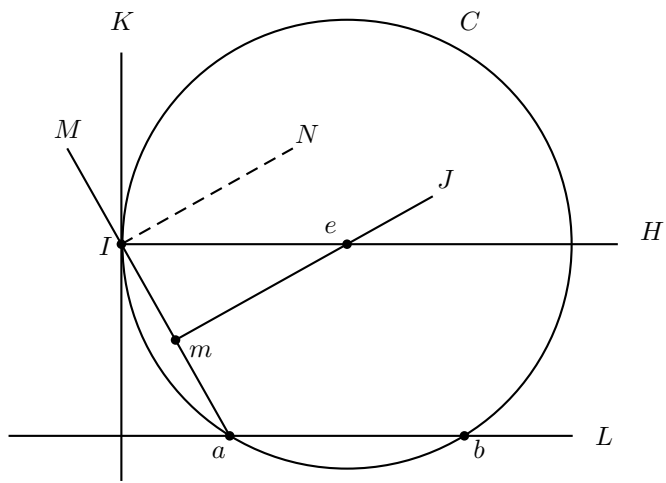
   
\center{\ReciprocalScriptFigure}
\caption{Construction of the reciprocal of $a$}
 \label{figure:ReciprocalScriptFigure}
\end{figure}

 According to the definition of $\HilbertMultiply$, 
a reciprocal of a nonzero point $a$ on $L = \Line(0,1)$ is a point $b$ on $L$ such that there is a  circle $C$ 
through $a$, $b$, and $I$, and $C$ is tangent to $K = \Line(0,I)$ at $I$.   The reciprocal of $x$ can be constructed
by the following script.  It uses the script $Other$ for the other intersection point of a line and circle, 
which was defined and proved correct in Theorem \ref{theorem:otherintersectionpoint}.  The script is illustrated 
in Fig.~\ref{figure:ReciprocalScriptFigure}.  The dashed line $N$ is used in the proof, but not in the construction.

\medskip

\begin{verbatim}
Point Reciprocal(Point a)
{ K = Perp(L,0)
  H = Perp(I,K)
  M = Line(a,I)
  m = Midpoint(a,I)
  J = Perp(m,M)
  e = IntersectLines(H,J)
  C = Circle(e,a)
  b = Other(a,L,C)
  return b
}
\end{verbatim}

\medskip

The following lemma proves the correctness of this script.  The hard part is to show that $\Reciprocal(x)$
is defined for $x \neq 0$.

\begin{Lemma} [Reciprocals] \label{lemma:reciprocal}
Let $\Reciprocal$ be defined as in the construction script above.   The strong parallel postulate 
implies that  for $a \neq 0$, 
$\Reciprocal(a)$ is defined, and if $b = \Reciprocal(a)$, then $\HilbertMultiply(a,b) = 1$.
Euclid 5 implies the same conclusion under the hypothesis $a > 0$ instead of $a \neq 0$.
\end{Lemma}

\noindent{\em Proof}. 
Assume $a \neq 0$.  We will show $\Reciprocal(x)$ is defined, by going through the script line by line. 
Since $\Perp$ is everywhere defined, $K$ and $H$ are defined;  since $I$ does not lie on $L$, $a \neq I$, 
so $M$ is defined.  Since $a \neq I$, $m$ is defined; since $\Perp$ is everywhere defined, $J$ is defined.
Now we have to prove that lines $H$ and $J$ meet in some point $e$.  That is where we use the 
strong parallel postulate.  Let $N = \Perp(I,M)$.  Then by Lemma~\ref{lemma:oneparallel}, $N$ and $J$ 
are parallel, since both are perpendicular to $M$.  Both $N$ and $J$ pass through $I$.
$M$ is not perpendicular to $H$ at $I$, since if it were, then $Ia$ and $I0$ would be two 
perpendiculars to $H$ at $I$, so they would coincide, and it follows that $a=0$, contrary to 
hypothesis.  Then let $x$ be any point on $H$ other than $I$, for example, $x=\IntersectLineCircleOne(H, \Circle(I,0,1)$.
  Then $x$ is not on  $N$, since if it were then $H$ and $N$ would coincide.
Hence we can apply the strong parallel postulate to conclude that $H$ meets $J$.  Then $e$ in the script 
is defined.  Since $\Circle$ is always defined, $C$ is defined.  Since $a$ lies on both $L$ and $C$, 
$b$ in the penultimate line is defined.  Hence $\Reciprocal(a)$ is defined.

By Lemma~\ref{lemma:tangentcircle}, $C$ is tangent to $K$ at $I$.
Then by definition of $\HilbertMultiply$, we have $\HilbertMultiply(a,b) = 1$, since circle $C$
passes through $a$ and $b$ and meets $K$ in $I$, and the ``fourth point of intersection'' is 
the point $I$,  since $C$ is tangent to $K$ at $I$.  That completes the proof.

I do not know of a direct verification of the distributive law for Hilbert multiplication.  Hilbert
showed that his definition is equivalent to that of Descartes and verified the distributive law for that definition.   Once the terms for $\Add(a,b)$ and $\HilbertMultiply(a,b)$ 
are defined,  the laws of ring theory are quantifier-free and disjunction-free, and we wish 
to use a general metamathematical tool to ``import'' classical proofs of such formulas, thus 
avoiding the necessity of checking a long proof for constructivity line-by-line.  The tool 
to be used is G\"odel's double-negation interpretation.  
 We therefore need only check that there is 
a classical proof, using axioms of geometry to which the double-negation interpretation applies.
The double-negation intepretation has been worked out for  Tarski's
axioms for geometry in \cite{beeson2015b}; using that result, we need only to check 
that the laws of ring theory are classically provable in Tarski's geometry.  That is done 
in \cite{schwabhauser}, for Descartes's multiplication; so it only remains to check in 
Tarski's geometry that Hilbert's multiplication agrees with Descartes's multiplication.  
Hilbert's multiplication is not specifically mentioned in \cite{schwabhauser}, but the
equivalence is proved directly on pp.~52--54 of \cite{hilbert1899}, 
using the theorem of Pappus-Pascal,
which {\em is} proved in \cite{schwabhauser}; this proof can be carried out in Tarski's 
geometry since all the theorems mentioned are proved in \cite{schwabhauser}.  That completes 
the proof.
\medskip

\noindent{\em Remarks}.  The construction of $\HilbertMultiply(a,b)$
 works in the degenerate cases, including $a=b=0$.   It is instructive
to trace through what happens when $a=b=0$, and we recommend the reader to think about the 
case $a=b=0$ and the cases when $a$ and $b$ are small, but not zero, and possibly of unknown signs.
You will see why we need the {\em uniform} perpendicular bisector of $ab$, and why we need
{\em uniform} reflection in constructing the ``other intersection point.''

 While the verifications of the ring-theory laws are not trivial, the point of the above proof is
 that it suffices to check them {\em classically},
since once we have given uniform constructions for addition and multiplication, the statements
of the ring-theory laws become quantifier-free and negative.  The heart of the constructive
treatment of this subject is giving {\em uniform} constructions. 
\medskip

\subsection{Square roots} \label{section:squareroots}
\begin{figure}[h]
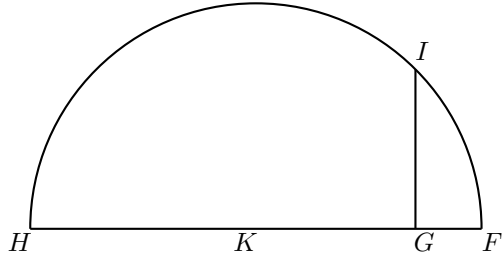
 
\center{\SquareRootFigure}
\caption{Square roots according to Descartes}
 \label{figure:SquareRootFigure}
\end{figure}

Now we take up the geometrical construction of square roots.   Fig.~\ref{figure:SquareRootFigure} shows Descartes' construction for 
finding the square root of $HG$.  His answer is the length of segment $IG$.  Here is a geometrical 
construction term (encoded as a program) to carry out Descartes' construction uniformly, as long 
as $G$ is not to the left of $H$ in the illustration,  regardless of whether $G=H$ or not.  We assume 
that $H$ in the diagram is $0$ and the horizontal line is $\Line(0,1)$.  The point 1 is not shown in 
the figure. 

\begin{verbatim}
SquareRoot(Point G) 
{ // H in Descartes' diagram is 0
  F = Add(0,G,0,1)  // so FG has unit length
  K = Midpoint(F,0)
  C = Circle(K,F)
  L = Perp(G,Line(0,F))
  I = IntersectLineCircle1(L,C)
    // next rotate unit length to line L
  U = IntersectLineCircle1(L,Circle(G,F))  
  R = Rotate(U,G,F,I) // so now RG = IG
    // but R is on Line(H,G), on the same side of G as F
    // now we need N so that N0 = RG
  MinusOne = IntersectLineCircle2(1,0,Circle(0,1))
  N = Extend(MinusOne,0,G,R)
  return N
}
\end{verbatim}

Descartes stops when he has constructed $I$.  What we have to do extra is to construct a point $N$ such 
that $0N = GI$.  In order to do that uniformly, we must not assume that $I \neq G$.  In order to get 
a non-degenerate angle we cannot use $IGF$; instead we lay off a unit length on the perpendicular $GI$,
which has been correctly constructed even if $0=H=G=I$.
Thus we do not need to assume $G \neq 0$ for this construction to work;  we only need that 
$0$ is not between $G$ and 1;  loosely speaking, $G \ge 0$.  This works because $\Perp$ is total.

We now check the validity of Descartes's construction of square roots.  Let $\SquareRoot(x)$ 
be the term defined by the construction script in Section \ref{section:squareroots}. Then we have 

\begin{Lemma} \label{lemma:HilbertDescartes} Playfair's parallel axiom implies that iff $x$ lies on $\Ray(0,1)$, and $z = \SquareRoot(x)$, then $\HilbertMultiply(z,z) = x.$
\end{Lemma}

\begin{figure}[ht]
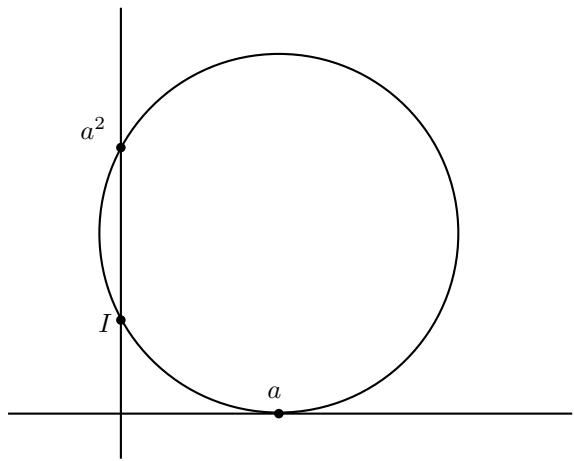
    % 
\center{\HilbertSquaringFigure}
\caption{Hilbert multiplication in the special case of squaring.}
\label{figure:HilbertSquaringFigure}
\end{figure}

\begin{figure}[h]
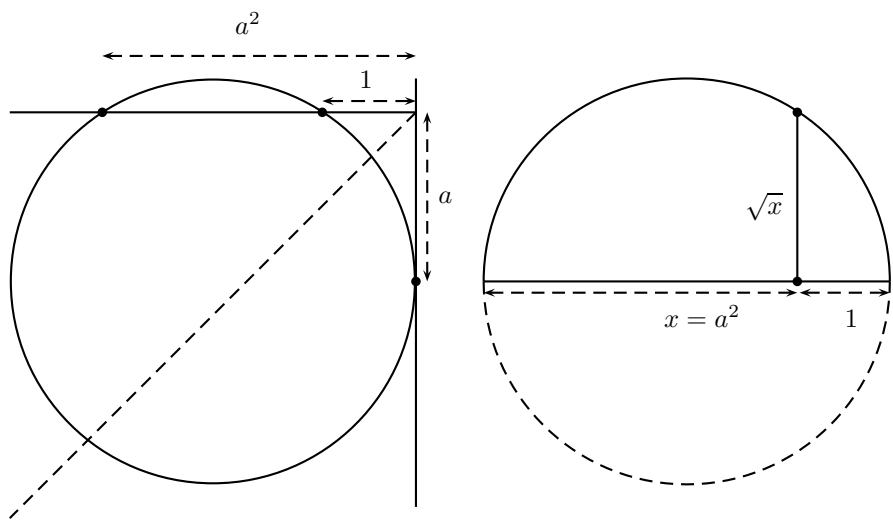
  
\center{\HilbertDescartesFigure}  
\caption{If  $z = \SquareRoot(x)$, then $\HilbertMultiply(z,z) = x$.  The left circle is Fig.~\ref{figure:HilbertSquaringFigure} reflected in the diagonal line.  The right circle is a relabeling that is recognized as Descartes's square root.}
\label{figure:HilbertDescartesFigure}
\end{figure}

\noindent{\em Proof}.  Since what is to be proved is quantifier-free, we can appeal to the 
double-negation interpretation in \cite{beeson2015b}.   But there is also a beautiful 
``proof by diagram'', which we now exhibit.    
 In the case of squaring,  Hilbert multiplication has the circle 
tangent to the horizontal axis, since $a$ and $b$ coincide, as illustrated in  
Fig.~\ref{figure:HilbertSquaringFigure}.
Replacing $a$ by $\sqrt x$ and reflecting Fig.~\ref{figure:HilbertSquaringFigure} in a diagonal line from upper left to lower right, 
we obtain the first diagram in 
 Fig.~\ref{figure:HilbertDescartesFigure}.  Relabeling that diagram, and omitting the bottom half of the circle, 
 we obtain Descartes's method of 
calculating the square root of $x$, shown in the second diagram in Fig.~\ref{figure:HilbertDescartesFigure}.  
%It is routine to turn this diagram into a rigorous proof, but since the double negation interpretation 
%already guarantees success, we omit the details.

\FloatBarrier
 
\section{The arithmetization of geometry} \label{section:arithmetization}
Our aim in this section is to lay the groundwork for a constructive version of the classical representation 
theorem:  the models of   
Euclidean geometry are exactly planes $\F^2$ where $\F$ is a Euclidean field. Classically, 
a Euclidean field is an ordered field in which positive elements have square roots. 
The groundwork in question consists in developing the theory of Euclidean fields constructively,
and verifying constructively that the axioms of line-circle continuity and circle-circle continuity hold in each plane over a
Euclidean field.

\subsection{Euclidean fields in constructive mathematics}
 We discuss the axiomatization of Euclidean fields with intuitionistic logic. 
 We use a language with symbols $+$ for addition and $\cdot$ for multiplication, and a unary predicate $P(x)$ for 
``$x$ is positive''.  
We take the usual axioms for fields, except the axiom for multiplicative inverse, which 
says that positive elements have multiplicative inverses.  If positive elements have inverses, it is an easy exercise 
to show that negative elements do too.  We define a {\em Euclidean field} to be a commutative ring satisfying the 
following additional axioms:
\begin{eqnarray*}
&\neg x \neq y \implies x = y&\mbox{\qquad stability of equality} \\
&0 \neq 1&  \mbox{\qquad EF0} \\
&x \neq 0 \implies \exists y\,( x \cdot y = 1) &\mbox{\qquad EF1}  \\
&P(x) \land P(y) \implies P(x+y) \land P(x\cdot y)  &\mbox{\qquad EF2} \\
&x+y = 0 \implies \neg(P(x) \land P(y))   &\mbox{\qquad EF3}\\
&x+y=0 \land \neg P(x) \land \neg P(y) \implies x = 0  & \mbox{\qquad EF4} \\
& x + y = 0 \land \neg P(y) \implies \exists\, z( z \cdot z = x) &\mbox{\qquad EF5}\\
& \neg \neg P(x) \implies P(x) &\mbox{\qquad EF6, or Markov's principle} 
\end{eqnarray*}

As usual, we define $x < y$ to mean $\exists z(P(z) \land x+z = y)$, or informally, $y-x$ is positive;
and $x \le y$ means $\neg(y < x)$.   Then Markov's principle is equivalent to $\neg (x \le 0) \implies 0 < x$.  It can also be written as $\neg \neg 0 < x \implies 0 < x$.
  One can extend theories including EF2 and EF3 conservatively by adding symbols $x<y$, $x \le y$, $-x$,
$x-y$ and $\sqrt x$, and then define $\vert x \vert = \sqrt{x^2}$.  Informally we make use of 
these symbols, but they are eliminable if desired.  We also informally abbreviate $x\cdot x$ as $x^2$.

Axiom EF5 says that non-negative elements have square roots.   This is a stronger axiom, intuitionistically,  than 
simply specifying that positive elements have square roots.  Using the abbreviations mentioned above,
EF5 can be 
  written $\neg P(-x) \implies \exists z\, (z^2 = x)$, or  as $0 \le x \implies \exists z\, (z^2 = x)$,
  or as $0 \le x \implies (\sqrt x)^2 = x$.  
  
We note that ``all positive elements have square roots'' implies ``all nonzero elements have 
square roots'',  in view of the fact that 
$$ \frac 1 \alpha = \alpha \cdot \bigg(\frac 1 {\alpha^2}\bigg)$$
That is, if $y \cdot \alpha^2 = 1$, then $(\alpha y)\cdot \alpha = 1$.  Hence, contrary 
to \cite{beeson2012}, it is not a different theory of Euclidean fields if we require 
only that all positive elements have square roots.

In commutative ring theory an element $x$ is called {\em invertible} if it has a multiplicative inverse; that is, for some $y$ we have $x\cdot y = 1$. We also say ``$x$ has a reciprocal.''   The weakest version 
of Euclidean field theory that we consider replaces Axiom EF1 with the following two axioms, which 
say that elements without reciprocals are zero, and that ``bounded quotients exist.''
\begin{eqnarray*}
& (\forall y (x \cdot y \neq 1)) \implies x = 0& \mbox{\qquad EF7} \\
& \vert a \vert \le \vert b \vert y \implies \exists z(a = bz) & \mbox{\qquad EF8, bounded quotients exist}
%& P(y) \land P(z) \land y+z = x  \land (y\cdot v = 1) \implies \exists\, w\, (w\cdot x = 1) &\mbox{\qquad EF8}
\end{eqnarray*}
Fields that satisfy EF0, EF2--EF10 are called ``Playfair rings'', because (as shown
in \S\ref{section:fieldtheory}) they correspond to models of the Playfair parallel axiom.  EF7 enables us to verify the 
Playfair axiom, and EF8 enables us to verify Pasch's axiom and circle-circle continuity.

To recap:
In a Euclidean field, all nonzero elements have multiplicative inverses,
In a Playfair ring, elements without reciprocals are zero,
and $a/b$ exists if it is bounded; that is,
 if for some $y$ we have $\vert a \vert \le \vert  b \vert y$.  
With classical logic, all these concepts coincide, since EF7 is classically equivalent to EF1.

We check that some common calculations still work without EF1:

\begin{Lemma} \label{lemma:Playfaircancel} 
In a Playfair ring the following hold:
\smallskip

(i) If  $x > 0$ and $y < 0$, then $xy < 0$.

(ii)  $\neg \neg\, x \le y  \implies x \le y $  (stability of $\le$).

(iii) If $ab \le ac$, and $0 < a$, then $b \le c$.

(iv) If $0 < a$ and $ab=0$ then $b=0$.
\end{Lemma}

\noindent{\em Proof}. Ad (i):  Suppose $x > 0$ and $y < 0$.  Then $P(-y)$ and $P(x)$, so by EF2,
$P(x \cdot(-y)$.  But $x \cdot (-y) = -xy$, so $P(-xy)$.  But that is, by definition,
$xy < 0$.  That proves (i).  

Ad (ii):   Suppose $\neg\neg\, x \le y$.   That is, $\neg \neg(\neg P(x-y))$.  Triple 
negation is the same as negation, so $\neg P(x-y)$.  That is, $x \le y$.  

 Ad (iii):  Suppose $0 < a$ and $ab \le ac$.  We must prove $b \le c$.   By (ii), it suffices to 
derive a contradiction from the assumption  $\neg b \le c$, that is, from $\neg \neg \,c < b$.  
For that it suffices (even without EF6) to derive a contradiction from $c < b$.
Suppose $c < b$; then 
$c-b < 0$, so by (i), $a(c-b) < 0$, contradicting $ab < ac$. 

Ad (iv):  $ab \le ab$ implies $ab/a$ exists (and is $b$),
 so we can divide both sides of $ab=0$ by $a$, obtaining $b=0$. That completes the proof. 
\smallskip

{\em Remark}. In a Playfair ring, there might be zero divisors ($ab=0$ with neither $a=0$ nor $b=0$),
but nevertheless by (iv) we can cancel positive factors.

\subsection{Line-circle continuity over Playfair rings} \label{section:continuityovereuclideanfields}

Let $\F$ be an ordered field.  Then ``the plane over $\F$'',  denoted by $\F^2$, is a geometrical structure
determined by defining relations of betweenness and equidistance in $\F^2$, using the given order of $\F$.
Namely, $b$ is between $a$ and $c$ if it lies on the interior of the segment $ac$.  That
relation can be expressed more formally in various ways, for example (using the cross product) 
by  $(c-b) \times (b-a) = 0$  and $(b-c) \cdot(a-b) > 0$.  Similarly, non-strict betweenness means that 
$b$ lies on the closed segment $ac$, or formally,  $(c-b) \times (b-a) = 0$  and $(b-c) \cdot(a-b) \ge 0$.

 The equidistance relation $E(a,b,c,d)$,
which means that segment $ab$ is congruent to segment $cd$,  can be defined by $(b-a)^2 = (d-c)^2$.
Note that no square roots were used, so these definitions are valid in any ordered field.  Though we 
have not given specific axioms for betweenness and congruence, it is reasonable to demand that 
any constructive axioms for betweenness and congruence should hold in $\F^2$.
By line-circle continuity
we understand this axiom:  let point $a$ be non-strictly inside circle $C$, and point $b$ be non-strictly outside
circle $C$.  Then $\Line(a,b)$ meets circle $C$ in a point non-strictly between $a$ and $b$.   By 
circle-circle continuity,  we mean that if circle $C$ has a point non-strictly inside circle $K$, and 
another point non-strictly outside circle $K$, then there is a point on both circle $C$ and circle $K$.
 
\begin{Theorem} [proved constructively] \label{theorem:hartshorne}
Let $\F$ be a Playfair ring;  that is, an ordered Euclidean ring in which bounded quotients exist
and elements without reciprocals are zero (EF7 and EF8 above).  Then
\smallskip

(i) $\F^2$ satisfies circle-circle continuity. 

(ii) $\F^2$ satisfies line-circle continuity.

(iii) Nonnegative elements have square roots in $\F$.
\end{Theorem}

\medskip

\noindent{\em Proof.}  A similar theorem (for classical Euclidean fields)
 is stated  in \cite{hartshorne}, p.~144, and a
 proof is sketched.  Here we have additional issues: we must use only 
 intuitionistic logic, and we must get by with weaker 
 assumptions about the existence of reciprocals.  In proving (i) implies (ii), no division 
 is required.  In proving (ii) implies (iii), we just use Descartes's square-root construction.
 Consider proving that (iii) implies (i).  We can assume (using translations and rotations 
 without any division) that one circle (call it $C$) is defined by $x^2 + y^2 = r^2$ and the other 
 (call it $K$) is defined by $(x-c)^2 + y^2 = R^2$, with $c > 0$.
   Subtracting, we have $2cx = r^2 - R^2 + c^2$, so 
 we only need to be able to divide the right side by $2c$.  Using EF8,  it suffices to 
 bound the quotient, that is, to bound the $x$-coordinate of the intersection points.  
 That seems quite intuitive:  $x$ must lie between $-r$ and $c+R$.  But to prove it, 
 we must use the hypotheses of circle-circle continuity.
One proves (by approximately one page of algebraic computations, here omitted)
 that $C$ has a point (non-strictly) inside $K$ and another (non-strictly) outside $K$
if and only if the following two inequalities hold:
\begin{eqnarray}
 r+R &\ge& c  \label{eq:circles1} \\
 R-c &\le& r \ \le \ R+c \label{eq:circles2}
 \end{eqnarray}
This is the algebraic expression of 
the hypotheses of circle-circle continuity.  
(These are the very details that are ``left to the reader'' in \cite{hartshorne}, p.~144.)
The conclusion is very believable without a proof, since (\ref{eq:circles1}) fails when the circles 
bound disjoint closed disks, and (\ref{eq:circles2}) fails if one lies inside the other. 
Inequality (\ref{eq:circles2}) implies that $\vert R-r \vert\le c$; this is obvious if $r \le R$
or $R \le r$, and since $\neg\neg\, (r \le R \lor R \le r)$, we have $\neg\neg\, \vert R-r \vert \le c$,
but then by the stability of $\le$ we can drop the double negation. 

Returning to the proof, in order to solve the equation $2cx = r^2-R^2 + c^2$ for $x$
in a Playfair ring using EF8, we need to justify dividing the right side by $2c$.  By EF8,
it suffices to show that for some $z$, we have
 $$\vert r^2-R^2 + c^2\vert \le 2cz $$
We take $z=r+R+c$.  Then the required inequality is  
  proved as follows:
\begin{eqnarray*}
\vert r^2- R^2 + c^2 \vert &\le& \vert r^2-R^2 \vert + c^2 \\
&\le& \vert (r-R) \vert\,  (r+R) + c^2 \\
&\le& (r+R)c + c^2  \mbox{\qquad since $\vert r-R \vert \le c$}\\
&\le& (r+R+c) c \ = \ zc \ \le \ 2zc
\end{eqnarray*}
Hence the $x$-coordinate of the intersection points exists.  Then 
 we also need to prove $y = \sqrt{r^2-x^2}$ is defined, i.e., $x^2 \le r^2$. 
 By (\ref{eq:circles2}), we have $ r \le R+c$. Hence $r-c \le R$.  If 
 $c \le r$ we have $\vert r-c \vert \le R$.  On the other hand, if $r \le c$
 then (\ref{eq:circles1}) implies $\vert c-r \vert \le R$.  Since $ \neg \neg\, (r \le c \lor c \le r)$, 
 we have $\neg\neg\, \vert r-c \vert \le R$, and hence by the stability of $\le$ we have
 $\vert r-c \vert \le R$. Hence
\begin{eqnarray*}
(r-c)^2 &\le& R^2 \\
r^2-2rc + c^2 &\le& R^2 \\
r^2 - R^2 + c^2 &\le& 2rc \\
2cx &\le& 2rc \mbox{\qquad since $2cx = r^2-R^2+c^2$} \\
x &\le& r  \mbox{\qquad \quad by Lemma~\ref{lemma:Playfaircancel} (iii), since $c>0$.}
\end{eqnarray*}
Now if $0 \le x$ we can conclude $x^2 \le r^2$ as desired. 

Next we claim $\vert R - c \vert \le r$.  If $c \le R$,  that follows from (\ref{eq:circles2}).
If $R \le c$, it follows from (\ref{eq:circles1}).  Since $\neg \neg\,(c \le R \lor R \le c)$,
the claimed inequality follows by the stability of $\le$.   Hence
\begin{eqnarray*}
\vert R-c \vert &\le& r \\
(R-c)^2 &\le& r^2 \\
R^2 - 2Rc + c^2 &\le& r^2 \\
-2Rc + 2c^2 &\le& r^2 -R^2 + c^2 \ = \ 2cx 
\end{eqnarray*}
As shown above it is legitimate to divide by $2c$; we obtain
\begin{eqnarray*}
-R + c &\le& x \\
R-c &\ge& -x
\end{eqnarray*}
If $x \le 0$ then we can square both sides, obtaining 
$ x^2 \le (R-c)^2$.  But $(R-c)^2 \le r^2$; hence $x^2 \le r^2$ (still under 
the assumption $x \le 0$).

Now we have proved $x^2 \le r^2$ follows from $0 \le x$, and also from $x \le 0$.
Since $\neg\neg\,0 \le x \lor x \le 0$,  we conclude by the stability of $\le$ that
$x^2 \le r^2$ holds without assumptions.  Hence $y = \sqrt{r^2-x^2}$ is defined.
That completes the proof.
\medskip

{\em Remark}.  We have proved that the existence of ``bounded quotients'' is enough
to guarantee the existence of intersection points of circles.   That corresponds intuitively
to the fact that these intersection points are known in advance to lie within the bounds of 
the two circles.  By contrast, the point of intersection of two lines asserted to exist
by Euclid 5 is not bounded in advance, unless we have a theory of 
similar triangles, but we need Euclid 5 to develop such a theory.  If we consider a line
in $\F^2$ passing through $(0,1)$ with slope $t$,  it will meet the $x$-axis at $(1/t,0)$.
But without an {\em a priori} bound on $1/t$,  we can't prove that point exists on the 
basis of the axioms of Playfair rings.

\subsection{Pasch's axiom in planes over Playfair rings}  
We first consider the use of analytic geometry to calculate the intersection point
of two lines in a plane over a Playfair ring $\F$.
 In $\F^2$,  if $\Line(u,v)$ and $\Line(s,t)$ intersect (and do not coincide),
 the formula for the 
 intersection point has denominator given by the cross product $D=(t-s) \times (v-u)$.
 To prove this, let $x$ be the intersection point.  Then $x$ satisfies the equations of 
 both lines:
 \begin{eqnarray*}
 (x-s) \times (t-s) &=& 0 \\
 (x-u) \times (v-u) &=& 0 
 \end{eqnarray*}
 Then 
 \begin{eqnarray*}
 x \times (t-s) &=& s \times (t-s) \ = \ s \times t\\
 x \times (v-u) &=& u \times (v-u) \ = \ u \times v
 \end{eqnarray*}
 Writing out the cross product on the left, we obtain linear equations for the 
 components $x_1$ and $x_2$ of $x$.   Cramer's
 rule yields expressions for $x_1$ and $x_2$ with $D$
 in the denominator, as claimed above; so if $D$ is invertible, we can 
 solve for $x$.   Even if $D$ is not 
 invertible, we can still use Cramer's rule to give a formula for $x_1$ and $x_2$
 as long as we know in advance that the lines intersect in $x$.   That is,
 if Cramer's rule gives the formula $x_1 = R/D$, then at least we know $Dx_1 = R$,
 so $R \le Dx_1$, so $R/D$ exists by EF8.  
 
 On the other hand, suppose we do not know that the lines intersect.  Then we
 can still calculate the formula that Cramer's rule would give, say $x_1 = R/D$
 and $x_2 = S/D$.  If we can show that these quotients are bounded, i.e., $R \le Dy$
 for some $y$ and $S \le Dz$ for some $z$,  then $R/D$ and $S/D$ are defined in $\F$,
 and we can verify that the formulas given by Cramer's rule do provide a solution.
 In essence, then, to prove two lines intersect we can proceed as follows:
 solve the equations for the intersection point as if $\F$ were a Euclidean field;
 having arrived at expressions for the solutions, show they are bounded by something 
 that we know exists in $\F$;  then the formal quotients defining the solution actually
 exist.
 
\begin{figure}[ht]
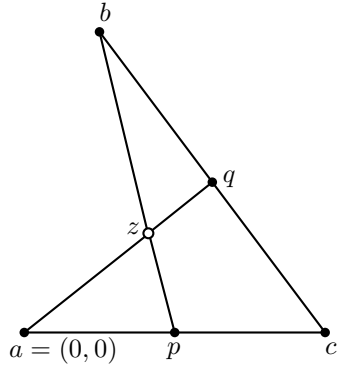
   
\center{\VerifyingPaschFigure}
\caption{Verification of   Inner Pasch}
\label{figure:VerifyingPaschFigure}
\end{figure}

We now turn to the verification of Pasch's axiom in planes $\F^2$ over a Playfair ring $\F$.
Specifically, we consider Pasch's axiom in the form ``inner Pasch''. See
Fig.~\ref{figure:VerifyingPaschFigure}.

\begin{Theorem} \label{theorem:PlayfairPasch}
Let $\F$ be a Playfair ring.  Then the plane $\F^2$ satisfies inner Pasch.
\end{Theorem}

\noindent{\em Proof}.
The hypothesis of inner Pasch is
 $$\B(a,p,c)\land\,\B(b,q,c) \land p \neq b \land q \neq a$$
 We must verify the existence in $\F^2$ of a point $z$ such that $\B(a,z,q)$ and $\B(b,z,p)$.
 That is, we must show how to calculate the coordinates of $z$,  without assuming that 
 all nonzero or positive elements have reciprocals.   We want to simplify the 
 algebra by first applying an affine transformation to bring the points to some standard position.
Although affine transformations preserve betweenness (even though they may not preserve distance),
the nature of the argument requires that we use {\em invertible} transformations,
 because we are arguing that if we can verify Pasch in the transformed situation, then we 
 can apply the inverse transform to construct the point required by Pasch in the original situation.
 
 We suppose $p = (p_1,p_2)$ with $p_1$ and $p_2$ in $\F$, and we use a similar notation for the
 coordinates of the other points. 
 By a translation, we can bring $a$ to $(0,0)$.  All translations are invertible.
   Now we want to bring $p$ to a point 
 on the positive $x$-axis.  Let $d=p_1^2 + p_2^2$; then $d \ge 0$ so $\sqrt d$ exists.
 Since $p \neq a$ and $a = (0,0)$, we have $p > 0$ by Markov's principle EF6.  
 We claim that $p_1/d$ and $p_2/d$ exist; that follows from the bounds $p_1/d \le 1$ 
 and $p_2/d \le 1$; that is, $p_1 \le d$ and $p_2 \le d$, which in turn follow from 
 $p_1^2 \le p_1^2 + p_2^2$ and $p_2^2 \le p_1^2 + p_2^2$.  Then
   by the transformation with the matrix
 $$ \matrixtwo {p_1/d} {p_2/d} {p_2/d} {-p_1/d} $$
 we can bring $p$ to the point $(d,0)$ on the positive $x$-axis.   
 This transformation is invertible, since its determinant is 1.  
 Because of the reduction by this transformation, we can assume 
 without loss of generality that $p_2 =0$.  Then because $\B(a,p,c)$, $c$ lies on the $x$-axis. 
 
 Fig.~\ref{figure:VerifyingPaschFigure} illustrates the situation. If $\F$ were a field, we 
 could arrange for $b$ to have the same first coordinate as $p$ 
using  another linear transformation, but we cannot show that the 
 required transformation is invertible, so we proceed without knowing anything about the 
 position of $b$.  We proceed as described above, to solve for $z$ as if $\F$ were a field.
 We have
 \begin{eqnarray*}
 z \times q &=& 0 \\
 (z-p) \times (b-p) &=& 0
 \end{eqnarray*}
 The first equation becomes $z_2 q_1 = z_1 q_2$, and the second 
 becomes 
 \begin{eqnarray*}
 z \times(b-p) &=& p \times b \\
 z_1 b_2  -z_2 (b_1 - p_1) &=& p \times b \mbox{\qquad since $p_2 = 0$} 
 \end{eqnarray*}
 In order to eliminate $z_2$ we multiply first by $q_1$, and the 
 replace $z_2q_1$ by $z_1q_2$:
 \begin{eqnarray*}
 z_1(q_1b_2 - q_2(b_1-p_1)) &=&  q_1(p \times b) \\
 z_1(q \times (b-p)) &=& q_1(p \times b) 
 \end{eqnarray*}
 Now, if the divisions could be carried out, we could solve for $z_1$ and $z_2$:
 \begin{eqnarray*}
 z_1 &=& q_1 \frac{p \times b}{q \times (b-p)} \\
 z_2 &=& q_2\frac{p \times b}{q \times (b-p)} \mbox{\qquad since $z_2q_1 = z_1q_2$}
 \end{eqnarray*}
 In order to justify these divisions, showing
  that $z_1$ and $z_2$ exist, we must bound the expressions on the right
 in the last two equations. 
 We claim they are bounded by $\vert q_1\vert$ and $\vert q_2\vert$, respectively.
   It suffices to show that the fraction 
 on the right is bounded by 1; that is,
 $$\vert p \times b \vert \le \vert q \times (b-p) \vert $$
 We have $p \times b = p_1b_2 > 0$. We have
 $q \times (b-p) = q \times b + p \times q$ and
 $p \times q = p_1q_2 > 0$, and $q \times b > 0$ because $\B(b,q,c)$.
 Hence it suffices to show 
 \begin{equation} p \times b \le q \times b + p \times q. \label{eq:5215} 
 \end{equation}
 That is ``evident'' from the diagram, if we knew that cross products measured
 area in the usual way.  But technically we have to compute it.  
 Let $\alpha = (c-p) \times(q-p)$. 
 Then $0<\alpha$, since $\alpha = (c_1-p_1)q_2$, and both factors are positive.
 
 It will suffice to 
 prove the following  assertions:
 \begin{eqnarray}
 c \times b &=& p \times q + q \times b + \alpha  \label{eq:5225}\\
 c \times b &=& p \times b + \alpha + (q-p) \times (b-p) \label{eq:5226} \\
 0 &<& (q-p) \times (b-p) \label{eq:5227}
 \end{eqnarray}
 for then, subtracting (\ref{eq:5225}) from (\ref{eq:5226}),
 we get 
 $$p \times b \le p \times q + q \times b -(q-p) \times (b-p) $$  
 and then, using (\ref{eq:5227}), we obtain the desired equation (\ref{eq:5215}).
 Again, each of these assertions is ``obvious'' if we interpret the cross products
 as areas of triangles, but we must give an algebraic proof.

 To prove (\ref{eq:5227}): We have
 \begin{eqnarray*}
 (q-p)\times (b-p) &=& (q_1-p_1)b_2 -q_2(b_1-p_1) 
 \end{eqnarray*}
 and the right side is positive because $q_1-p_1 > b_1-p_1$ and $b_2 > q_2$,
 both of which follow from $\B(b,q,c)$.
 
 Since $\B(b,q,c)$, we have 
 \begin{eqnarray*}
 0 &=&(c-b)\times (q-b) \\
 &=& (c_1-b_1)(q_2-b_2) + b_2(q_1-b_1) \\
 c_1 b_2 &=& c_1q_2 - b_1q_2 + b_2q_1 \\
 c \times b &=& c_1q_2 - b_1q_2 + b_2q_1 \\
 &=& c_1q_2 + q \times b
 \end{eqnarray*}
 That is the left side of (\ref{eq:5225}).  Now we calculate the right side:
 \begin{eqnarray*}
p \times q + q \times b + \alpha &=& p \times q + q \times b + (c-p) \times (q-p)\\
&=&  p_1q_2 + q_1b_2 - q_2b_1 + c_1q_2 - p_1q_2  \\
&=&q_1b_2- q_2b_1  +c_1q_2  \\
&=&c_1q_2 + q \times b
 \end{eqnarray*}
which is equal to what we found for the left side.  That proves (\ref{eq:5225}).

It remains to prove (\ref{eq:5226}).  The right side of (\ref{eq:5226}) is 
\begin{eqnarray*}
& &p \times b + \alpha + (q-p) \times (b-p)   \\
&=&p \times b + (c-p) \times(q-p) + (q-p) \times (b-p) \\
&=& p_1b_2  + (c_1-p_1)q_2 + (q_1-p_1)b_2 - q_2(b_1-p_1) \\
&=& p_1b_2 + c_1q_2 - p_1q_2 + q_1b_2 - p_1b_2 - q_2b_1 + q_2 p_1 \\
&=& c_1q_2 + q \times b \\
&=& c \times b  \mbox{\qquad as shown   above}
\end{eqnarray*}
That proves (\ref{eq:5226}), and completes the proof that $z$ is defined.

  It remains to verify the conclusion of the axiom, namely
  $\B(p,z,b) \land \, \B(q,z,a)$.  By construction $z$ lies on $\Line(p,b)$
  and $\Line(a,q)$.  Since $a_2 = p_2 = 0$, it suffices to show $0 < z_2 < b_2$
  and $0 < z_2 < q_2$.  Recall
  $$ z_2 = q_2 \frac {p \times b}{q \times(b-p)}$$
  and just before (\ref{eq:5215}) we proved that the numerator and denominator
  are both positive; hence $0 < z_2$;  and we proved that the fraction is less 
  than one, so $z_2 < q_2$.   Because $\B(b,q,c)$ we have $q_2 < b_2$.  
  Hence $z_2 < b_2$ as well.  That completes
 the proof of the theorem.

\subsection{The representation theorem}
In this section we will prove that every model of neutral constructive geometry
plus Euclid 5  and line-circle continuity is a plane over a Euclidean field.

\begin{Theorem}[Representation Theorem] \label{theorem:representation}
Euclid 5 implies that the field operations defined by $\Add$ and 
$\HilbertMultiply$ turn the $x$-axis 
into a Euclidean field; and every point arises as $\MakePoint(a,b)$ for some 
$a$ and $b$.   Hence (constructively) every model of neutral constructive geometry 
plus Euclid 5 is isomorphic to a (constructive) Euclidean field.  In particular, 
the points on a given line are exactly the points $\MakePoint(x,y)$ where $x$ and $y$
satisfy a linear equation.   
\end{Theorem}

\noindent{\em Proof}.  By Lemma~\ref{lemma:coordinatization}, every point arises as $\MakePoint(a,b)$.
By Theorem~\ref{theorem:signedmultiplication}, $\HilbertMultiply$ is defined. 
By Lemma~\ref{lemma:ringlaws}, the laws of commutative rings are satisfied.  
By Lemma~\ref{lemma:HilbertDescartes}, nonnegative elements have square roots.  
By Lemma~\ref{lemma:reciprocal}, positive elements have reciprocals.  But this implies
that all elements have reciprocals, since  $1/t = t/t^2$.  More explicitly, if $y$
is the reciprocal of the positive element $t^2$ then $ty$ is the reciprocal of $t$.
Hence all the laws of Euclidean field theory are satisfied.  Now suppose $L = \Line(a,b)$
is a line.  Then $a = \MakePoint(a_1,a_2)$ and $b = \MakePoint(b_1,b_2)$, and the points 
on the line are just the points $x$ such that the cross product $(x-a) \times (b-a)$ is zero, that is,
$(x_1-a_1)(b_2-a_2)-(x_2-a_2)(b_1-a_1) = 0$.  By the stability of equality, we are allowed to prove 
that by cases.  If $L$ is not vertical, then each vertical line $x_1 = p$ cuts $L$ exactly once,
and since there is a point with $x_1 = p$ satisfying the equation, we are done.  If $L$ is 
vertical, then $b_1-a_1 =0$ and $b_2-a_2 \neq 0$, so the equation becomes $x_1=a_1$, and indeed
the points of $L$ are exactly the solutions of that equation.
That completes the proof.

\begin{Corollary} \label{lemma:Euclid5suffices} Euclid 5 implies the strong
parallel postulate in neutral constructive geometry.
\end{Corollary}

\noindent{\em Proof}.  Let line $L$ and point $p$ not on $L$ be given.  We may
choose $L$ as the $x$-axis and $p$ as the point $I$ in setting up coordinates.
Thus $0$ is the foot of the perpendicular from $I$ to $L$, and $1$ is a point on $L$
such that $0I = 01$.  
Let line $K$ be perpendicular to the $y$-axis at $I$, and let $M$ be any other 
line through $I$.  We must prove that $M$ meets $L$.  According to the Representation
Theorem, the points on $M$ are exactly the points $(x,y)$
 satisfying some linear equation $dx + by = c$.  Here $(x,y)$ abbreviates
 $MakePoint(x,y)$.  Specifically, the proof of the 
 representation theorem yields $d=(a_2-1)$ and $c=a_1$, where $a=(a_1,a_2)$ is not on 
$K$, so $a_2 \neq 1$ and $d \neq 0$.  
or more officially, $\Add(\HilbertMultiply(d,x),\HilbertMultiply(b,y)) = c$.
Setting $y=0$ and solving for $x$ we find $x = c/d$, or officially
$x = \HilbertMultiply(c,\Reciprocal(d))$.  Then $x$ is the desired intersection
point of $L$ and $M$.
That completes the proof.
\smallskip

{\em Remark}.  The proof works even if $M$ is perpendicular to $L$, so it provides a 
continuous construction for the intersection point as $M$ rotates through all angles
(except of course when $M$ is parallel to $L$).  The proof  provides a 
geometrical construction proceeding from the two points determining $M$
to the coefficients of the linear equation for $M$ and from there to $x$.  The construction 
is explicit in the proof,  and it is not difficult to see exactly what is involved, and 
the key role that is played by uniform rotation.  Namely, $x = c/d = a_1/(a_2-1)$, 
which is constructed by drawing a circle through $I$, $(0,a_1)$, and $(a_2-1,0)$, and finding 
the other intersection point with the $x$-axis. The points $(0,a_1)$ and $(a_2-1,0)$ have
to be constructed by uniform rotation, since we are not certain whether $a_1 =0$ or $a_2 = 1$.

\begin{Corollary} \label{lemma:uniformreflectioninpoint} 
Euclid 5 implies the existence of a uniform construction of the reflection of point $x$ 
in point $p$.  
\end{Corollary}

\noindent{\em Remarks}.  The proof uses Euclid 5. 
 Is there a uniform reflection construction provable in 
neutral geometry?  
\medskip

\noindent{\em Proof}.  To construct the reflection of point $x$ in point $p$, set up 
coordinates on any line $L$ containing $p$, with $p$ at $(0,0)$.  Then the reflection 
of $x$ in $p$ has coordinates $(-x_1,-x_2)$.  The construction in question is to use the 
uniform perpendicular construction to project $x$ on each of the two axes; then use
uniform reflection in the axes to get $a =(0,-x_2)$ and $b=(-x_1,0)$.   Then there 
is a point whose projections on the axes are $a$ and $b$; that is the desired point.
It is not exactly $MakePoint(a,b)$ since $MakePoint$ wants arguments on the $x$-axis.
Technically, we first rotate $a$ to $c = (-x_2,0)$.  Then the desired point is  
$MakePoint(c,b)$ is the desired point.  That completes the proof.

\begin{Corollary} \label{lemma:stability} [Stability of definedness]  In elementary 
constructive geometry (that is, with circles and Euclid 5),  every term $t$ for a 
ruler and compass construction satisfies $ \neg \neg t \defined \implies t \defined$.
\end{Corollary}

\noindent{\em Proof}.  By induction on the complexity of the term $t$.   The problematic 
case is when $t$ has the form $\IntersectLines(L,M)$.  By Lemma~\ref{lemma:linestability},
this case follows from the strong parallel axiom, and hence by Corollary~\ref{lemma:Euclid5suffices},
from Euclid 5. That completes the proof.

\subsection{Models and interpretations}
When studying independence results relative to theories based on classical logic,
we have the completeness theorem. That enables us to replace proof theory by 
model theory in some situations.  For example, to transfer results about real-closed
fields to results about Tarski's geometry (with first-order completeness) we use the fact
that the models of this geometry are all of the form $\F^2$ for $\F$ a real-closed field.
Constructively, we can often convert a model-theoretic argument to a proof-theoretic
{\em interpretation}.  An interpretation is a map that associates to each formula $\phi$
in one theory, another formula $\bar \phi$ in another theory, preserving provability.
Sometimes one has an inverse interpretation, going ``the other way''.   When the two 
theories have different languages, the technical details are often intimidating, or 
at least lengthy, and are seldom written out.  

In the case at hand,  the fact is that there is an interpretation from almost 
any reasonable theory of  Euclidean geometry 
  to a corresponding theory in the language of Euclidean fields.   
Each ``point'' variable becomes a pair of variables $(x,y)$, and one uses the 
algebraic definitions of collinearity and betweenness to interpret the predicates of 
incidence of a point on a line, and betweenness.  Equidistance (or segment congruence)
is interpreted using the ``distance formula.''    We write the resulting
interpretation of the formula $\phi$ as $\bar \phi$.  Essentially $\bar \phi$ is 
nothing but ordinary analytic geometry.  

There is one technical detail that should be mentioned.  In an interpretation, it
is convenient if there are ``enough terms'' in the target theory to interpret the 
terms of the source theory.  For example, if we have function symbols for the 
intersection points of circles and lines,  since the coordinates of those points 
involve square roots, it will be convenient to formulate Euclidean field theory with 
a symbol for the square root.   We can of course add symbols for the square root and 
multiplicative inverse, conservatively over the theories as formulated above; no 
matter that the square root of negative elements is not defined.  Since this 
paper is quite long as it is, we omit these details and just assume that the interpretation $\phi \mapsto \bar \phi$ can be defined
and is sound.

The reader will notice that we have not specified a precise formal theory 
for Euclidean geometry.  Indeed one of our points is that we have a lot of latitude
to make choices about the details of the theory.  We now suppose, however, that 
``neutral \ECG''  is a theory of Euclidean constructive geometry (with no 
parallel axiom) such that 
the interpretation $\phi \mapsto \bar \phi$ is sound into 
the theory of Playfair rings.   This is true of the version 
of \ECG\  presented in \cite{beeson-kobe} and the version mentioned in \cite{beeson2012}
and presented in detail in \cite{beeson2015b}, and of many variations on these theories.

The following theorem 
states the soundness of this interpretation for two versions of the parallel 
postulate:

\begin{Theorem} \label{theorem:soundnessofbar}
(i)  If $\phi$ is provable in \ECG\ with Playfair's axiom 
instead of the strong parallel postulate, then $\bar \phi$ is provable
in the theory of Playfair rings.

(ii) if $\phi$ is the strong parallel postulate, then $\bar \phi$ is equivalent
to the statement that nonzero elements have multiplicative inverses (using the 
other axioms of Euclidean field theory).
\end{Theorem}

{\em Remark}. The theorem states more than just the faithfulness of the interpretation.
Faithfulness would only require that $\bar \phi$ be provable from, rather than equivalent to,
the corresponding field theory axiom in (ii) and (iii).
\medskip

\noindent{\em Proof}.  We are assuming that $\bar\phi$ is provable in the
theory of Playfair rings for all the axioms of neutral $\ECG$.  It remains to 
check the interpretations of the parallel axioms.

{\em Ad (i)}.  Recall that Playfair's axiom hypothesizes a line $L$ (which we may
interpret as the $x$-axis),  a point $p$ not on $L$ (which we may interpret as $(0,b)$),
and  two lines $K$ and $M$ through $p$, that are both parallel to $L$.  The conclusion
is that $K$ and $M$ coincide.  We may assume without loss of generality that 
one of the lines, say $K$, makes a right angle with the perpendicular from $p$ to $L$,
so all points on $K$ have the same $y$-coordinate $b$.
Let $(u,v)$ be another point on $M$.  We must show $v=b$.  
 If $v-b$ has a multiplicative inverse,
then ordinary analytic geometry allows us to compute the coordinates of a point 
that lies on $M$ and $L$, contradiction.  Hence $v-b$ has no inverse.  Hence, by
Playfair's axiom in field theory, $v-b=0$. That completes the proof of (i).

{\em Ad (ii)}.  We  coordinatize the diagram in 
Fig.~\ref{figure:StrongParallelRawFigure}. 
We may assume $s = (0,0)$ and $L$ is the $x$-axis.  Let the equation of $M$ be $Ax+By=C$.  
Setting $y=0$ we find $x=C/A$, so $M$ and $L$ will meet if and only if $A$
has a multiplicative inverse.  Since the point $a$ on $M$ is assumed not to lie on $K$, $A \neq 0$.
Hence, the axiom in field theory that nonzero elements have inverses
is sufficient to prove $\bar \phi$.   Conversely, given a nonzero field element $A$,
we take $s=(0,0)$, $q=(2,0)$, $t = (1,1)$, $p=(0,2)$, $r=(2,2)$, and $a = (2,2-4A)$.
The equation of $M$ is then $2Ax + y = 2$.   Since $1 < 2$,  $1/2$ is a bounded quotient,
so $2$ has a multiplicative inverse in every Playfair ring.   
Assume $\bar \phi$.  Then $K$ and $L$ meet; hence there is an $x$ such that 
$2Ax = 2$. Since 2 has a multiplicative inverse, we have $Ax=1$, so $x$ is the inverse of $a$. 
That completes the proof of the lemma. 

We can also define an ``inverse'' interpretation $\phi^\circ$ from Euclidean field
theory to geometry.   The details of this interpretation and its
soundness are lengthy and technical, and the result is not needed anyway for our independence proofs,
so we omit those details (although we have in fact written them out).

  \section{Independence  in Euclidean field theory} \label{section:fieldtheory}
One might also ask what the field-theoretic version of Playfair's axiom is.  Recall that Playfair says, if $p$ is not 
on $L$ and $K$ is parallel to $L$ through $p$, that if line $M$ through $P$ does not meet $L$ then $M=K$.
Since $\neg \neg M=K \implies M=K$, Playfair is just the contrapositive of the strong parallel postulate,
which says that if $M \neq K$
then $M$ meets $L$.  Hence it corresponds to the contrapositive of $x \neq 0 \implies 1/x \defined$; that 
contrapositive says that if $x$ has no multiplicative inverse, then $x = 0$.   Thus Playfair geometries correspond 
to ordered fields in which elements without multiplicative inverses are zero.  Of course, we need
{\em bounded quotients} to exist, to make Pasch's axiom be satisfied in $\F^2$, but the Playfair 
axiom is satisfied in $\F^2$ when elements without multiplicative inverses are zero.

Constructively, the strong parallel postulate  implies Playfair;  we wish to show that the implication 
is not reversible.  Since classically, the implication {\em is} reversible,  we cannot hope to give a counterexample.
In terms of field theory, we won't be able to construct a Euclidean ring in which all elements
without reciprocals are zero, but not all nonzero elements have reciprocals.
We must use some tools of logic.

Among the possible techniques for proving that a constructive implication is not reversible is the method 
of Kripke models.  We shall not  explain this technique in its full generality, but only in the case
of theories based on (ordered) ring theory.  A full explanation can be found in \cite{vanDalen}, but the presentation here is self-contained.

\subsection{What is a Kripke model of ring theory?}
It will be essential to understand the notions  ``term'' and ``formula'' (of ordered ring theory) as logicians use them.
A term is, intuitively, an expression meant to denote an element of a ring.  A variable (such as $x$, $y$, etc.) 
is a term, and so are the constants 0 and 1.  If $t$ and $s$ are terms then so are $(t \cdot s)$ and $(t + s)$,
as well as $(-t)$ and $(1/t)$.  In informal usage, many parentheses are left unwritten, according to the usual 
conventions.  For example,  $(1/((x+1) \cdot (y + 1)))$ is a term.   Note that $1/0$ is technically a term; 
not every term necessarily denotes something  (``is defined'').    

Next we explain the notion ``formula''.  This notion is defined 
recursively:  formulas are built up by combining smaller formulas using the logical connectives $\land$ (and), $\lor$ (or),
$\neg$ (not), and $\implies$ (implies), as well as the quantifiers $\forall \, x$ and $\exists \, x$,
according to some rules which we will not spell out in detail.  The base case of this definition 
is the ``atomic formula'',  which is either of the form $P(t)$ for some term $t$, or of the form $t=s$ for some 
terms $t$ and $s$, or of the form $t \defined$  for some term  $t$.   The formula $t \defined$ is read ``$t$ is defined.''

   An (ordered) ring is a model of (ordered) ring theory;  a Kripke model of 
(ordered) ring theory is a more complicated thing.    It is a collection of rings, $R_\alpha$, where the 
subscripts $\alpha$ come from some partial ordering $(K, <)$.   Each $R_\alpha$ must have a notion of ``positive''; that is, 
  a subset of so-called positive elements to interpret the predicate symbol $P(x)$, but it is not required that
$R_\alpha$ be an ordered ring.   The ordering on the index set $K$ has nothing 
to do with the ordering on the ring, which is given by a predicate $P(x)$ defining the positive elements.
It is often required that if $\alpha \le \beta$, then $R_\alpha$ is a subring of $R_\beta$.  It will be 
convenient to generalize this requirement by allowing $R_\alpha$ to be embedded in $R_\beta$ by means of a one-to-one function $j_{\alpha\beta}$ (which would be the identity if $R_\alpha \subseteq R_\beta$).
These functions must compose according to the law $j_{\alpha\beta} j_{\beta\gamma} = j_{\alpha\gamma}$.
Using some abstract nonsense, we could replace each $R_\alpha$ by a suitable copy to ensure that the $j_{\alpha\beta}$
are all the identity and $R_\alpha \subseteq R_\beta$, but it is convenient not to require that.   There 
is also a requirement of ``persistence'':  if $x$ is positive in $R_\alpha$, then 
$j_{\alpha\beta}(x)$ has to 
be positive in $R_\beta$.  But note, there can also be positive elements in $R_\beta$ that do not arise in that way.

The Kripke model 
$\mathcal R$ is technically the function $\alpha \mapsto R_\alpha$ with domain $K$,  though one often thinks
of it as the collection of the $R_\alpha$.  Usually the index set $K$ has a least element, the ``root''.
The elements of $K$ are called ``nodes''. Usually (and in all the models we will use)  the set $K$ is a tree, i.e., the set of nodes less than 
any given $\alpha$ is linearly ordered.   What we have defined so far is a ``Kripke structure'';  to be a Kripke model
of (ordered) ring theory,  or of ordered field theory, the structure must ``satisfy the axioms'' of the theory.
We next define that concept.

We consider {\em valuations} $\sigma$; these are functions that assign an element 
$x \sigma$ of $R_\alpha$ to each variable $x$.  (Logicians  write valuations on the right, as $x \sigma$,
rather than $\sigma(x)$.)   
If $\sigma$ is a valuation, it starts out as a function defined on variables, but is easily extended to 
a function defined on terms.  For example, $(t + s)\sigma = t\sigma + s\sigma$, where the $+$ on the 
right is  addition in $R_\alpha$, and the $+$ on the left is just a symbol.   This extended function 
is a partial function, because $(1/t)\sigma$ is undefined if $t\sigma$ is 0.   If $t$ belongs to the 
domain of the extended valuation $\sigma$ then we say $t$ is defined in $R_\alpha$.  If $\alpha < \beta$ 
and $\sigma$ is a valuation into $R_\alpha$ and $\tau$ is a valuation into $R_\beta$ then we say $\tau$ 
agrees with $\sigma$ at $x$ if $x \sigma = (j_{\alpha\beta} x) \tau$.   

We next define the notion ``$R_\alpha$ satisfies formula $A$ under valuation $\sigma$'', which is written
$R_\alpha \models_\sigma A$.  The rules for this definition are as follows:
\begin{eqnarray*}
R_\alpha &\models_\sigma& P(x) \qquad \mbox{if and only if $x\sigma$ is positive in  $R_\alpha$ } \\
R_\alpha &\models_\sigma&  t =s  \qquad \mbox{if and only if terms $t\sigma$ and $s\sigma$ are equal (equivalent) elements of $R_\alpha$} \\
R_\alpha &\models_\sigma&  t\defined  \ \ \ \qquad \mbox{if and only if $t\sigma$ is defined in $R_\alpha$} \\
R_\alpha &\models_\sigma& (A \land B) \quad \mbox{if and only if $R_\alpha \models_\sigma A$ and $R_\alpha \models_\sigma B$}\\
R_\alpha &\models_\sigma& (A \lor B) \quad \mbox{if and only if $R_\alpha \models_\sigma A$ or $R_\alpha \models_\sigma B$}\\
R_\alpha &\models_\sigma& (\neg A) \qquad \mbox{ if and only if for all $\beta \ge \alpha$  and valuations $\tau$ on $R_\beta$, } \\
&&\mbox{\qquad\qquad \ \  if $\tau$ agrees with $\sigma$ then not $R_\beta \models_\sigma A$} \\
R_\alpha &\models_\sigma& (A \implies B) \qquad \mbox{ if and only if for all $\beta \ge \alpha$}\\
           &&\mbox{\hskip 1in
                                                  $R_\beta \models_\sigma A$ implies $\exists \gamma \ge \beta
                                                  R_\gamma \models_\sigma B$} \\
R_\alpha &\models_\sigma& (\exists x\, A) \quad \mbox{if and only if for some $\tau$ that agrees with $\sigma$ }\\
           && \mbox{ \hskip 0.5in except perhaps on $x$, $R_\alpha \models_\tau A$} \\
R_\alpha &\models_\sigma& (\forall x\, A) \qquad \mbox{if and only if for all $\beta \ge \alpha$, if $x \in R_\beta$ and} \\
&& \mbox{\hskip 0.5in $\tau$ is a valuation on $R_\beta$ that agrees with $\sigma$} \\
&& \mbox{\hskip 0.5in except perhaps on $x$, then
                                             $R_\beta \models_\tau A$}
\end{eqnarray*}

A Kripke model of a theory $T$ is one such that all the axioms of $T$ are satisfied at every node of the model, 
i.e.,  $R_\alpha \models_\sigma A$ for every axiom $A$ and every valuation.  The Kripke completeness theorem 
says that the formulas $A$ that are provable from $T$ with intuitionistic logic are exactly the formulas
satisfied in all Kripke models of $T$.   In particular, if one can construct a Kripke model of $T$ in 
which some formula $B$ is not satisfied, then $B$ is not a consequence of $T$ with intuitionistic logic.%
\footnote{See \cite{vanDalen} for a proof of the Kripke completeness theorem, and Exercise 4, p. 99 of \cite{beeson-book}
for the extension to the logic of partial terms.}
Such a model is called a Kripke countermodel to $B$.  We will apply this technique to settle
the question of the reversibility of the implications between the different forms of 
the parallel postulate.

\subsection{A Kripke model whose points are functions}
\label{section:independence1}
The following concepts will be used to develop Kripke models in which 
the ``points'' are functions.

\begin{Definition} \label{definition:semidefinite}
 A function $f$ from 
$\R$ to $\R$ is called {\bf positive definite} if $f(x) > 0$ for every $x$.  

$f$ is called
{\bf positive semidefinite} if $f(x) \ge 0$ for all $x$. 

$f$ is called {\bf strongly positive 
semidefinite} if it is positive semidefinite and is not zero on any open interval.
\end{Definition}

\begin{Definition} \label{definition:pusieux}
 A {\bf Pusieux series} in $t$ is a power series in a rational power of  $t$, convergent in a neighborhood of 0.
 
 A {\bf Pusieux series at $a$ in $t$} is a power series in a rational power of $(t-a)$, convergent in a neighborhood of $a$.
 
A {\bf generalized Pusieux series in $t$} is a Pusieux series in $t$ or a power series in a rational power of $t$,
  times $\vert t \vert$, convergent in a neighborhood of 0.
  
A {\bf generalized Pusieux series in $t$ at $a$} is a Pusieux series in $(t-a)$ or a power series in 
a rational power of $(t-a)$, times $\vert t -a \vert$, convergent in a neighborhood of $a$.

A {\bf Pusieux series at $\infty$} is a Pusieux series convergent in a neighborhood of $\infty$; similarly 
for a {\bf generalized Pusieux series at $\infty$}.
\end{Definition}

For example, $\sqrt{t^3} = \vert t \vert\sqrt t$ has a generalized Pusieux series at 0, but not a Pusieux series.

\begin{Theorem} \label{theorem:kripke1}
   In constructive Euclidean field theory, the following are not provable:
   
(i)  ``two-sides'', namely $x \neq 0 \implies x > 0 \lor x < 0$;
 
(ii)  $x \le 0 \lor 0 \le x$

(iii)  ``apartness'', namely $0 < x \lor x < 1$.
 \end{Theorem}

\noindent{\em Proof}.
 Let $\K$ be the field of ``constructible numbers'',
which is the least subfield of $\R$ closed under taking square roots. 
\ Let $\A_0$ be the ring of 
polynomial functions from $\R$ to $\R$ with coefficients in $\K$.
  $\A_0$ is not a field, since for example $1/t$ does not belong to $\A_0$.
For each nonnegative integer $n$, we define the ring $\A_{n+1}$ to be the least set of 
real-valued functions containing $\A_n$ together 
with all sums, differences, and products of members of $\A_n$, together with all  square roots of positive semidefinite 
members of $\A_n$, and reciprocals of all strongly
positive semidefinite  members of $\A_n$.  These square roots and reciprocals are defined on dense 
subsets of $\R$, as shown below. 
For example, the functions $\sqrt{1+t^2}$  and $1/(1+t^2)$ are in $\A_1$, and 
$$ \sqrt{ \sqrt{1+t^2} + \sqrt{1+t^4}} +\frac 1 {1+t^2}$$
is in $\A_2$.  Also $\sqrt{x^2} = \vert x \vert$ is in $\A_2$,  and $\vert x \vert -x$ is in $\A_3$; 
that function is zero on the positive real axis.  
Now define $\A$ to be the union of the $\A_n$.  Then $\A$ is a ring of functions from $\R$ to $\R$.
Another way of describing $\A$ is to say that it is the least ring of functions containing $\K[x]$
and closed under taking square roots of positive semidefinite functions and reciprocals of
strongly positive semidefinite functions.  (If a function in $\A$ is positive semidefinite,
but vanishes on some interval, its square root is in $\A$, but not its reciprocal.)
 
 We claim  that each $f$ in $\A_n$ 
(i)  is defined at all but a finite number of points in $\R$, and has a    Pusieux series at all but a finite number of  
 points $a$, and also has a Pusieux series at $\infty$, and a generalized Pusieux series at the remaining points of its domain,
and (ii) $f$  is zero on a finite number of closed intervals and a finite number of isolated points.
  
For example $\vert x \vert = \sqrt{x^2}$ does not have a Pusieux series at 0; but such points
(of functions in $\A_n$) 
occur only at the zeroes of functions belong to $\A_{n-1}$, as we will show in detail below.
 
We prove (i) and (ii) simultaneously by induction on $n$.
First we show that (i) at $n$ follows from (ii) at $n-1$.  Case 1, $f$ is a square root, $f = \sqrt{g}$ and 
$g(x) \ge 0$ for all $x$.
Then at points where $g(x) > 0$,  $f$ has a Pusieux series, and at the finitely many points (by (ii) for $n-1$)
where $g(x) = 0$,  $f$ has a generalized Pusieux series; and $f$ is everywhere defined and has a Pusieux series 
at $\infty$. Case 2,   $f$ is a reciprocal of a strongly positive semidefinite function $g$ in $\A_{n-1}$.
Then by (ii) for $n-1$, 
$g$ has only finitely many zeroes (since there are no intervals of zeroes because $g$ is strongly positive
definite), so the domain of $f$ omits only finitely many points, 
and at each of those points where $g$ is positive, $f$ is defined and has a Pusieux series, and at the finitely 
many points where $g$ is zero, it is undefined; so (i) for $n$ follows from (ii) for $n-1$
also when $f$ is a reciprocal.  When $f$ is a sum or difference, (i) for $n$ follows from (i) for 
$n-1$.     Now we show that (ii) at $n$ 
follows from (i) at $n$.   By (i) at $n$,  $f$ in $\A_{n}$   has a  Pusieux series at $a$
for all but finitely many points $a$; and at the points where $f$ has a Pusieux series, if $f(a) = 0$ then the zero 
is isolated, or the function is identically zero in an interval about $a$. The endpoints of these intervals
of zeroes will be among the finitely many points where $f$ does not have a Pusieux series, so there are 
finitely many of these intervals.   Hence the zeroes of $f$ consist of (at most) finitely many points where 
$f$ does not have a Pusieux series, plus finitely many intervals, plus finitely many points where $f$ does 
have a Pusieux series--altogether finitely many.  So (ii) for $n$ does follows from (i) for $n$, as claimed.
That completes the inductive proof of (i) and (ii).

We call a zero of $f$  ``half-isolated''  if it is the endpoint of one of the intervals on which $f$ is zero.
By (i) and (ii), 
  there is a countable set  of reals that includes all the isolated or half-isolated 
  zeroes and all the singularities of all the functions in $\A$.
Define $\Omega$ to be the complement of that set; thus for each $f$ in $\A$,  if  $f(x)$ is zero
for any $x$ in $\Omega$, then $f$ is identically zero on some interval about $x$.
   Note that, since the complement of $\Omega$ 
is countable, $\Omega$ is dense in $\R$.

We will exhibit a Kripke model ${\mathbb K}$ of Euclidean field theory.
  The partially ordered set $K$ (of nodes of ${\mathbb K}$) is $\{0\} \cup \Omega$,
ordered so that $0$ is the root node, and $0 < \alpha$ for each $\alpha \in \Omega$, and the 
different elements of $\Omega$ are incomparable.  The root node
of this model is the ring $\A$, with $P(x)$ interpreted to mean ``$x$ is strongly positive semidefinite.''
We interpret equality as ``equal almost everywhere'';  that is, $x=y$ holds in the root node
if and only if $x(\alpha) = y(\alpha)$ except for a finite set of numbers $\alpha$.  For example,
$x \cdot y = 1$, where $x$ is the   function $x(t) = \vert t \vert$, and $y(t) = 1/\vert t \vert$, and $1$ is the 
constant function with value 1, because except for $\alpha = 0$, we have $x(\alpha)y(\alpha) = 1$.

$\A$ is not even an
ordered ring,  because there are polynomials $x$ such that neither
$x$ nor $-x$ is strongly positive semidefinite.  But we can still use $\A$ as the root of our Kripke model.

   For $\alpha \in \Omega$,  the structure $\A_\alpha$  at  the  node $\alpha$ is 
the quotient field of $\A$, with $P(x)$ interpreted to mean that $x(\alpha) > 0$, and $x=y$ 
interpreted to mean that $x$ and $y$ are equal on some neighborhood of $\alpha$. In other words,
the elements of $\A(\alpha)$ are equivalence classes of quotients of members of $\A$,  where 
$f/g$ is equivalent to $u/v$ if $vf$ and $ug$ agree in some neighborhood of $\alpha$.  
Since $\Omega$ does not include any singularities or zeros of members of $\A$, $x(\alpha)$ is defined 
for each $x$ in $\A_\alpha$. 

We note that $\A_\alpha$ is isomorphic to
the least Euclidean subfield of $\R$ containing $\alpha$.  The isomorphism is given by $x \mapsto x(\alpha)$.
It is an isomorphism because its kernel is trivial, since $x(\alpha) = 0 $ only if $x$ is identically zero
on some neighborhood of $\alpha$.
 Moreover, if $P(x) $ holds in $\A$, then $x$ is strongly positive semidefinite and not identically zero; 
hence for $\alpha$ in $\Omega$, we have $x(\alpha) > 0$, since $x(\alpha) \neq 0$ for $\alpha$ in $\Omega$.
Hence $P(x)$ holds in $\A_\alpha$.    

 The ring axioms are satisfied in this Kripke structure, since all these structures 
 are rings.  The ordered field axioms are satisfied at the leaf nodes $\A_\alpha$,
since these structures are classical ordered fields.   We therefore only need to verify the reciprocal and order and square root
 axioms at the root node $\A$. 

Consider Axiom EF2, which says that sums and products of positive elements are positive.  This holds
 at $\A$ since the sum and 
product of strongly positive semidefinite functions are also strongly positive semidefinite.    

Consider Axiom EF3,which says that not both $x$ and $-x$ are positive.  Suppose both $x$ and $-x$ are 
strongly positive semidefinite members of $\A$.  
Then for each $\alpha \in \Omega$ we have $x(\alpha) > 0$ and $-x(\alpha) > 0$, 
since $x(\alpha)$ is not zero for $\alpha$ in $\Omega$.  But that is a contradiction;
 hence Axiom EF3 holds at $\A$. 

Consider Axiom EF4, which says that if both $x$ and $-x$ are not positive, then $x$ is zero.  
  Suppose both $x$ and $-x$ are satisfied at $\A$ to 
be not positive.  That means that for every node $\A_\alpha$, $x$ and $-x$ are not positive at $\A_\alpha$.
That is,
$x(\alpha)\le 0$ and $-x(\alpha) \le 0$.  Hence,
$x(\alpha) = 0$.  But by construction, $\Omega$ excludes the zeroes of $x$, 
contradiction.   Hence $\A$ satisfies
Axiom EF4.

Consider Axiom EF5, which says that if $-x$ is not positive, then $x$ has a square root.   If $x$ is
identically zero there is nothing to prove, so we may assume that $x$ is not identically zero.
If $\A$ satisfies that $-x$ is not positive, i.e., $\A$ satisfies $\neg\, P(-x)$,
 that means that $-x$ is not positive in any $\A_\alpha$;
that is, $-x(\alpha) \le 0$ for all $\alpha \in \Omega$.  Then $x(\alpha) \ge 0$.   Since this is true 
for every $\alpha \in \Omega$, and since $\Omega$ is dense in $\R$,  and $x$ is continuous, it follows that $x$ is 
positive semidefinite.   Hence $\sqrt{x}$
belongs to $\A$, by construction of $\A$.   Hence $\A$ satisfies Axiom EF5.  

We now consider Markov's principle EF6.  Suppose that $\neg\neg P(x)$ is satisfied at the root node $\A$.
Then for every $\alpha$ in $\Omega$,  $P(x)$ is satisfied at the leaf node $\A_\alpha$;  that means that $x(\alpha) > 0$
for each $\alpha$ in $\Omega$.  As shown in the verification of E5, this implies that $x$ is positive semidefinite; and 
it is not identically zero on a neighborhood of $\alpha$,  since $x(\alpha) > 0$.  
Hence $P(x)$ is satisfied at the root node $\A$.

Axioms EF2 to EF7 imply that if $x \neq 0$, then $P(x^2)$.  We prove this as follows:
  if $P(x)$ then $P(x^2)$ by EF2.
If $-P(x)$ and $x\neq 0$ then $\neg\neg\,P(-x)$ by EF4, so $P(-x)$ by EF6, so $P((-x)^2)$
by EF2, so $P(x^2)$ since $(-x)^2 = x^2$ by the commutative ring laws.  Assume $x \neq 0$.
We have $\neg \neg (P(x) \lor P(-x))$, and either disjunct implies $P(x^2)$.  Hence 
$\neg\neg\, P(x^2)$.  Hence by EF6, $P(x^2)$.   

Now we consider Axiom EF1, which says nonzero elements have reciprocals.  It suffices to 
show that positive elements have reciprocals, since if $x \neq 0$, then $P(x^2)$,
and $1/x = x/x^2$.  More explicitly, if $y$ is a multiplicative inverse of $x^2$ then $xy$
is a multiplicative inverse of $x$.

Suppose $P(x)$ is satisfied at the 
root node; then $x$ is strongly positive semidefinite and not identically zero, and $x$ belongs to some $\A_n$; so
$1/x$ belongs to $\A_{n+1}$.  Then $x \cdot (1/x) = 1$ is satisfied, since 
$x(\alpha) \cdot(1/x(\alpha)) = 1$ except when $x(\alpha) = 0$, and the set of zeroes of $x$
is a finite set, because $x$ is strongly positive semidefinite.  Hence EF1 is satisfied.

The commutative ring axioms are satisfied, since the structure at every node of ${\mathbb K}$
is a commutative ring.
Thus ${\mathbb K}$ is a Kripke model of Euclidean field theory.

Therefore any statements that are not satisfied in ${\mathbb K}$ are not provable in 
Euclidean field theory.   For example, $P(x) \lor \neg\, P(x)$ is not satisfied, since 
when $x$ is the identity polynomial $i(t) = t$,  $P(i)$ fails at the root node because
$i$ is not positive semidefinite,  and $\neg\, P(i)$ also fails at the root node, since 
$i$ is positive at every leaf node $\alpha$ with $\alpha > 0$.

Ad (i), namely that ``two-sides'' is not satisfied:   Again let $i$ be the identity
function.  Then ${\mathbb K}$ satisfies $i \neq 0$, since $i=0$ is never satisfied at a 
leaf node $\alpha$ in $\Omega$, since $0$ does not belong to $\Omega$.  By definition 
$x < 0$ means $P(-x)$.  At the root node, neither $i$ nor $-i$ is positive definite,
so $0 < i \lor i < 1$ is not satisfied.  Hence $x \neq 0 \implies x > 0 \lor x < 0$
(which is ``two-sides'') is not satisfied.

Ad (ii), namely that $x \le 0 \lor 0 \le x$ is not satisfied:  $y \le x$ means 
$\neg P(y-x)$, so the principle under consideration is $\neg P(x) \lor \neg P(-x)$.
Again we take $x$ to be the identity function $i$.  Since at some nodes $\alpha$ we have
$P(i)$ (when $\alpha > 0$) and at other nodes we have $P(-i)$  (when $\alpha < 0$),
neither disjunct is satisfied at the root node.

Ad(iii), namely that $0 < x \lor x < 1$ is not satisfied:  in primitive syntax 
this principle is $P(x) \lor P(1-x)$.  Again taking $x$ to be the identity function, 
neither disjunct is satisfied at the root node, since neither $i$ nor $1-i$ is positive 
semidefinite.  

That completes the proof.
\smallskip

{\em Remark}.  In \cite{beeson2015b}, we use 
cut-elimination to give a completely different proof of 
this theorem.  In some respects that proof is better, as it shows a more general conclusion,
and it exposes the ``reason for'' these results: namely, the axioms are quantifier-free
and disjunction-free.   

The proof presented here is also interesting, because the 
``points'' of this model give us a better intuition into the possible interpretations 
of constructive geometry.  A point $(i,0)$ (where $i$ is the identity function) represents
a point about which all we know is that it is on the $x$-axis, but we don't know where.
It is not ``located''.   A point $(x,0)$,  where $\vert x \vert$ is a piecewise linear
function, equal to $-1$ for $t \le -1$ and to 1 for $t \ge 1$,  represents a point 
of which we know that it is between $-1$ and 1.   The model shows that there is a consistent
interpretation of this vague idea.

\subsection{A Kripke model with fewer reciprocals}

Recall that ``Playfair ring theory'' is the theory of ordered fields satisfying Axioms EF0, EF2-EF6, and 
EF7-EF8.  That is, without the axiom of reciprocals (EF7), but with 
 the axiom (EF7) that elements without reciprocals are zero, and the axiom (EF8) that says 
 that bounded quotients exist. 

\begin{Theorem}  \label{theorem:kripke2}
 Playfair ring theory  does not imply that
reciprocals of positive elements exist.
\end{Theorem}

\noindent{\em Proof}.  
To say that $\neg \exists y\,(y\cdot x = 1)$ holds at a node $\A$ of a Kripke model is to say that 
no node above $\A$ contains an inverse of $x$.   If one of the leaf nodes above $\A$ is a (classical) field,
then $x$ must be zero in that field, and hence in $\A$ also.  Hence the axiom that elements without reciprocals
are zero will hold in any Kripke model, all of whose leaf nodes are fields.  What we need, then, is a
Kripke model in which all the leaf nodes are ordered fields, and the root node has a positive element without
a reciprocal.  

Recall that $x$ is ``strongly positive semidefinite'' if it is positive semidefinite, and it is not zero on 
any open interval.  
We construct a model similar to the one in the preceding proof, but differing as follows: 
When constructing $\A$,
we throw in square roots of positive semidefinite functions as before, but instead of throwing in reciprocals
of all positive semidefinite functions, we throw in (only) bounded quotients.
  
   More precisely,  if $u$ and $v$ have Pusieux series at a point $\alpha$ where $v$ is zero, and 
$u/v$ is bounded in a neighborhood of $\alpha$, then $u/v$ can be 
extended (using the Pusieux series) to be defined at $\alpha$ as well.  When we write
$u/v$ in this context, we mean the extended function.
 Let $\A_0$ be the ring of 
polynomial functions from $\R$ to $\R$ with coefficients in $\K$.
  For each nonnegative integer $n$, we define the ring $\A_{n+1}$ to be the least set of 
real-valued functions containing $\A_n$ together 
with all sums and differences of members of $\A_n$, together with all  square roots of positive semidefinite 
members of $\A_n$, together with $u/v$ for each three functions $u,v,y$ in $\A_n$ such that 
$v$ has at most finitely many zeroes, and 
$$\vert u \vert \le \vert v \vert y.$$
wherever both sides are defined.
 
The domains of these functions are dense 
subsets of $\R$, as we will show below. 
Now define $\A$ to be the union of the $\A_n$.  Then $\A$ is a ring of functions from $\R$ to $\R$.
Another way of describing $\A$ is to say that it is the least ring of functions containing $\K[x]$
and closed under taking square roots of positive semidefinite functions,  and closed under the rule
that if $\vert u/v \vert \le y$ and $u$, $v$, and $y$ are in $\A$ and $v$ has only finitely 
many zeroes, then $u/v$ (extended to the zeroes of $v$) is in $\A$.

  As in the previous model, all the members of $\A$ have Pusieux series except at finitely 
many points, and their zero sets have the same structure as before, so we can again find a 
countable set containing all the singularities and all the isolated and half-isolated zeroes.
 Let $\Omega$ be the complement of this countable set, and define a
Kripke model as before,  with index set $\{0,\Omega\}$ and $R_0 = \A$ and for $\sigma \in \Omega$,
$R_\sigma = A_\sigma$, that is, $\A$ ordered by evaluation at $\sigma$. 
We claim 
 that all members of $\A$ are continuous on the real line.  To prove this, we proceed by 
 induction on $n$ to prove that the members of $\A_n$ are continuous
  on the whole real line.   That is true for $n=0$, since the elements of $\A_0$ are polynomials.
 The square root of a positive semidefinite function continuous on $\R$ is also continuous on $\R$,
 as are sums and differences of such functions.  Suppose $\vert u/v \vert \le y$ and 
 $u$, $v$, and $y$ are in $\A_{n}$.   Then $u/v$ will be in $\A_{n+1}$.  By induction hypothesis, 
 $u$, $v$, and $y$ are defined 
 everywhere.  Then as explained above, $u/v$ is extended by its Pusieux series to be defined
 at the zeroes of $v$, so it too is defined everywhere. 
 That completes the proof that all members of $\A$ 
 are continuous on $\R$. 
 
 We define $P(x)$ to hold at the root node if and only if $x$ is strongly positive semidefinite.  The 
 verification of axioms E2 through E6 is the same as in the previous proof.
 
 EF7 holds, because if $\forall y \neg\, (xy = 1)$ holds at the root node, then at every 
 leaf node $\alpha$ and for every $y$ in $\A$, we have $x(\alpha)y(\alpha) \neq 1$.  Hence
 $x(\alpha) = 0$ for every $\alpha$ in $\Omega$.  But $x$ is continuous and $\Omega$ is dense;
 hence $x$ is identically zero; hence the root node satisfies $x=0$. 
  
 Now we verify Axiom EF8.   It is automatic at the leaf nodes since there the structures are
 Euclidean fields.  We check the root node.  Suppose $\vert u \vert \le \vert v \vert y $.
 Then (the Pusieux extension of) $z$ of $u/v$ belongs to $\A$,   and satisfies
 $u = z v$.  Hence EF8 holds at the root node.

It only remains to show that EF1 is not satisfied at the root node.  If EF1 were satisfied at $\A$,
then every nonzero element of $\A$ would have its reciprocal in $\A$.  To refute this, take
$i$ to be the identity function.  Suppose $iy=1$ holds at the root node.  Then it holds
at every leaf node $\alpha$, so $y(\alpha) = 1/\alpha$ for each $\alpha$ in $\Omega$. 
But $y$ is continuous and $\Omega$ is dense, so $y$ is everywhere defined and $y(t) = 1/t$ for every $t$,
which is impossible when $t=0$.
 That completes the proof.
\medskip

 \subsection{Independence of Markov's principle }

\begin{Theorem} Markov's principle is independent of the other axioms of Euclidean field theory.
\end{Theorem}

{\em Proof}.  We construct a model similar to the ones used above, in which Markov's principle EF6
fails, but the other axioms hold.
Let $\A$ be the same class of functions used in the first model ${\mathbb K}$ 
 above, but this time 
interpret $P(x)$ at the root $\A$ as ``$x$ is positive definite'',  instead of as ``$x$ is positive semidefinite and
not identically zero.''   Then the function $x(t) = t^2$ is not satisfied to be positive at the root node, 
but since $0$ is not in $\Omega$, we have $x(\alpha) > 0$ for all $\alpha$ in $\Omega$, so $P(x)$ is satisfied
at every leaf node.  Hence $\neg \neg P(x)$ is satisfied at the root node.  Hence $\neg \neg P(x) \implies P(x)$
is not satisfied at the root.  Hence Markov's principle EF6 is not satisfied in this model.

The other axioms of Euclidean field theory are verified as in the proof of Theorem~\ref{theorem:kripke1}.
That completes the proof.
\medskip

\section{Independence of Euclid 5 }

In this section we give two different proofs that Playfair's axiom does not imply Euclid 5.  
Our first proof can be given without committing to 
a specific axiomatization of constructive geometry!  All we need to assume is that the 
axioms of the theory (other than the parallel axioms) are soundly interpreted by ordinary
analytic geometry in the theory of Playfair rings.  The second proof, on the other hand, 
does not depend on the theory of Playfair rings at all, but does depend on the existence of 
a quantifier-free axiomatization of constructive neutral geometry.  We present such 
an axiomatization in \cite{beeson2015b}.   Given that axiomatization, the second proof is 
perhaps simpler, but the first proof has the advantage of not depending on the details of the 
axiomatization.   It is worth giving two proofs, as each of the models we give illustrates
different aspects of constructive geometry.  The second proof, for example, shows
that Playfair does not imply Euclid 5, even with the aid of decidable equality and order,
while the first shows Playfair does not imply Euclid 5, even with the aid of the 
negations of decidable equality and order.

\subsection{Playfair does not imply Euclid 5}
 In this section, $\bar \phi$ denotes
the interpretation of the geometric formula $\phi$  in (the language of) Euclidean field theory,
as discussed above.  ``Neutral geometry'' means geometry without any parallel postulate.

\begin{Theorem}[Playfair does not imply Euclid 5] 
Suppose given a theory $T$ of neutral constructive geometry such that 
the analytic-geometry interpretation $\phi \mapsto \bar \phi$ is sound from $T$ to 
the theory of Playfair rings.  Then in $T$, 
  the strong parallel postulate implies Playfair's postulate, but not conversely. 
\end{Theorem}

\begin{Corollary} Playfair's postulate does not imply Euclid 5,  or any of the 
triangle circumscription principles.
\end{Corollary}
\medskip
  
\noindent{\em Proof}.  The corollary follows from the theorem, since we have proved
in Corollary \ref{lemma:Euclid5suffices} that Euclid 5 is equivalent to the strong 
parallel postulate, and in Theorem~\ref{theorem:trianglecircumscription} that the 
strong parallel postulate is equivalent to the triangle circumscription principles.

In Theorem~\ref{theorem:Euclid5impliesPlayfair}, we proved that Euclid 5
implies Playfair's axiom.  It therefore suffices to show that 
 Playfair's axiom does not imply the strong parallel postulate.
 Let $\phi$ be Playfair's axiom, and $\psi$ be the strong parallel postulate.
 Suppose $T$ prove $\phi \implies \psi$.  According to Theorem~\ref{theorem:soundnessofbar}, part (i),
  the theory of Playfair
fields proves $\bar \phi \implies \bar \psi$, and also 
$\bar \phi$.    Hence 
 $\bar \psi$ is provable in Playfair ring theory.  But, according to Theorem~\ref{theorem:soundnessofbar},
 part (ii), 
$\bar \psi$ is equivalent to EF1, the axiom that nonzero elements have reciprocals.
Hence the theory of Playfair rings proves EF1.  But that contradicts Theorem~\ref{theorem:kripke2}.
That completes the proof of the theorem.

\subsection{Another proof that Playfair does not imply Euclid 5}
In this section, we give a different Kripke model, providing a second proof of this independence
result.  Although technically we have defined Kripke models in this paper 
only for ordered ring theory, the 
definition of Kripke models for geometry is similar; only the atomic sentences are different. 
The model in question we call ${\mathcal D}$, after Max Dehn.  It has only two 
nodes.  Let $\F$ be a non-Archimedean Euclidean field.  A member $x$ of $\F$ is 
{\em infinitesimal} if $1/x$ is greater than every integer.  A member $x$ of $\F$ 
is {\em finitely bounded}, or just {\em finite}, if $\vert x \vert < n$ for some 
integer $n$.    Let ${\mathbb D}$
be the set of finitely bounded elements of $\F$.   Let the root node of ${\mathcal D}$ be ${\mathbb D}^2$, and 
the leaf node be $\F^2$.  

\begin{Theorem} The Kripke model ${\mathcal D}$ satisfies Playfair's axiom
and neutral constructive geometry, including line-circle and circle-circle continuity,
but does not satisfy Euclid 5. Moreover, it also satisfies the trichotomy law
$$x < y \lor x = y \lor x=y,$$ when expressed in geometric language. 
\end{Theorem}

\noindent{\em Proof}.  We start by pointing out that $\F^2$ is a model 
of classical Euclidean geometry.  Next we observe that ${\mathbb D^2}$ is  a model 
of classical neutral geometry.  (See Example 18.4.3 and Exercise
18.6 of \cite{hartshorne}.)
Specifically,  the points asserted
to exist by Pasch's axiom and by line-circle and circle-circle continuity have
finite coordinates, if the points given in the axioms have finite coordinates.
Euclid 5, however, fails in ${\mathbb D^2}$, since if $t$ is infinitesimal,
there is a line passing through $(0,1)$ with slope $t$,  whose intersection 
point with the $x$-axis is at $(1/t,0)$, which does not belong to ${\mathbb D}$.

Next we observe that it is possible to axiomatize neutral constructive geometry 
so that all the axioms are quantifier-free. For an 
example of such an axiomatization, see \cite{beeson2015b}.  Then all the formulas
in the axioms (with parameters from ${\mathbb D^2}$) hold at node ${\mathbb D}$ if and only
if they hold at node $\F^2$,  since that is true of equality, betweenness, and equidistance
statements.  Hence the Kripke model ${\mathcal D}$ satisfies constructive neutral 
geometry.   Since Euclid 5 fails classically in ${\mathbb D^2}$, it also fails in 
the Kripke model ${\mathcal D}$.

Lines in this model are specified by two points.  If $a \neq b$ then the points on 
$\Line(a,b)$ are the points satisfying $(x-a) \times (b-a) = 0$.  Thus if $a$ 
and $b$ are finite, there are points in $\F^2$ that are satisfied to be on $\Line(a,b)$
that are not in ${\mathbb D^2}$.   The interpretation of $\Line(a,b)$ is the extension 
of that line in $\F^2$.  

We next show that Playfair's axiom holds in ${\mathcal D}$.  Let line $L$ 
be given in ${\mathbb D}^2$, and let point $p$ in $\mathbb D^2$ not lie on $L$.
 Suppose 
$\mathbb D^2$ satisfies that line $M$ is parallel to $L$.  That is, 
$\neg \exists x\, (on(x,L) \land on(x,M)$.  Then $\F^2$, as a node above ${\mathbb D^2}$ in ${\mathcal D}$,
does not satisfy  $\exists x\, (on(x,L) \land on(x,M)$.   That is, $M$ and $L$ 
do not meet, even in $\F^2$.   Since $\F^2$ satisfies Euclid 5,  there is only one 
such line in $\F^2$.  Hence there is only one parallel to $L$ through $p$ in ${\mathcal D}$,
though there are many parallels to $L$ through $p$ in ${\mathbb D^2}$.
Hence Playfair's axiom holds in ${\mathcal D}$.

We now show that equality, equidistance, and betweenness are stable  in ${\mathcal D}$.  Let $A$ be an 
atomic statement (equality, betweenness, or equidistance) with parameters from ${\mathbb D}$
and suppose that ${\mathcal D}$ satisfies $\neg\neg\, A$.  Then $\F^2$ satisfies $A$.
But then ${\mathbb D^2}$ satisfies $A$, since ${\mathbb D^2}$ is a submodel of $\F^2$.

Finally, we observe that trichotomy holds in ${\mathcal D}$ since  betweenness relations
hold in ${\mathbb D^2}$ if and only if they hold in $\F^2$.
That completes the proof of the theorem.
\smallskip

{\em Remarks}. Similarly, the Kripke model with nodes $\F$ and ${\mathbb D}$ satisfies 
the axioms of Playfair rings, but not the Euclidean field axioms.  So we could have used
that model, together with the interpretation of geometry in Playfair rings, to prove
the theorem, instead of relying on rings of functions.
The two models cast different lights on constructive geometry: in the Dehn model,
Euclid 5 is unprovable because there are different lines whose intersection point cannot 
be bounded in advance, because it's ``too far away''.  On the other hand, in the 
model whose ``points'' are functions, the ``point'' given by the identity function 
represents a point of which we do not know anything at all about its position:  in the 
future, it might turn out to be {\em anywhere} on the line $L$, so if we use it as the 
slope of a line through $(0,1)$,  of course we can't compute the intersection point of 
that line with the $x$-axis.  In the Dehn model, there is no such ``uncertainty'' about 
the location of a point; points are located just as they are classically.  What
is uncertain is the extent of lines--just how far away is ``infinity''?  
The two models illustrate two different reasons why Euclid 5 is not constructively valid.

\section{Conclusion}
We have made the following main points in this paper:

\begin{itemize}
\item 
It is possible, indeed elegant, to do Euclidean geometry with intuitionistic logic, i.e., 
without arguments by cases.   One is permitted to prove statements about betweenness and congruence by contradiction,
using the stability axioms, as long as only atomic statements are proved that way.  It is not necessary to introduce
new,  purely intuitionistic concepts, such as apartness.
\smallskip

\item  In constructive geometry with Euclid 5,
all of Euclid's results are still provable.   Moreover, the geometric 
definitions of addition, multiplication, and square root can be given constructively, i.e., 
without case 
distinctions, so the models of constructive geometry are all planes over Euclidean fields.  This is the 
``representation theorem''.
\smallskip

\item Euclid 5 turns out to imply the strong parallel postulate.  This  is proved  by 
verifying that Euclid 5 suffices for coordinatization and for geometric definitions of addition 
and multiplication.  

\item There is a representation theorem for models of Playfair's axiom:  they are all planes $\F^2$ 
over Playfair rings, which are like Euclidean fields except instead of reciprocals of nonzero elements,
we only have the existence of bounded quotients.

\item Playfair's axiom, which Hilbert and most modern authors use in place of Euclid's 
parallel postulate,  is constructively weaker: It does not imply Euclid 5 using the other 
axioms of constructive geometry.
This result is  not sensitive to the exact choice of language or of the axioms of neutral geometry; all one needs is that the 
axiomatization should be quantifier-free and disjunction-free, or alternatively, that the 
representation theorem holds.
\end{itemize}
\medskip

These results complement our other work on constructive geometry:
In \cite{beeson-kobe} and \cite{beeson2015b} we have given formal theories
in the styles of Hilbert and Tarski, respectively.  Since these axiomatizations are 
 quantifier-free and  disjunction-free theory, we obtain from cut-elimination
 and   G\"odel's double-negation interpretation that

\begin{itemize}
\item
 Points, lines, and circles that one can prove to exist constructively
are constructible with ruler and compass,  uniformly (and hence continuously) in parameters.

\item For quantifier-free theorems (with Skolem functions corresponding to existential axioms),
non-constructive proofs can be converted to constructive proofs.  Hence we can import 
may of the  theorems from \cite{schwabhauser} into constructive geometry.
\end{itemize}

%\bibliographystyle{siam}
%\bibliography{../geometry} 

\end{document}